\documentclass{amsproc}
\usepackage{amsmath, amsthm, amssymb, slashed, framed, tikz}

\usepackage{epstopdf}
\usepackage{url}
\newtheorem{theorem}{Theorem}[section]
\newtheorem{lemma}[theorem]{Lemma}

\theoremstyle{definition}
\newtheorem{definition}[theorem]{Definition}
\newtheorem{example}[theorem]{Example}

\newtheorem{conjecture}[theorem]{Conjecture}

\newtheorem{proposition}[theorem]{Proposition}
\newtheorem{corollary}[theorem]{Corollary}

\theoremstyle{remark}
\newtheorem{remark}[theorem]{Remark}

\newdimen\tableauside\tableauside=1.0ex
\newdimen\tableaurule\tableaurule=0.4pt
\newdimen\tableaustep
\def\phantomhrule#1{\hbox{\vbox to0pt{\hrule height\tableaurule width#1\vss}}}
\def\phantomvrule#1{\vbox{\hbox to0pt{\vrule width\tableaurule height#1\hss}}}
\def\sqr{\vbox{%
\phantomhrule\tableaustep
\hbox{\phantomvrule\tableaustep\kern\tableaustep\phantomvrule\tableaustep}%
\hbox{\vbox{\phantomhrule\tableauside}\kern-\tableaurule}}}
\def\squares#1{\hbox{\count0=#1\noindent\loop\sqr
\advance\count0 by-1 \ifnum\count0>0\repeat}}
\def\tableau#1{\vcenter{\offinterlineskip
\tableaustep=\tableauside\advance\tableaustep by-\tableaurule
\kern\normallineskip\hbox
    {\kern\normallineskip\vbox
      {\gettableau#1 0 }%
     \kern\normallineskip\kern\tableaurule}%
  \kern\normallineskip\kern\tableaurule}}
\def\gettableau#1 {\ifnum#1=0\let\next=\null\else
  \squares{#1}\let\next=\gettableau\fi\next}

\numberwithin{equation}{section}

\title{Quadruply-graded colored homology of knots}

\author{Eugene Gorsky}
\address{Department of Mathematics, Stony Brook University, Stony Brook, NY 11733, USA}
\email{egorsky@math.sunysb.edu}

\author{Sergei Gukov}
\address{California Institute of Technology, Pasadena, CA 91125, USA\\
\newline
Max-Planck-Institut f\"ur Mathematik, Vivatsgasse 7, D-53111 Bonn, Germany.}
\email{gukov@theory.caltech.edu}
\author{Marko Sto$\check{\text{s}}$i$\acute{\text{c}}$}

\address{Instituto de Sistemas e Robotica and CAMGSD, Instituto Superior Tecnico, Torre
Norte, Piso 7, Av. Rovisco Pais, 1049-001 Lisbon, Portugal\\
\newline
Mathematical Institute SANU, Knez Mihailova 36, 11000 Beograd, Serbia.}
\email{mstosic@math.ist.utl.pt}

\font\teneurm=eurm10 \font\seveneurm=eurm7 \font\fiveeurm=eurm5
\newfam\eurmfam
\textfont\eurmfam=\teneurm \scriptfont\eurmfam=\seveneurm
\scriptscriptfont\eurmfam=\fiveeurm

 \font\teneusm=eusm10 \font\seveneusm=eusm7 \font\fiveeusm=eusm5
\newfam\eusmfam
\textfont\eusmfam=\teneusm \scriptfont\eusmfam=\seveneusm
\scriptscriptfont\eusmfam=\fiveeusm

\font\tencmmib=cmmib10 \skewchar\tencmmib='177
\font\sevencmmib=cmmib7 \skewchar\sevencmmib='177
\font\fivecmmib=cmmib5 \skewchar\fivecmmib='177
\newfam\cmmibfam
\textfont\cmmibfam=\tencmmib \scriptfont\cmmibfam=\sevencmmib
\scriptscriptfont\cmmibfam=\fivecmmib

\def\example#1{\bgroup\narrower
\baselineskip\footskip\bigbreak
\hrule\medskip\nobreak\noindent {\bf Example}. {\it #1\/}\par\nobreak}
\def\endexample{\medskip\nobreak\hrule\bigbreak\egroup}

\newcommand{\be}{\begin{equation}}
\newcommand{\ee}{\end{equation}}
\newcommand{\bea}{\begin{eqnarray}}
\newcommand{\eea}{\end{eqnarray}}

\newcommand{\C}{\mathbb{C}}
\newcommand{\Z}{\mathbb{Z}}
\newcommand{\R}{\mathbb{R}}

\def\tilde{\widetilde}
\def\hat{\widehat}
\def\bar{\overline}

\def\CA{{\mathcal A}}
\def\CB{{\mathcal B}}

\def\CD{{\mathcal D}}
\def\CE{{\mathcal E}}
\def\CF{{\mathcal F}}

\def\CH{{\mathcal H}}
\def\CI{{\mathcal I}}

\def\CL{{\mathcal L}}
\def\CM{{\mathcal M}}
\def\CN{{\mathcal N}}
\def\CO{{\mathcal O}}
\def\CP{{\mathcal P}}

\def\CR{{\mathcal R}}

\def\CT{{\mathcal T}}
\def\CU{{\mathcal U}}
\def\CV{{\mathcal V}}

\def\t23{\tableau{2 2 2}}
\def\t32{\tableau{3 3}}
\def\d23{d^{\t23}}

\def\a{a}
\def\q{q}
\def\t{t}

\def\unknot{{\,\raisebox{-.08cm}{\includegraphics[width=.4cm]{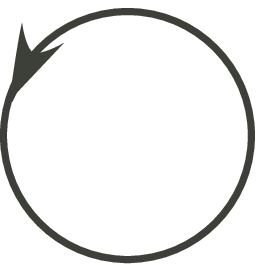}}\,}}

\def\cp{{\mathbb{C}}{\mathbf{P}}}

\def\leadsto{\rightsquigarrow}
\def\p{\partial}

\DeclareMathOperator{\Hilb}{Hilb}
\DeclareMathOperator{\Hom}{Hom}
\DeclareMathOperator{\Ext}{Ext}
\DeclareMathOperator{\Sym}{Sym}

\DeclareMathOperator{\Coef}{Coef}
\DeclareMathOperator{\Spec}{Spec}

\begin{document}


\begin{abstract}
We conjecture the existence of four independent gradings in the colored HOMFLY homology.
We describe these gradings explicitly for the rectangular colored homology of torus knots
and  make qualitative predictions of various interesting structures and symmetries in the colored homology of general knots.
We also give a simple representation-theoretic model for the HOMFLY homology of the unknot colored by any representation.
While some of these structures have a natural interpretation in the physical realization of knot homologies based on
counting supersymmetric configurations (BPS states, instantons, and vortices), others are completely new.
They suggest new geometric and physical realizations of colored HOMFLY homology as the Hochschild homology
of the category of branes in a Landau-Ginzburg B-model or, equivalently, in the mirror A-model.
Supergroups and supermanifolds are surprisingly ubiquitous in all aspects of this work.
\end{abstract}

\maketitle

{\tt {CALT-68-2903}}

\tableofcontents


\section{Introduction}

The categorification of quantum knot invariants started with Khovanov's seminal paper \cite{Khfirst}, where he defined a doubly-graded homology theory, whose homotopy type is an invariant of a knot, and such that its graded Euler characteristic is equal to the Jones polynomial. This opened a fast-growing field in low-dimensional topology and started work on categorification of other quantum knot invariants. Unlike the quantum polynomial invariants of a knot $K$, denoted by $P^{g,\lambda}(K)(q)$, which are explicitly defined for any representation $\lambda$ of Lie algebra $g$, their categorifications are known only in very few cases. In particular, in the case of the Lie algebras $sl(N)$ --- which is the case on which we focus in this paper --- the explicit categorification is known in the case of the fundamental representations $\lambda=\tableau1$ \cite{KRa}, and in a very few other cases. In addition, apart from the case of the Jones polynomial (fundamental representation of $sl(2)$), the definition of the homology theory is very complicated, which makes the computation of such invariants extremely hard.

For a fixed Young diagram $\lambda$ and a knot $K$, the $sl(N)$ quantum polynomials $P^{sl(N),\lambda}(K)(q)$ can be organized in a single two-variable $\lambda$-colored HOMFLY polynomial $P^{\lambda}(K)(a,q)$, so that the $sl(N)$-quantum polynomials are obtained as $a=q^N$ specializations, i.e.
\be
P^{sl(N),\lambda}(K)(q)=P^{\lambda}(K)(a=q^N,q).
\label{PPN}
\ee
A triply-graded homology categorifying uncolored ($\lambda=\tableau1$) HOMFLY polynomial was defined in \cite{KRb}. Throughout the paper, we are mainly focusing on the reduced polynomials, i.e. normalized such that the value of the unknot is equal to 1, and unless otherwise stated, the (colored) HOMFLY polynomial means the reduced (colored) HOMFLY polynomial.

A remarkable property of the quantum knot invariants and their categorifications is that they are related to many different areas of mathematics and physics, all of them bringing new viewpoints to this topic. The relationship between quantum field theory, in particular the Chern-Simons theory, and the quantum knot invariants was discovered in celebrated paper by Witten \cite{Witt}, and the physics insights have been extremely fruitful ever since.

In the case of the categorifications, the realization of knot homologies as the space of certain BPS states gave rise to various prediction on the structure of the (colored) HOMFLY homologies, see e.g. \cite{GSaberi} for a friendly introduction and a review.
Thus, in \cite{DGR} N. Dunfield, the second author and J. Rasmussen predicted the existence of a triply-graded knot homology theory $\CH^{\tableau1}(K)$, whose Euler characteristic is given by the HOMFLY polynomial, and which comes equipped with the collection of anti-commuting differentials $\{d_N\}_{N\in\Z}$ lifting the $sl(N)$ specialization property of the HOMFLY polynomial. More precisely, the differentials $d_N$, for $N>0$, are such that the homology of $\CH^{\tableau1}(K)$ with respect to $d_N$ is isomorphic to the homology $\CH^{sl(N),\tableau1}(K)$ that categorifies $P^{sl(N),\tableau1}(K)(q)$. Additionally, $\CH^{\tableau1}(K)$ has an involution $\phi$ extending the $q\leftrightarrow q^{-1}$ symmetry of the HOMFLY polynomial:
\be
P^{\tableau1}(K)(a,q)=P^{\tableau1}(K)(a,q^{-1}),
\ee
It exchanges the positive and negative differentials: $\phi d_N = d_{-N} \phi$.

Such a rigid structure enabled the computation of triply graded homology of various knots. In addition such a triply-graded homology theory was constructed mathematically in \cite{KRb}, and it was shown in \cite{Ras} that it also carries  the differentials $d_N$ with the wanted properties.

In \cite{GS} the second and the third authors extended this picture and described various structures in the {\em colored} HOMFLY homology in the case of the symmetric ($S^r$) and anti-symmetric ($\Lambda^r$) representations. For such representation $\lambda$ and a knot $K$ one can associate an $(a,q,t)$-graded vector space $\CH^{\lambda}(K)$ such that its Euler characteristic with respect to the $t$-grading equals to the $\lambda$-colored HOMFLY polynomial $P^{\lambda}(K)(a,q)$ of $K$. Like in the uncolored case, these homology theories come with the collection of the differentials corresponding to the $sl(N)$ specializations. The main new feature is the existence of another collection of the so-called {\em colored} differentials which give ``dynamics'' in the sequence of the homology theories $\CH^{S^r}(K)$, for various $r$. For every pair of nonnegative integers $r$ and $k$ with $r>k$ there are two different differentials  $d^{+}_{S^{r}\to S^{k}}$ and $d^{-}_{S^{r}\to S^{k}}$ on $\CH^{S^r}(K)$ such that the homology of $\CH^{S^r}(K)$ with respect to any of the two differentials $d^{\pm}_{S^{r}\to S^{k}}$ is isomorphic to $\CH^{S^k}(K)$.

The involution $\phi$  becomes the so-called {\em mirror symmetry} in colored homology, that relates the homologies $\CH^{\lambda}(K)$ and $\CH^{\lambda^t}(K)$: there exists an isomorphism between $\CH^{S^r}(K)$ and $\CH^{\Lambda^r}(K)$  preserving the $a$-grading and reversing the $q$-grading. It lifts the following relation between the colored HOMFLY polynomials
\be
P^{\lambda}(K)(a,q)=P^{\lambda^t}(K)(a,q^{-1}).
\label{PPmir}
\ee

Surprisingly, in \cite{GS}
it was realized that the {\em exponential growth property} holds for certain classes of knots, i.e. the size of the $S^r$-colored homology grows exponentially in $r$:
\be
\dim \CH^{S^r}(K)=\left(\dim \CH^{\tableau1}(K) \right)^{r}.
\ee
Although this property does not hold for all knots, it holds for various large classes of knots, like 2-bridge knots and torus knots. Moreover, the refined exponential growth is valid  for the two-variable $q=1$ specializations of the Poincar\'e polynomials of the corresponding homology theories. 

Such a large predicted structure made possible the explicit computation of colored HOMFLY homologies for various knots and colors, that in turn have shown the consistency of all
the conjectured properties.\\


In present article we give more quantitative description of the results of \cite{GS}.
 We extend the description of the colored differentials and the exponential growth conjecture to rectangular diagrams $\lambda$. This allows us to give a more unified treatment of the $\Lambda^r$ and $S^r$ colored homology.

The major novelty is endowing the space $\CH^{\lambda}(K)$ with the fourth grading. In other words, we conjecture the colored HOMFLY homology $\CH^{\lambda}(K)$ of a knot $K$ to carry four independent gradings, which we will denote by $a,q,t_{r},t_{c}$: apart from the $a$  and $q$ gradings from the polynomial invariants, we introduce two homological (or refined) $t$-gradings, denoted $t_r$ and  $t_c$.  Both $t$-gradings on the HOMFLY homology give a categorification of $P^{\lambda}(K)$. In other words, if $\CP^{\lambda}_r(K)$ and $\CP^{\lambda}_c(K)$ are three-variable Poincar\'e polynomials of $\CH^{\lambda}(K)$ with respect to $t_r$ and $t_c$ gradings, respectively, then
\[
\CP^{\lambda}_r(K) (a,q,t_{r}=-1) \; = \; \CP^{\lambda}_c(K) (a,q,t_{c}=-1) \; = \; P^{\lambda}(K)(a,q),
\]
The appearance of the fourth grading is yet mysterious for us from geometric point of view, but it seems to be inevitable. Let us list some of the evidences.

First of all,
in such a way we managed to reconcile the two different conventions
for the homological grading in the case of symmetric representations:
 $t_r$ is the $t$-grading assigned to a generator of $\CH^{\lambda}(K)$
in the grading conventions of \cite{GS}, whereas one can interpret $t_{c}$ as
the $t$-grading assigned to a generator of $\CH^{\lambda}(K)$ in grading conventions of \cite{AShakirov,DMMSS,DGR,GWalcher}.
The ``mirror symmetry'' exchanges the two grading conventions. Therefore, one practical application of
``mirror symmetry'' is that it allows to convert from one set of grading conventions used in the literature
to the other one by looking at the HOMFLY homology of the knot colored by transposed Young tableaux $\lambda^t$.

One advantage of the quadruply-graded theory is that it makes
all of the structural properties and isomorphisms completely explicit. This cannot be achieved within a triply-graded theory: for example in the case of the symmetric representations (see \cite{GS}) only a few of the isomorphisms and symmetries could be made explicit.
Now, quadruply-graded theory $\CH^{\lambda}(K)$ has remarkable properties: it gives an explicit re-gradings of all colored isomorphisms, enables
explicit expression in all gradings for the mirror symmetry, and makes the exponential growth property manifest as a fully refined exponential growth property of three-variable polynomials. In addition, miraculously, we discovered a new symmetry on $\CH^{\lambda}(K)$ that we call {\em self-symmetry}. This symmetry requires all four gradings, and unlike the mirror symmetry it does not lift any known relation
of the colored HOMFLY polynomials.

What is even more surprising, is that all these symmetries, isomorphism re-gradings and exponential growth property have
rather simple and pleasant expressions. Moreover, all of them
become particularly elegant when expressed in terms of an auxiliary
grading, called $Q$-grading, defined in the following simple way when $\lambda$ is a rectangular Young diagram with $R$ rows:
\be\label{Qgrad}
Q=\frac{q+t_r-t_c}{R}.
\ee
This new grading, which can be considered as a certain ``corrected'' $q$-grading,  cannot be seen on the decategorified, polynomial level, and as the formula indicates, both $t$-gradings are
needed for its definition. Therefore we need the quadruply-graded theory in order to have all properties manifest,
while only a small fraction of this rich structure can be seen at the polynomial, ``decategorified'' level.

Finally, the symmetries take the most elegant form when the homology is written in $(a,Q,t_r,t_c)$-gradings. Because of a linear relation between the gradings (\ref{Qgrad}), this is just a simple re-grading of $\CH^{\lambda}(K)$. However, due to its importance, we give it a special name, $\tilde{\CH}^{\lambda}(K)$, and refer to it as the {\it tilde-version} of colored HOMFLY homology. Thus, it is given by
\be
\tilde{\CH}^{\lambda}_{i,j,k,l}(K) \; := \; {\CH}^{\lambda}_{i,Rj-k+l,k,l}(K).
\ee

We note that only in the uncolored case the $t_{r}$ and $t_{c}$ gradings coincide (as well as $Q$ and $q$ gradings), and the resulting homology is triply graded in agreement with \cite{DGR}.\\

The second major novelty is that we extend the definition of the differentials $d_N$ to a two-dimensional family of differentials $d_{n|m}$ labeled by Lie superalgebras ${sl}(n|m)$. Then, the value of $N$ from the familiar one-parameter family of the differentials $\{ d_N \}$ corresponds to the difference $n-m$ which is precisely the super-rank of ${sl}(n|m)$. This gives a much more natural treatment of the differentials $d_N$, for non-positive $N$, together with their gradings. We show that these new differentials are nontrivial even in the uncolored homology (for sufficiently large knots) and bring an interesting structure both to colored and uncolored knot homology. We also describe explicitly the interaction between the differentials  $d_{n|m}$ and the colored differentials for the rectangle-colored homology.

While it still deserves a much deeper understanding, the appearance of $sl(n|m)$ is natural for a number reasons. Primarily, the representation theory of $sl(n|m)$ explains the behavior of the colored differentials. For each rectangular Young diagram $\lambda$ we define colored differentials removing any number of columns or rows from $\lambda$. These differentials naturally appear from the identification of the representations of the superalgebra $sl(n|m)$ labeled by two different rectangles $\lambda$ and $\mu$. It turns out that these are such that $\mu$ is obtained from $\lambda$ by erasing some of its rows or some of its columns.  Then, the corresponding colored differential,
$d_{\lambda \to \mu}$, closes the following commutative diagram:

\begin{center}
\begin{tabular}[c]{ccc}
$\CH^{\lambda}(K)$ & $\xrightarrow{d_{\lambda\to \mu}}$ & $\CH^{\mu}(K)$\cr
 $\vline$ & & $\vline$ \cr
${d_{n|m}}$ & & ${d_{n|m}}$\cr
 $\downarrow$ & & $\downarrow$ \cr
$\CH^{sl(n|m),\lambda}(K)$ & $\xrightarrow{\qquad \cong\qquad}$ & $\CH^{sl(n|m),\mu}(K)$
\end{tabular}
\end{center}

Furthermore, the differentials $d_{n|m}$ on $\CH^{\lambda}(K)$ exist for all pairs of nonnegative integer $(n,m)$ such that the superalgebra $sl(n|m)$ has a representation corresponding to the Young diagram $\lambda$. In such a way, all differentials have natural interpretation, unlike the single-parameter family of differentials $\{d_N\}_{N\in\Z}$ for which only the ones with $N>0$ have clear interpretation as $sl(N)$ specialization differentials.
In addition, supergroups and the corresponding Lie superalgebras are familiar in the study of brane/anti-brane systems \cite{VafaUNM,WittenK}:
much like a collection of $N$ coincident D-branes carries a gauge bundle with the structure group $U(N)$, a similar system of $n$ branes and $m$ anti-branes carries $U(n|m)$ gauge symmetry. Therefore, since many physical realizations of knot homologies that we encounter
in sections \ref{sec:geomphys}, \ref{sec:supermn}, \ref{sec:mirror}, and \ref{sec:inst-surf} are based on branes, it is not inconceivable that appearance of supergroups and superalgebras is rooted there (although we will not try to pursue this interpretation in the present paper).\\

Another direction of our study is the construction of the models for colored homology
for certain knots satisfying all above properties. We give a combinatorial model for the unreduced homology of the unknot and study various differentials on it. The HOMFLY homology is conjectured to be a free supercommutative algebra,
\be
\CA \; = \; \CH^{\lambda} (\unknot)
\label{Aunknot}
\ee
and the differentials are described by the Koszul construction. We compare this approach with other mathematical and physical approaches (potentials and matrix factorizations, categorification of the Jones-Wenzl projectors, algebra of BPS states).
We give an algebraic model for the $S^r$ colored homology of torus knots and compare it with the approaches of \cite{GORS} and \cite{EGL}. We study  various differentials in this model as well and show how the differentials $d_{n|m}$ naturally appear from the rational Cherednik algebra. Finally, to demonstrate the full strength and robustness of our approach, we describe the colored homology for the figure eight knot, which is a non-torus knot.


\subsection{Structure of the paper}

In Section \ref{sec:bosonic} we describe a model for the unreduced colored homology of the unknot.
It is a a free supercommutative algebra $\CA$ with one even and one odd generator per each box of the diagram $\lambda$.
For the representation $S^n$, this algebra can be naturally identified with the algebra of $S_n$-invariant differential forms on $\C^n$. Using the physics of BPS states, we describe some of these even and odd generators geometrically. We also check that the Poincar\'e polynomial for this algebra agrees with the evaluation formulas for the Macdonald polynomial, that various authors \cite{AShakirov,Ch,DMMSS} recently assign to the unknot in the refined Chern-Simons theory.

Both $d_{m|n}$ and the colored differentials can be interpreted as the Koszul differentials on this algebra, sending odd generators to certain polynomials in even generators. We show that the approach of \cite{GWalcher}, where the homology of the unknot is interpreted as the Jacobi ring of a certain potential $W$, fits into this framework if the Koszul complex is associated with the partial derivatives of $W$:
$$d(\xi_i)=\frac{\partial W}{\partial x_i},\ d(x_i)=0.$$
Using this model, we also derive the equations for the colored differentials: if the representations of $sl(m|n)$ labelled by $\lambda$ and $\mu$ are isomorphic, we show that in some  examples the isomorphism between the corresponding Jacobi rings follows from the equation
$$W_{m|n}^{\lambda}(x_1,\ldots,x_{|\lambda|})=W_{m|n}^{\mu}(x_1,\ldots,x_{|\mu|})+W_{\lambda\to \mu}(x_{|\mu|+1},\ldots,x_{|\lambda|}),$$
where $W_{\lambda\to \mu}$ is a non-degenerate quadratic function. We conjecture that the colored differential
$d_{\lambda\to \mu}$ is a Koszul differential associated with the partial derivatives of $W_{\lambda\to \mu}$.

In Section \ref{sec:gradings} we describe the general structures on the colored HOMFLY homology.
We give a detailed list of the properties of the quadruply-graded  colored HOMFLY homology $\CH^{\lambda}(K)$ of a general knot $K$, together with the $sl(n|m)$ and colored differentials. Furthermore we describe explicitly the degrees of all colored differentials and re-gradings in the corresponding colored isomorphisms. We end the section with the two particular features in the case of the symmetric representations: the  $\delta$-grading and the existence of the HFK-like differential on the $S^r$-colored homology. Also, we discuss similarities and differences between the case of rectangle-labeled  representations, which are the main topic of our present paper, and the case of non-rectangular representations.

Section \ref{sec:examples} contains the computations of the quadruply-graded homologies for various knots and representations. Due to a large number of predicted properties such computations present a highly non-trivial check and confirmation of the theory.

Section \ref{sec:algmodel} is devoted to torus knots.
In Section  \ref{sec:torus} we recall the statement about the stabilization of the HOMFLY polynomials of $(p,q)$ torus knots in the limit $q\to \infty$. We relate such stable $\lambda$-colored invariant of $(p,\infty)$ torus knot
to the invariant of the $p\lambda$-colored unknot. Then, we conjecture a similar relation between the corresponding colored homology theories of stable torus knots and that of the unknot. The $\lambda$-colored homology of the finite $(p,m)$ torus knot can be realized as a certain quotient of this stable homology. This allows us to use the algebraic description of the homology of the unknot, and in Section \ref{sec:alggrad} we give a precise description of the generators in the rectangular homology of a torus knot and their quadruple gradings.
We describe and check the structural properties of colored homology using this algebraic model.
We show the existence of a symmetry between $\lambda$- and $\lambda^{t}$-colored homology, exchanging $t_r$ and $t_c$ gradings and preserving the $Q$-grading. We check the refined exponential growth conjecture for the stable homology of torus knots. We also describe respectively the self-duality involution  and the action of differentials in the rectangular homology of torus knots.

Generalizing the results of \cite{GORS} and \cite{EGL}, we identify the
colored homology of torus knots with the certain representations of the rational Cherednik algebra. For $\lambda=S^r$, we make this description precise and describe the space $\CH^{S^r}(T(m,n))$ as a space of differential forms on a certain non-reduced scheme
$$\CM_{m,n,r}:=\{(u_{r+1},\ldots, u_{nr};v_{r+1},\ldots, v_{mr}\ |\ U(z)^m=V(z)^n\},$$
where the polynomials $U(z)$ and $V(z)$ are defined by the equations
$$U(z)=1+z^{r+1}u_{r+1}+\ldots+z^{nr}u_{nr},\ V(z)=1+z^{r+1}u_{r+1}+\ldots+z^{mr}u_{mr}.$$
This description is ideally suited for a realization of colored HOMFLY homology via Landau-Ginzburg B-model,
which is the subject of Section \ref{sec:mirror}.
Using mirror symmetry we also reformulate it in terms of the A-model and Lagrangian Floer homology.
We also provide some explicit formulas for the differentials acting on the HOMFLY homology
and describe $\CH^{S^r}(T(m,n))$ as a Jacobi ring of a certain potential on a supermanifold.

Finally, section \ref{sec:bottom} is devoted to the ``bottom row'' of the HOMFLY homology, i.e. $a=0$ limit.
We compare it to certain combinatorial models of Haglund et. al., generalizing the models for $q,t$-Catalan numbers
which played important role in \cite{GORS}, and to the geometric models for coupled instanton-vortex systems proposed in \cite{DGH}.

\subsection{Acknowledgements}

We are grateful to
M. Abouzaid, M. Aganagic, J.M. Baptista, I. Cherednik,
K. Costello, R. Elliot, P. Etingof, A. Gorsky, K. Hikami, M. Khovanov, B. Kim, A. N. Kirillov and A. A. Kirillov Jr.,
I. Losev, A. Morozov, H. Nakajima, A. Negut,  N.~Nekrasov, A. Oblomkov, A. Okounkov, J. Rasmussen, L. Rozansky,
S. Shakirov, V. Shende, C. Vafa, O. Viro, E. Witten, and C. Woodward
for the useful discussions.
E. G. would like to thank the California Institute of Technology and Kyoto Research Institute for Mathematical Sciences for the hospitality.
The research of E. G. is partially supported by the grants RFBR-10-01-678, NSh-8462.2010.1 and the Simons foundation.
S.G. would like to thank Instituto Superior Tecnico in Lisbon and the Simons Center for Geometry and Physics at Stony Brook for hospitality during the key stages of this work.
The work of S. G. is supported in part by DOE Grant DE-FG03-92-ER40701FG-02 and in part by NSF Grant PHY-0757647.
M. S. would like to thank the California Institute of Technology for hospitality during which part of this work was done. 
The work of M. S. was partially supported by the Portuguese Funda\c c\~ao para a Ci\^encia e a Tecnologia through ISR/IST plurianual funding and through the project number PTDC/MAT/101503/2008, New Geometry and Topology, and by the Ministry of Science of Serbia, project no. 174012.
Opinions and conclusions expressed here are those of the authors and do not necessarily reflect the views of funding agencies.

\section{Bosonic and Fermionic generators}
\label{sec:bosonic}

\subsection{The algebra of BPS states and knot cobordisms}
\label{sec:geomphys}

While most of this paper is probably aimed at a more mathematical reader, here we offer a geometric
interpretation of some of our observations in the framework more familiar to physicists.
Therefore, the purpose of this section is twofold. First, it will help physicists to
understand better some of our observations (which, we believe, may go beyond the scope of the present
paper and are worth studying further in the context of BPS state counting).
Second, the physical interpretation / motivation provides further evidence
for the structures of knot homologies discussed in this paper.
In Section \ref{sec:inst-surf}
we also comment on the connection with the equivariant instanton counting in the presence of a surface operator.

We begin with a lightning review of the physical interpretation of knot homologies, in which
$\CH_{\text{knot}}$ is identified with a space of supersymmetric particles (the so-called BPS states) \cite{GSV,GukovRTN,fiveknots}:
\be
\CH_{\text{knot}} \; = \; \CH_{\text{BPS}} \,.
\label{HHBPS}
\ee
In particular, it will allow us to see bosonic and fermionic generators of homological knot invariants
and will also help us in section \ref{sec:inst-surf} to understand the limit $\a \to 0$.
As we shall see, the latter allows one
to focus on the ``bottom row'' of the superpolynomial
which turns out to have a very elegant and simple geometric interpretation in the physical framework.

In one of the duality frames, the physical system involves a compactification on a very particular (non-compact) Calabi-Yau 3-fold $X$
in the presence of five-branes supported on a Lagrangian submanifold $L_{K} \subset X$:
\bea
\text{space-time} & : & \qquad \R \times X \times M_4 \label{theoryB} \\
\text{M5-branes} & : &  \qquad \R \times L_K \times D \nonumber
\eea
Specifically, $X$ is the total space of the $\CO (-1) \oplus \CO (-1)$ bundle over $\cp^1$,
while the Lagrangian submanifold $L_K$ is determined by a knot $K \subset \mathbf{S}^3$, {\it cf.} \cite{OoguriV,Taubes,Koshkin}.
Furthermore, the number of five-branes in \eqref{theoryB} is equal to $|\lambda|$ (= the total number of boxes in the Young diagram $\lambda$),
whereas the K\"ahler modulus of $X$ is related to the variable $\a$ in the superpolynomial:

\begin{equation}
\a \; = \; \exp \left( - {\rm Vol} (\cp^1) \right) \,.
\label{NvolCP1}
\end{equation}

This summary was very telegraphic, of course, but the reader can always consult e.g. \cite{GS} for further details, references, and the outline of the relation between different ways of looking at this physical system \eqref{theoryB}.
As reviewed there, the system \eqref{theoryB} was studied in physics as well as in math literature from
many different angles, e.g. from the target space point of view \cite{GSV,GWalcher} that leads to a description of \eqref{HHBPS}
in terms of enumerative invariants of $(X,L_K)$, or from the viewpoint of the five-brane \cite{fiveknots}, or from the viewpoint
of the $\CN=2$ gauge theory on a 4-manifold $M_4$ with ramification \cite{DGH}.
However, only recently \cite{FGS2,FGSS} the system \eqref{theoryB} has been considered from
the viewpoint of the two-dimensional theory on $D$, that will be yet another
motivation for the discussion in Section~\ref{sec:inst-surf}.

Now we are ready to discuss the BPS states that contribute to the reduced and unreduced superpolynomials,
$\CP^{\lambda} (\a,q,t)$ and $\bar \CP^{\lambda} (\a,q,t)$, respectively.

As discussed e.g. in \cite{OoguriV,GSV,NagaoN}, there are two kinds of BPS particles that contribute to \eqref{HHBPS}
in the problem at hand: M2-branes supported on bordered surfaces embedded in $X$ with boundary on $L_K$,
and momentum modes ({\it a.k.a.} D0-branes).
The charge of the former is precisely the $\a$-grading, while the charge of the latter is the $q$-grading.
Both types of states can mix together (forming stable BPS particles that carry both charges) and both can carry
a non-trivial spin $j_3$, which becomes the third $t$-grading on \eqref{HHBPS}.
To summarize, the dictionary between three gradings on the colored HOMFLY homology
and the charges of the BPS particles is as follows:
\be
\begin{array}{rcll}
``\a-\text{grading}" & = & \beta \in H_2 (X,L_K) \cong \Z \; & \quad \text{(``winding number'')} \\[.2cm]
``q-\text{grading}" & = & n \in \Z \;  & \quad  \text{(``D0-brane charge'')} \\[.2cm]
``t-\text{grading}" & = & 2 j_3 \in \Z & \quad  \text{(spin)}
\end{array}
\label{aqtgradings}
\ee
Using this dictionary, we can identify every generator of \eqref{HHBPS} with the corresponding BPS particle.
For example, if we see a term $\a^m q^n t^k$ in the superpolynomial, the proper translation to the physics
language is
\be
\begin{array}{ccl}
\text{monomial} & \longleftrightarrow & \text{``a BPS particle with winding number } m \,, \\[.1cm]
\a^m q^n t^k \subset \CP^{\lambda} & & \; \; \; \text{D0-brane charge } n \,, \text{ and spin } \frac{k}{2} \text{ ''}
\end{array}
\ee
The particles that carry no winding number ($m=0$) are usually called simply D0-branes.

One well-known property of the BPS states in the present setup is that particles that carry non-trivial winding number $m=1$
are fermions, whereas D0-brane states (with $m=0$) are bosonic \cite{Refined,JafferisM}.
When translated to the language of knot homologies, this statement says that statistics of the generators
in the colored HOMFLY homology is correlated with the $\a$-grading!

Namely, the generators with trivial $\a$-grading are bosonic, whereas the generators with $\a$-grading equal to 1 are fermions.
This will play a key role in what follows.\\

\begin{figure}[htb]
\centering
\includegraphics[width=4.0in]{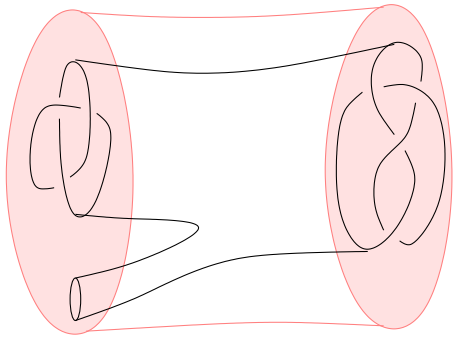}
\caption{A real isotropic knot cobordism $\Sigma \subset \R \times S^3$ defines a coassociative 4-manifold $L_{\Sigma} \subset \R \times T^* S^3$.}
\label{fig:cobordism2}
\end{figure}

Another useful lesson that we can learn from the connection with enumerative (BPS) invariants \eqref{HHBPS}--\eqref{theoryB}
is the behavior of colored HOMFLY homology with respect to cobordisms.
First, however, it is convenient to consider a dual brane system related to \eqref{theoryB}
by a geometric transition\footnote{sometimes also called ``large $N$ duality,'' which is somewhat of a misnomer
since the duality holds for all values of $N$, as long as one is careful with differentials and wall crossing \cite{DGR,GS}} \cite{GopakumarV,OoguriV}:
\bea
\text{space-time} & : & \qquad \R \times T^* W \times M_4 \nonumber \\
N~\text{M5-branes} & : & \qquad \R \times W \times D \label{theoryA} \\
|\lambda|~\text{M5-branes} & : &  \qquad \R \times L_K \times D \nonumber
\eea
where $W = S^3$ in most of our applications (or, its close cousin $W = \R^3$).
This system is relevant to $sl(N)$ homological knot invariants and has been studied from various points of view in \cite{AShakirov,IK,fiveknots}.
One advantage of the duality frame \eqref{theoryA} that will serve us well in building knot cobordisms is that the Lagrangian submanifold $L_K \subset T^* W$ can be defined simply as the conormal bundle to the knot $K \subset W$.

Now, it is easy to describe how knot cobordisms can be implemented in the physical setup.
Since the factor of $\R$ in \eqref{theoryB} and \eqref{theoryA} plays the role of the ``time'' direction,
one can replace $\R \times S^3$ by a general 4-manifold cobordism $V$ and replace $\R \times K$ embedded in it by
a more general knot cobordism $\Sigma \subset V$, as in Figure~\ref{fig:cobordism2}:
\bea
\text{space-time} & : & \qquad X_7 \times M_4 \nonumber \\
N~\text{M5-branes} & : & \qquad V \times D \label{theoryAcob} \\
|\lambda|~\text{M5-branes} & : &  \qquad L_{\Sigma} \times D \nonumber
\eea
Supersymmetry imposes certain constraints, which ensure that the physical setup here continues to
maintain the relation between supersymmetric configurations (BPS states) and functors on knot homologies.
Specifically, the 7-manifold $X_7$ must be a special holonomy manifold with holonomy group $G_2$,
while the 4-manifold $L_{\Sigma}$ must be a coassociative submanifold in $X_7$.
Both can be constructed in a fairly simple way, starting with the data of $V$ and $\Sigma$, respectively.

Just like $T^* W$ is a ``canonical'' Calabi-Yau 3-fold associated to a 3-manifold $W$, from a 4-manifold $V$
one can build a ``canonical'' $G_2$-holonomy manifold $X_7 = \Lambda^{2,+} (V)$ as the total space of the bundle of self-dual 2-forms.
This 7-manifold $X_7 = \Lambda^{2,+} (V)$ replaces $\R \times T^* W$ and reduces to it when $V = \R \times W$.
In what follows, we restrict our attention only to this case and, moreover, assume $W = S^3$.
Then, given a knot cobordism $\Sigma \subset V = \R \times S^3$,
we wish to construct a coassociative 4-manifold $L_{\Sigma} \subset \Lambda^{2,+} (V)$ with the property
\be
L_{\Sigma} \cap V = \Sigma
\ee
which is a suitable $G_2$-manifold version of the property $L_K \cap S^3 = K$ in \eqref{theoryA}.
By analogy with the conormal bundle construction of $L_K$, the desired coassociative 4-manifold $L_{\Sigma}$
can be defined as the subbundle of $\Lambda^{2,+} (V)$ over $\Sigma \hookrightarrow V$
such that restriction of the fiber to $\Sigma$ is trivial, provided that $\Sigma$ is a real isotropic minimal surface \cite{IKM,KMO}.
Details of this construction and its implications for functors on colored knot homologies will be discussed elsewhere.
Here, we simply point out that the existence of such geometric construction implies the existence of a map
\be
\CH^{\lambda} (\unknot) \otimes \CH^{\lambda} (K) \; \to \; \CH^{\lambda} (K)
\label{HHHunknot}
\ee
Indeed, applying geometric transition one more time to the setup \eqref{theoryAcob} with $V = \R \times S^3$ we get
a ``cobordism version'' of \eqref{theoryB}:
\bea
\text{space-time} & : & \qquad \R \times X \times M_4 \label{theoryBcob} \\
\text{M5-branes} & : &  \qquad \quad L_{\Sigma} \times D \nonumber
\eea
where $L_{\Sigma}$ is a coassociative 4-manifold in $\R \times X$.
When $\Sigma$ is a link cobordism from $K \bigsqcup \unknot$ to $K$,
the setup \eqref{theoryBcob} defines a functor on the space of refined open BPS states and, via \eqref{HHBPS},
leads to a functor \eqref{HHHunknot} on the colored HOMFLY homology.

\begin{figure}[htb]
\centering
\includegraphics[width=3in]{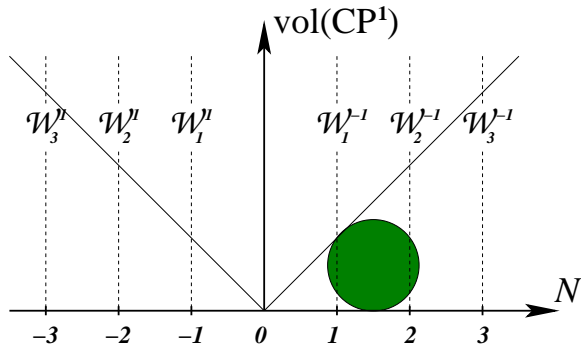}
\caption{Walls and chambers for the conifold $X$ in \eqref{theoryB}.}
\label{fig:walls}
\end{figure}

Note, that \eqref{HHHunknot} is very similar --- and, in fact, possibly related --- to the conjecture \cite{GS}
that $\CH^{\lambda} (K) = \CH_{\text{BPS}}^{\text{open}}$
is a module over the algebra of {\it closed} BPS states $\CH_{\text{BPS}}^{\text{closed}}$:
\be
\begin{array}{rcl}
\text{refined open BPS states}: &  & \CH_{\text{BPS}}^{\text{open}} \\
 & & \circlearrowleft \\
\text{refined closed BPS states}: &  & \CH_{\text{BPS}}^{\text{closed}}
\end{array}
\label{HHrep}
\ee
Even though we did not pursue the detailed comparison of
$\CA = \CH (\unknot)$ and the algebra $\CH_{\text{BPS}}^{\text{closed}}$ here,
it is worth mentioning that the two share some similarities.
For instance, much like $\CH (\unknot)$ has $sl(N)$ specializations indexed by $N$
and realized by the action of differentials $d_N$ (see below), the space $\CH_{\text{BPS}}^{\text{closed}}$
for the problem at hand also has an infinite set of chambers, conveniently labeled by $N$, see Figure \ref{fig:walls}.
Moreover, in the $N$-th chamber the space $\CH_{\text{BPS}}^{\text{closed}}$
has $N$ stable BPS objects \cite{NagaoN}: D0/D2 bound states of charge $(n,\beta) = (j,1)$, where $0 \le j < N$,
nicely matching the fact that $\CH^{sl(N), \tableau{1}} (\unknot)$ is also $N$-dimensional.

\subsection{Colored polynomials of the unknot}
\label{sec:unknot}

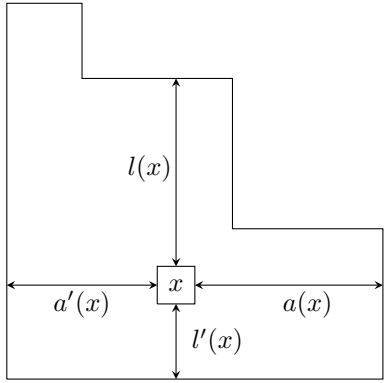
\begin{figure}
\begin{tikzpicture}
\draw (0,0)--(0,5)--(1,5)--(1,4)--(3,4)--(3,2)--(5,2)--(5,0)--(0,0);
\draw (2,1)--(2,1.5)--(2.5,1.5)--(2.5,1)--(2,1);
\draw (2.25,1.25) node {$x$};
\draw [<->,>=stealth] (2.5,1.25)--(5,1.25);
\draw [<->,>=stealth] (0,1.25)--(2,1.25);
\draw [<->,>=stealth] (2.25,1.5)--(2.25,4);
\draw [<->,>=stealth] (2.25,0)--(2.25,1);
\draw (4,1) node {$a(x)$};
\draw (1,1) node {$a'(x)$};
\draw (1.9,2.8) node {$l(x)$};
\draw (2.8,0.5) node {$l'(x)$};
\end{tikzpicture}
\caption{Arm, leg, co-arm and co-leg}
\label{armlegs}
\end{figure}

For a Young diagram $\lambda$ and a box $x$ in it we define the arm $a(x)$, leg $l(x)$, co-arm $a'(x)$ and co-leg $l'(x)$ by the Figure \ref{armlegs}. The hook-length of $x$ is defined as
$h(x)=a(x)+l(x)+1,$ and the content of $x$ is defined as $c(x)=a'(x)-l'(x)$.

Let $\lambda$ be a Young diagram with $k$ boxes, let $V_{\lambda}$ denote the irreducible representation of $S_k$ labelled by $\lambda$,
and let $U_{\lambda}=\mathbb{S}_{\lambda}(\C^{N})$ denote the irreducible representation of $sl(N)$ labelled by $\lambda$.
The $sl(N)$ quantum invariant of the $\lambda$-colored unknot coincides with the  character of the representation $U_{\lambda}$ and is given by (\ref{plambda}). This character is given by the hook formula:
$$\mbox{\rm ch} U_{\lambda}=\prod_{x\in \lambda}\frac{(1-q^{N+c(x)})}{(1-q^{h(x)})}.$$
If we substitute $A=q^{N}$, we get
\begin{equation}
\label{plambda}
P^{\lambda}(a,q)=\prod_{x\in \lambda}\frac{(1-A q^{c(x)})}{(1-q^{h(x)})}.
\end{equation}
Let $M_{\lambda}(q,t;p_i)$ be a Macdonald polynomial corresponding to a diagram $\lambda$. Recall that its evaluation can be computed by the formula \cite{macd19}:
\begin{equation}
\label{mstar}
M^{*}_{\lambda}(A,q,t)=M_{\lambda}\left(p_i=\frac{1-A^{i}}{1-t^{i}}\right)=\prod_{x\in\lambda}\frac{t^{l'(x)}-Aq^{a'(x)}}{1-q^{a(x)}t^{l(x)+1}}.
\end{equation}
Note that at $q=t$ the Macdonald polynomials degenerate to Schur polynomials and (\ref{mstar}) degenerates to the equation (\ref{plambda}).
In all approaches based on the refined BPS invariants the (unreduced) colored superpolynomials of the unknot
are given by the Macdonald dimensions
\be
\bar \CP^{\lambda} (\unknot) \; = \; M^*_{\lambda}(A,q,t) \,.
\label{mcdim}
\ee
In order to write these expressions with positive coefficients, one needs to convert to the variables
used e.g. in \cite{DGR}:
$$
A = -\a^2 t \,, \qquad
q = \q^2 \t^2 \,, \qquad
t = \q^2.
$$
The first few values of \eqref{mcdim} read as
\be
\bar \CP^{\tableau{1}} (\unknot) \; = \;
\frac{1 - A^{2}}{1 - t^{2}} =
\frac{1 + \a^2 \t}{1 - \q^2}
\ee
\be
\bar \CP^{\tableau{2}} (\unknot) \; = \;
\frac{(1 - A^{2})(1 - A^{2}q^{2})}{(1 - t^{2})(1 - q^{2} t^{2})}
= \frac{(1 + \a^2 \t) (1 + \a^2 \q^2 \t^3)}{(1 - \q^2) (1 - \q^4 \t^2)}
\ee
where we normalized by an overall powers of $A$, $q$, and $t$, so that
all expressions \eqref{mcdim} have the form of a power series $\bar \CP^{\lambda} (\unknot) = 1 + \ldots$.

{}From now on we will work in the variables $\a, \q, \t$. The variable $\t$ will be also denoted by $t_c$. The equation
(\ref{mcdim}) in these variables has a form:
\be
\label{mcdim1}
\bar \CP^{\lambda} (\unknot)=\prod_{x\in \lambda}\frac{1+\a^2\q^{2c(x)}t^{2a'(x)+1}}{1-\q^{2h(x)}\t^{2a(x)}}
\ee
Note that the right hand side is a polynomial in $\a$ and a series in $\q$ with all nonnegative coefficients. For $\lambda=S^r$, these answers match the stable homology of $(r,\infty)$ torus knots defined in \cite{DGR}.

\subsection{Schur-Weyl duality}

Following \cite{KRb}, let us model the fundamental representation of $sl(N)$ as $U_{\square}=\C[x]/(x^N).$
The Schur-Weyl duality (e.g. \cite{FH}) says that
$$U_{\square}^{\otimes k}=\bigoplus_{|\mu|=k}V_{\mu}\otimes U_{\mu},$$
or, equivalently,
\begin{equation}
\label{hom}
U_{\lambda}=\mbox{\rm Hom}_{S_k}(V_{\lambda},U_{\square}^{\otimes k})=\mbox{\rm Hom}_{S_k}(V_{\lambda},\C[x_1,\ldots,x_k]/I_{N}),
\end{equation}
where the ideal $I_N\subset \C[x_1,\ldots,x_k]$ is generated by $x_1^N,\ldots, x_k^N$. If $\lambda$ corresponds to the symmetric power of the fundamental representation, the construction
(\ref{hom}) gives the identity
$$S^{k}U_{\square}=(\C[x_1,\ldots,x_k]/I_{N})^{S_k}=\mathbb{C}[x_1,\ldots,x_k]^{S_k}/(I_{N})^{S_k}.$$
The symmetrized ideal $I_{N}^{S_k}$ is generated by $p_{N+s}=\sum x_i^{N+s}$ for $0\le s\le k-1$. If $\lambda$ corresponds to the exterior power of the fundamental representation, we have to study the anti-symmetric part of
$I_n$. One can check (\cite{GWalcher}) that $(I_N)^{\epsilon}$ is generated by $\Delta\cdot h_{N+1-j}, 1\le j\ge k$, where $h_i$ denote the complete symmetric polynomials, and $\Delta$ denotes the Vandermonde determinant. Therefore
$$\Lambda^{k}U_{\square}=(\mathbb{C}[x_1,\ldots,x_k]/I_{N})^{\epsilon}=\Delta \cdot (\C[x_1,\ldots,x_k])^{S_k}/(h_{N+1-j}).$$
Remark that the latter equation allows us to identify the exterior power with the cohomology of the Grassmannian:
$$\Lambda^{k}U_{\square}(sl(N))=\Delta\cdot H^{*}(Gr(k,N)),$$
since $h_{i}$ can be naturally identified with the Chern classes of the quotient bundle.

The setup of \cite{DGR} and \cite{GS} prescribes the knot homology to be obtained by a certain reduction procedure from the
``large $N$ limit''. Let us describe this construction for the case of unknot.
Given a Young diagram $\lambda$, there supposed to exist a triply graded space $\CH^{\lambda}$ modelling the $\lambda$-colored HOMFLY homology of the unknot. This space is equipped with the differentials $d_N$, and the homology of $\CH^{\lambda}$ with respect to $d_N$ should be isomorphic to the $sl(n)$ homology of the unknot colored by $\lambda$. Let us describe $\CH^{\lambda}$ along the lines of the above Schur-Weyl duality.
We start with the fundamental representation
$U_{\square}(sl(N))=\C[x]/(x^N).$
The corresponding triply graded space will be
$$\CH^{\square}=\Omega^{\bullet}(\C)=\C[x,dx],$$
and the differential comes from the Koszul complex:
$$d_{N}(x)=0,\quad d_{N}(dx)=x^{N}.$$
It is clear that $H^{*}(\CH^{\square}, d_N)\simeq U_{\square}(sl(N)).$
The $k$-th tensor power of this construction will give
$$\CH^{\square^{\otimes k}}=\Omega^{\bullet}(\C^{k})=\C[x_1,\ldots,x_k;dx_1,\ldots,dx_k],$$
and the differential is given by the formula
$d_N(x_i)=x_i^{N}$.
Therefore a natural candidate for $\CH^{\lambda}$ will be \cite{GWalcher}:
\begin{equation}
\label{hlambda}
\CH^{\lambda}:=\mbox{\rm Hom}_{S_k}\left(V_{\lambda},\Omega^{\bullet}(\C^{k})\right).
\end{equation}
Note that this space (as in the uncolored case) differs from the naive limit of $U_{\lambda}(sl(N)).$
It is naturally equipped with two gradings: $a$-grading is given by the degree of a differential form and $q$-grading is a standard
polynomial grading. In other words, the $(a,q)$-degree for $x_i$ equals to $q^{2i}$ and for $dx_i$ equals to $a^2$.

Let us check that the $(a,q)$-character of the space $\CH^{\lambda}$ agrees with the colored HOMFLY polynomial of the unknot. It follows from the results of \cite{KirPak} (see also \cite{Molch}) that the character of $\CH^{\lambda}$ equals 
\begin{equation}
\label{charUnknot}
\mbox{\rm ch}\CH^{\lambda}(a,q)=\prod_{(i,j)\in \lambda}\frac{(q^{i-1}+aq^{j-1})}{(1-q^{h(i,j)})}.
\end{equation}
Therefore
$\mbox{\rm ch}\CH^{\lambda}(a,q)=q^{\nu(\lambda)}P^{\lambda}(-a,q),$ where $\nu(\lambda)=\sum (i-1)\lambda_{i}.$
At $a=0$ the character (\ref{charUnknot})  degenerates to:
\begin{equation*}
P_{\text{bottom}}(\CH^{\lambda})=\prod_{(i,j)\in \lambda}\frac{q^{i-1}}{1-q^{h(i,j)}}.
\end{equation*}
For a symmetric representation we have
$$\CH^{S^{k}}=\left(\Omega^{\bullet}(\C^{k})\right)^{S_k}=\Omega^{\bullet}(\C^{k}/S_k)=\C[u_1,\ldots,u_k;du_1,\ldots,du_k],$$
where $u_i$ denote the elementary symmetric polynomials. Remark that the $(a,q)$-degree for $u_i$ equals to $q^{2i}$ and
for $du_i=\sum \frac{\partial u_i}{\partial x_j}dx_j$ equals to $a^2q^{2i-2}$. Therefore the bigraded Poincar\'e series $\CH^{S^{k}}$ has a form:
$$P(\CH^{S^{k}})=\prod_{i=1}^{k}\frac{1+a^2q^{2i-2}}{1-q^{2i}},$$
in agreement with (\ref{mcdim1}).

Let us describe the ``bottom'' $a=0$ part of the homology for $\lambda=\tableau{2 1}$.
The corresponding representation of $S_3$ is the standard two-dimensional representation.
We have to compute $\mbox{\rm Hom}_{S_k}\left(V_{\lambda},\C[x_1,x_2,x_3]\right)$.
Let $\mathfrak{a}\subset \C[x_1,x_2,x_3]$ denote the  ideal generated by the symmetric polynomials of
positive degree.
It is known that $\C[x_1,x_2,x_3]/\mathfrak{a}$ is isomorphic to the regular representation of $S_3$, and it contains two copies of $V_{\lambda}$
in gradings 1 and 2. These copies can be multiplied by arbitrary symmetric polynomials, so the character of the $a=0$ part of the homology
equals to
$$\frac{q+q^2}{ (1-q)(1-q^2)(1-q^3)}=\frac{q}{(1-q)^2(1-q^3)}.$$
This results agrees with the hook formula (\ref{plambda}) at $a=0$, since the hook-lengths for the boxes of our diagram are equal to 1, 1 and 3.

\subsection{Colored homology of the unknot}
\label{sec:alggrad}

Consider a free supercommutative algebra $\CA$ with even generators $u_1,\ldots,u_r$ and odd generators $\xi_1,\ldots,\xi_{s}$. Suppose that the algebra is graded such that
$\deg u_i=\alpha_{i},\ \deg \xi_i=\beta_{i}.$
It is clear that the Hilbert series of this algebra is given by the formula
\begin{equation}
\label{hilbert free}
H_{\alpha,\beta}(q)=\frac{\prod_{i=1}^{s}(1+q^{\beta_{i}})}{\prod_{i=1}^{r}(1-q^{\alpha_{i}})}.
\end{equation}
The similarity of this formula to (\ref{plambda}) and (\ref{mcdim1}) suggests the following construction.
We conjecture that the homology $\CH^{\lambda}$ of the unknot colored by the representation $\lambda$ is isomorphic to the
free graded supercommutative algebra $\CA$. To every box $x$ of $\lambda$ we associate an even generator $u_{x}$ and an odd generator $\xi_{x}$. The algebra has three gradings $q,a,t_{c}$ such that
$$(q,a,t_{c})[u_{x}]=(2h(x),0,2a(x)),\quad (q,a,t_{c})[\xi_{x}]=(2c(x),2,2a'(x)+1).$$
By (\ref{hilbert free}) and (\ref{mcdim1}) the Hilbert series of this algebra equals to
$$
\bar \CP^{\lambda}(\unknot)(a,q,t_{c})=\prod_{x\in \lambda}\frac{1+a^2q^{2c(x)}t^{2a'(x)+1}_{c}}{1-q^{2h(x)}t_{c}^{2a(x)}}.
$$

As we will mostly discuss the rectangle-shaped Young diagrams, let us fix some notations for them.
Let $\lambda$ be a rectangle with $R$ rows and $S$ columns. For a box $x=(i,j),\ 1\le i\le S, 1\le j\le R$ in it we have:
$$a(x)=S-i,\qquad l(x)=R-j,\qquad h(x)=R+S-i-j+1,$$
$$a'(x)=i-1,\qquad l'(x)=j-1.$$
To this box we associate a bosonic generator $u_{ij}$ and a fermionic generator $\xi_{ij}$ such that
$$(q,a,t_{c})[u_{ij}]=(2(R+S-i-j+1),0,2(S-i)),\quad (q,a,t_{c})[\xi_{ij}]=(2(i-j),2,2i-1).$$

We define a fourth grading using the mirror symmetry between $R\times S$ and $S\times R$ rectangles.
To write down the isomorphism $M_{\lambda} : \CH^{\lambda} \to \CH^{\lambda^t}$, one can define it on the level of generators.
We define
$$M_{\lambda} (u_{ij}):=u_{R-j+1,S-i+1},\quad M_{\lambda} (\xi_{ij}):=\xi_{ji}.$$

\begin{figure}
\begin{tikzpicture}
\draw (0,4)--(0,6)--(3,6)--(3,4)--(0,4);
\draw (1.2,5) node {$\bullet$};
\draw [<->,>=stealth] (0,5)--(1.1,5);
\draw [<->,>=stealth] (1.2,4)--(1.2,4.9);
\draw (0.6,5.2) node {$i$};
\draw (1.4,4.5) node {$j$};
\draw (2.8,4.2) node {*};

\draw [<->,>=stealth] (3.5,5)--(4.5,5);
\draw (4,5.2) node {$M$};

\draw (5,4)--(5,7)--(7,7)--(7,4)--(5,4);
\draw (6,5.8) node {$\bullet$};
\draw [<->,>=stealth] (6.1,5.8)--(7,5.8);
\draw [<->,>=stealth] (6,5.9)--(6,7);
\draw (5.9,6.4) node {$i$};
\draw (6.5,5.6) node {$j$};
\draw (6.8,4.2) node {*};

\draw (0,0)--(0,2)--(3,2)--(3,0)--(0,0);
\draw (1.2,1) node {$\bullet$};
\draw [<->,>=stealth] (0,1)--(1.1,1);
\draw [<->,>=stealth] (1.2,0)--(1.2,0.9);
\draw (0.6,1.2) node {$i$};
\draw (1.4,0.5) node {$j$};
\draw (0.2,0.2) node {*};

\draw [<->,>=stealth] (3.5,1)--(4.5,1);
\draw (4,1.2) node {$M$};

\draw (5,0)--(5,3)--(7,3)--(7,0)--(5,0);
\draw (6,1.2) node {$\bullet$};
\draw [<->,>=stealth] (5.9,1.2)--(5,1.2);
\draw [<->,>=stealth] (6,1.1)--(6,0);
\draw (5.9,0.5) node {$i$};
\draw (5.4,1.4) node {$j$};
\draw (5.2,0.2) node {*};

\end{tikzpicture}
\caption{Mirror symmetry for $u_{ij}$ (upper pair) and for $\xi_{ij}$ (lower pair).}
\label{mirrorrect}
\end{figure}
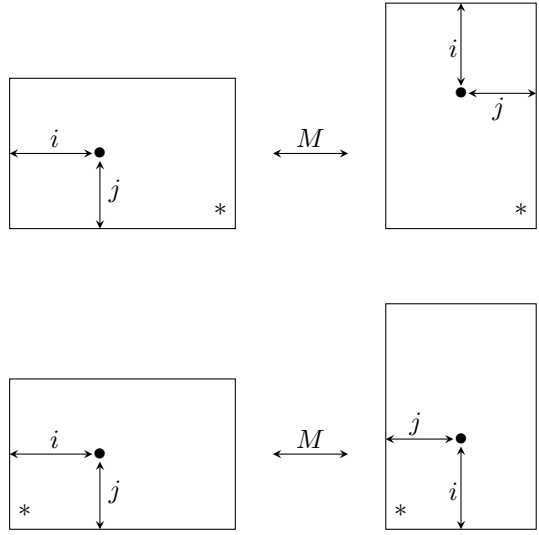

The correspondence between boxes is shown in Figure \ref{mirrorrect}.
Note that for bosonic generators we fix the southeast corner of the rectangle (marked with *) while
for the fermionic ones we fix the southwest corner.
We define the grading $t_{r}$ using the mirror symmetry relation $t_{r}(w)=t_{c}(M_{\lambda} (w))$, so
$$t_{r}(u_{ij})=2j-2,\qquad t_{r}(\xi_{ij})=2j-1.$$
Using the formula $$Q=\frac{q+t_{r}-t_{c}}{R},$$
we get
$$Q(u_{ij})=2,\quad Q(\xi_{ij})=0.$$

\subsection{Potentials and differentials: antisymmetric case}
\label{sec:potentanti}

It is known (\cite{DGR,Gornik,GWalcher,Lee}, see also Section \ref{sec:gradings}) that the reduction from HOMFLY homology to $sl(N)$ homology
is given by the following procedure. Let $\CH^{\lambda}(K)$ and $\CH^{sl(N),\lambda}(K)$
respectively denote the HOMLFY and $sl(N)$ homology of a knot $K$. There exists a differential $d_{N}$
of $(a,q)$-degree $(-1,N)$ such that
$$\CH^{sl(N), \lambda}(K)=H^{*}(\CH^{\lambda}(K),d_{N})$$
Remark that $d_{N}$ preserves the grading $Na+q$, and its Euler characteristic with respect to this grading equals
$P^{\lambda}(K;a=q^{N},q)=P^{sl(N), \lambda}(K).$

On the other hand, it is known that for some diagrams $\lambda$ one can construct a family of potentials
$W_{\lambda,N}(u_1,\ldots,u_d)$ such that $\CH^{sl(N), \lambda}(\unknot)$ coincides with the Jacobi ring $J=\C[u_i]/\left(\frac{\partial W_{\lambda,N}}{\partial u_i}\right)$.
More generally, $\CH^{sl(N), \lambda}(\unknot)$ is expected to coincide with the Hochschild homology of the category of the matrix factorizations of $W_{\lambda,N}$, that is, with the homology of the Koszul complex associated with the sequence of the partial derivatives of
$W_{\lambda,N}$. The latter homology contains the Jacobi ring $J$ and is equal to it if and only if $W_{\lambda,N}$ has only isolated singularities.
In other words, if the potentials $W_{\lambda,N}(u_1,\ldots,u_d)$ exist then one can consider the free supercommutative algebra $\CA$
with even generators $u_1,\ldots,u_d$,  odd generators $\xi_1,\ldots,\xi_{d}$ and the differential
$$d_{N}:\C[u_i,\xi_j]\to \C[u_i,\xi_j],\quad d_{N}(\xi_i)=\frac{\partial W_{\lambda,N}}{\partial u_i},\ d_{N}(u_i)=0.$$
This matches with the above construction of differentials, if one can identify $u_i$ and $\xi_i$ with the even and odd generators of
$\CH^{\lambda}$, so that $\C[u_i,\xi_j]=\CH^{\lambda}$.

The potential $W_{\Lambda^{k},N}$ was found in \cite{GWalcher}:
$$W_{\Lambda^{k},N}(u_1,\ldots,u_k)=\Coef_{N+1}\left[\ln (1+zu_1+\ldots+z^{k}u_k)\right].$$
It has an isolated singularity at the origin, and the Jacobi rings of $W_{\Lambda^{k},N}$ and $W_{\Lambda^{N-k},N}$ are isomorphic.
Moreover, one can check that up to some change of variables $(u_i)\to (v_i)$ the following identity holds:
\begin{equation}
\label{split of W}
W_{\Lambda^{k},2k-j}(u_1,\ldots,u_k)=W_{\Lambda^{k-j},2k-j}(v_1,\ldots,v_{k-j})+\frac{1}{2}\sum_{s=1}^{j}v_{k-j+s}v_{k+1-s}.
\end{equation}
For example, one can check that
$W_{\Lambda^{2},3}=-\frac{u_1^4}{4} + u_1^{2}u_2 - \frac{u_2^2}{2},W_{\Lambda^{1},3}=-\frac{u_1^{4}}{4},$
hence
$$W_{\Lambda^{2},3}=-W_{\Lambda^{1},3}-\frac{(u_2-u_1^2)^2}{2}.$$

\begin{remark}
Let us explain the construction of $v_{i}$ in general.
Consider supplementary variables $w_j$ defined by the equation
$$w_s=\Coef_j\left[\frac{1+zu_1+\ldots+z^{k}u_k}{1+zu_1+\ldots+z^{k-j}u_{k-j}}\right].$$
Clearly, $w_s=0$ for $s\le k-j$, and the change of variables
$$(u_1,\ldots,u_k)\to (u_1,\ldots, u_{k-j},w_{k-j+1},\ldots,w_{k})$$
is invertible. Now
$$W_{\Lambda^{k},2k-j}(u_1,\ldots,u_k)-W_{\Lambda^{k-j},2k-j}(u_1,\ldots,u_{k-j})=$$ $$\Coef_{2k-j+1}\ln\left[\frac{1+zu_1+\ldots+z^{k}u_k}{1+zu_1+\ldots+z^{k-j}u_{k-j}}\right]=$$
$$\Coef_{2k-j+1}\ln\left[(1+w_{k-j+1}z^{k-j+1}+\ldots+w_{k}z^{k}+\ldots)\right]=$$
$$\frac{1}{2}\sum_{s=1}^{j}w_{k-j+s}w_{k+1-s}+\ldots.$$
Since the quadratic part of this expression is non-degenerate, one can apply Morse lemma
with parameters and get rid of terms of higher degree in $w$'s.
See also \cite{GZV} for an alternative proof of this result and further discussions.
\end{remark}

The equation (\ref{split of W}) can be rewritten as $W_{\Lambda^{k},2k-j}=W_{\Lambda^{k-j},2k-j}+W'$, and we get the following diagram:
\begin{equation}
\begin{tabular}[c]{ccc}
$\CH^{\Lambda^{k}}$ & $\xrightarrow{d_{\Lambda^{k}\to \Lambda^{k-j}}}$ & $\CH^{\Lambda^{k-j}}$\cr
 $\vline$ & & $\vline$ \cr
${d_{2k-j}}$ & & ${d_{2k-j}}$\cr
 $\downarrow$ & & $\downarrow$ \cr
$\CH^{\Lambda^{k},2k-j}$ & $\xrightarrow{\qquad =\qquad}$ & $\CH^{\Lambda^{k-j},2k-j}$
\end{tabular}
\end{equation}

\begin{figure}
\begin{tikzpicture}
\draw (0,0)--(0,5)--(1,5)--(1,0)--(0,0);
\draw (0,4.5)--(1,4.5);
\draw (0,3)--(1,3);
\draw (0,2.5)--(1,2.5);
\draw (0,0.5)--(1,0.5);
\draw (-0.2,4.7) node {$2$};
\draw (-1,2.7) node {$2(k-j+1)$};
\draw (-0.3,0.3) node {$2k$};
\draw (0.5,4) node {$\vdots$};
\draw (0.5,1.5) node {$\vdots$};
\draw (0.5,4.7) node {$u_1$};
\draw (0.5,2.7) node {$u_{k-j+1}$};
\draw (0.5,0.3) node {$u_k$};
\draw (0.5,-1) node {$q(u_x)$};

\draw (4,0)--(4,5)--(5,5)--(5,0)--(4,0);
\draw (4,4.5)--(5,4.5);
\draw (4,3)--(5,3);
\draw (4,2.5)--(5,2.5);
\draw (4,0.5)--(5,0.5);
\draw (5.6,4.7) node {$2-2k$};
\draw (5.7,2.7) node {$2j-2k$};
\draw (5.3,0.3) node {$0$};
\draw (4.5,4) node {$\vdots$};
\draw (4.5,1.5) node {$\vdots$};
\draw (4.5,4.7) node {$\xi_1$};
\draw (4.5,2.7) node {$\xi_{j}$};
\draw (4.5,0.3) node {$\xi_k$};
\draw (4.5,-1) node {$q(\xi_x)$};

\draw [->,>=stealth] (3.8,4.7)--(1.2,2.7);
\draw (2.5,2.7) node {$\vdots$};
\draw [->,>=stealth] (3.8,2.7)--(1.2,0.3);
\end{tikzpicture}
\caption{Action of the differential $d_{\Lambda^{k}\to \Lambda^{k-j}}$ and the $q$-degrees of the generators $u_{x},\xi_{x}$.}
\label{dlambda}
\end{figure}
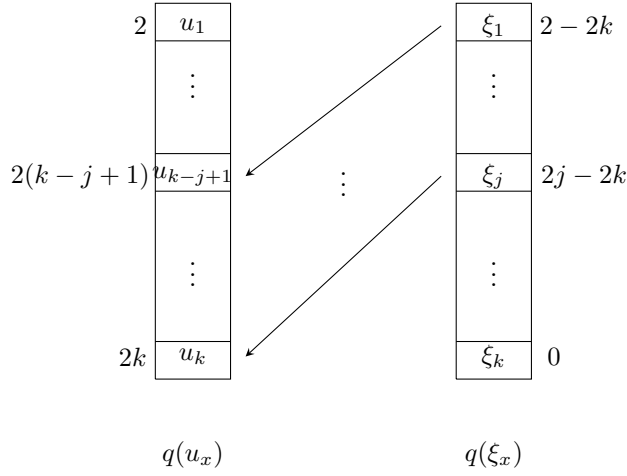

Here the colored differential $d_{\Lambda^{k}\to \Lambda^{k-j}}$ coincides with the one defined in \cite{GS}.
It is defined on HOMFLY homology colored by $\Lambda^{k}$, and its homology are isomorphic to the HOMFLY homology colored by $\Lambda^{k-j}$.
It corresponds to the ``remainder potential'' $W'=\frac{1}{2}\sum_{s=1}^{j}v_{k-j+s}v_{k+1-s}$ in (\ref{split of W}).
The corresponding Koszul complex is associated to the sequence of polynomials $v_{k-j+1},\ldots,v_{k}$:
$$d_{\Lambda^{k}\to \Lambda^{k-j}}(\xi_i)=v_{i+k-j},$$
and its homology is isomorphic to a free superpolynomial algebra generated by $v_1,\ldots,v_{k-j},\xi_{j+1},\ldots,\xi_n$.
We can easily see this differential on the picture with the bosonic and fermionic generators labeled by the boxes of a Young diagram.
Recall that their $q$-degrees are twice the contents and the hook-length of the boxes of $\Lambda^{k}$ respectively, as shown in Figure \ref{dlambda}. The differential $d_{\Lambda^{k}\to \Lambda^{k-j}}$ has $q$-degree $4k-2j$,
\be
(a,q) [d_{\Lambda^{k}\to \Lambda^{k-j}}] \; = \; (-2,4k-2j)
\ee

\subsection{Potentials and differentials: extension to $n\lambda$}
\label{sec:extension}

Let us first fix some notations. Let $\lambda=(\lambda_1,\ldots,\lambda_r)$ be a Young diagram. Then the boxes of
$n\lambda=(n\lambda_1,\ldots,n\lambda_r)$ are split into groups by $n$ boxes labelled by the boxes of $\lambda$.
If $x$ is a box of $\lambda$, then we label the corresponding boxes in $n\lambda$ by $x^{(1)},\ldots,x^{(n)}$.

For every $x$ in $\lambda$ we can also consider the generators $u_{x}$ and $\xi_{x}$ in
$\CH^{\lambda}$ and the generators  $u_{x}^{(i)}$ and $\xi_{x}^{(i)}$ in
$\CH^{n\lambda}$. For convenience, we can consider the generating functions
$$u_{x}(\tau)=\sum_{i=1}^{n-1} u_{x}^{(i)}\tau^{i-1},\ \xi_{x}(\tau)=\sum_{i=1}^{n-1} \xi_{x}^{(i)}\tau^{i-1}.$$

Suppose that the diagram $\lambda$ admits a potential $W_{\lambda,N}(u_1, \ldots, u_d)$. We conjecture that for all $n$
the diagrams $n\lambda$ admit potentials $W_{n\lambda,N}$ such that
\begin{equation}
\label{extension of W}
\hat{W}(\tau)=\sum_{n=1}^{\infty}\tau^{n-1}W_{n\lambda,N}=W_{\lambda,N}(u_1(\tau),\ldots,u_d(\tau)).
\end{equation}

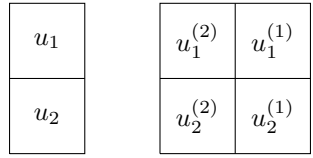
\begin{figure}
\begin{tikzpicture}
\draw (0,0)--(0,2)--(2,2)--(2,0)--(0,0);
\draw (0,1)--(2,1);
\draw (1,0)--(1,2);
\draw (0.5,0.5) node {$u_{2}^{(2)}$};
\draw (1.5,0.5) node {$u_{2}^{(1)}$};
\draw (0.5,1.5) node {$u_{1}^{(2)}$};
\draw (1.5,1.5) node {$u_{1}^{(1)}$};

\draw (-2,0)--(-2,2)--(-1,2)--(-1,0)--(-2,0);
\draw (-2,1)--(-1,1);
\draw (-1.5,0.5) node {$u_2$};
\draw (-1.5,1.5) node {$u_1$};
\end{tikzpicture}
\caption{Bosonic generators for $\Lambda^2$ and $2\Lambda^2$.}
\label{2x2}
\end{figure}

Since $W_{\Lambda^{2},3}(u_1,u_2)=-\frac{u_1^4}{4} + u_1^{2}u_2 - \frac{u_2^2}{2}$, we have
$$W_{2\Lambda^{2},3}(u_{1}^{(1)},u_{1}^{(1)},u_{2}^{(1)},u_{2}^{(2)})=\Coef_{1}\left[-\frac{u_1(\tau)^4}{4} + u_1(\tau)^{2}u_2(\tau) - \frac{u_2(\tau)^2}{2}\right]=$$
$$=-\left(u_{1}^{(1)}\right)^{3}u_{1}^{(2)} + \left(u_{1}^{(1)}\right)^{2}u_{2}^{(2)}+ 2u_{1}^{(1)}u_{1}^{(2)}u_{2}^{(1)} - u_{2}^{(1)}u_{2}^{(2)}.$$
Remark that $u_{i}^{(j)}$ correspond to the boxes of the $2\times 2$ rectangle (see Figure \ref{2x2}), so that their $q$-degrees are equal to the hook-lengths:
$$q(u_{1}^{(1)})=q(u_1)=2,\ q(u_{1}^{(2)})=4,\ q(u_{2}^{(1)})=q(u_{2})=4,\ q(u_{2}^{(2)})=6.$$
The potential $W_{2\Lambda^{2},3}$ is homogeneous of $q$-degree 10.

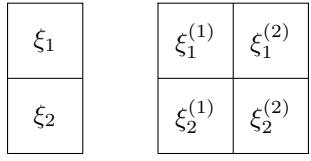
\begin{figure}
\begin{tikzpicture}
\draw (0,0)--(0,2)--(2,2)--(2,0)--(0,0);
\draw (0,1)--(2,1);
\draw (1,0)--(1,2);
\draw (0.5,0.5) node {$\xi_{2}^{(1)}$};
\draw (1.5,0.5) node {$\xi_{2}^{(2)}$};
\draw (0.5,1.5) node {$\xi_{1}^{(1)}$};
\draw (1.5,1.5) node {$\xi_{1}^{(2)}$};

\draw (-2,0)--(-2,2)--(-1,2)--(-1,0)--(-2,0);
\draw (-2,1)--(-1,1);
\draw (-1.5,0.5) node {$\xi_2$};
\draw (-1.5,1.5) node {$\xi_1$};
\end{tikzpicture}
\caption{Fermionic generators for $\Lambda^2$ and $2\Lambda^2$.}
\label{2x2ferm}
\end{figure}

Moreover, suppose that we have some differential $D$ on $\CH^{\lambda}$ given by the formula
$D(\xi_l)=p_{l}(u_1,\ldots,u_d).$ Such a differential can be naturally extended to $n\lambda$ by the formula
\begin{equation}
\label{extension of D}
D_n(\xi_{l}^{(i)})=\Coef_{n-1-i}[p_{l}(u_1(\tau),\ldots,u_d(\tau))]
\end{equation}
Remark that equations (\ref{extension of W}) and (\ref{extension of D}) agree with each other.
Indeed, every potential $W$ induced a differential $D_{W}$ such that $D_{W}(\xi_l)=\frac{\partial W}{\partial u_l}$,
if we match a bosonic variable $u_l$ with the fermionic variable $\xi_l$. We can match $u_{l}^{(i)}$ with $\xi_{l}^{(i)}$
so that
$$\frac{\partial W_n}{\partial u_{l}^{(i)}}=\Coef_{n-1}\left[\frac{\partial \hat{W}}{\partial u_{l}}\frac{\partial u_{l}(\tau)}{\partial u_{l}^{(i)}}\right]=\Coef_{n-1-i}[p_{l}(u_1(\tau),\ldots,u_d(\tau))].$$

As an example of this construction, consider $\lambda=(1)$, so that $n\lambda=(n)=S^{n}$.
We have $W_{\lambda,N}=x^{N+1}$, and it was suggested in \cite{GOR} (see also \cite{GORS})
that
$$W_{S^{n},N}=\Coef_{n-1}(x^{(1)}+x^{(2)}\tau+\ldots+x^{(n)}\tau^{n-1})^{N+1}.$$
This differential differs form the one proposed in \cite{GWalcher}, but its homology agree with the known categorifications of the Jones-Wenzl projectors \cite{CHKrushkal,CKrushkal,Stroppel}.

As a more interesting example, let us compute the colored differential for the rectangular $R\times S$ diagram.
Such a diagram can be presented as $S\cdot \Lambda^{R}$, and, as such, admits an extension of the colored differential
$d_{\Lambda^{R}\to \Lambda^{j}}$. We have
$d_{\Lambda^{R}\to \Lambda^{j}}(\xi_{l})=u_{l-j},$
so
$$d_{R\times S\to j\times S}(\xi_{l}^{(m)})=u_{l-j}^{(m)}.$$

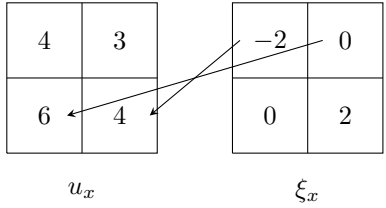
\begin{figure}
\begin{tikzpicture}
\draw (0,0)--(0,2)--(2,2)--(2,0)--(0,0);
\draw (0,1)--(2,1);
\draw (1,0)--(1,2);
\draw (0.5,0.5) node {$6$};
\draw (1.5,0.5) node {$4$};
\draw (0.5,1.5) node {$4$};
\draw (1.5,1.5) node {$3$};

\draw (3,0)--(3,2)--(5,2)--(5,0)--(3,0);
\draw (3,1)--(5,1);
\draw (4,0)--(4,2);
\draw (3.5,0.5) node {$0$};
\draw (4.5,0.5) node {$2$};
\draw (3.5,1.5) node {$-2$};
\draw (4.5,1.5) node {$0$};

\draw [->,>=stealth] (3.1,1.5)--(1.9,0.5);
\draw [->,>=stealth] (4.2,1.5)--(0.8,0.5);

\draw (1,-0.5) node {$u_{x}$};
\draw (4,-0.5) node {$\xi_{x}$};
\end{tikzpicture}
\caption{Colored differential $d_{2\Lambda^2\to 2\Lambda^1}$.}
\label{d22to21}
\end{figure}

In Figure \ref{d22to21} we show the action of the differential $d_{2\Lambda^2\to 2\Lambda^1}$.
The left and right rectangles illustrate the $q$-degrees of the bosonic and fermionic generators respectively.
The differential has $q$-degree 6 (the same as $d_{\Lambda^2\to \Lambda^1}$), and its homology is a free graded supercommutative algebra
isomorphic to $\CH^{S^2}$.

\section{Gradings, Differentials, Mirror Symmetry}
\label{sec:gradings}

\subsection{Introduction}
For any Young diagram $\lambda$, the $\lambda$-colored HOMFLY polynomial of a knot $K$, denoted $P^{\lambda}(K)(a,q)$, is a two-variable polynomial knot invariant whose $a=q^N$ one-variable specialization \eqref{PPN} equals the quantum knot invariant obtained by the representation $\lambda$ of the quantum $sl(N)$. Throughout the paper, $P^{\lambda}(K)$ denotes the reduced polynomial, i.e. it is normalized so that the value of the unknot is equal to 1.

We conjecture the existence of the homology theory $\CH^{\lambda}(K)$ categorifying the polynomials $P^{\lambda}(K)$, i.e. such that its graded Euler characteristic is equal to $P^{\lambda}(K)$. Although we do not give
the explicit definition of these homology theories, the insight from physics, representation theory and the double affine Hecke algebras (DAHAs) predict a very rich structure of differentials on $\CH^{\lambda}(K)$, together
with various symmetries of $\CH^{\lambda}(K)$. We focus on rectangular Young diagrams $\lambda$, i.e. Young diagrams which have the form of a rectangle with $R$ rows and $S$ columns, and which will be denoted $\lambda= R\times S$. This is a large class of representations, that includes as particular cases all symmetric representations ($R=1$) and antisymmetric ones ($S=1$).\footnote{Although many of the properties and results still hold for non-rectangular representations, there are certain subtleties in that case. We discuss this in section \ref{hooks}.}

Our main results for $\CH^{\lambda}(K)$ can be split in two groups, which match together perfectly.
First, we conjecture that certain, grading-independent structures exist and the following properties hold on $\CH^{\lambda}(K)$:
\begin{itemize}
\item $\CH^{\lambda}(K)$ is finite-dimensional.

\item There are two symmetries on $\CH^{\lambda}(K)$: one is the {\it mirror symmetry} that switches between $\lambda$ and its transpose partition $\lambda^t$:
\be
\CH^{\lambda}(K)\cong \CH^{\lambda^t}(K).
\ee
The other one is the {\it self-symmetry} on $\CH^{\lambda}(K)$.
While the mirror symmetry categorifies the known relation for the colored HOMFLY polynomials \eqref{PPmir}:
\be
P^{\lambda}(K)(a,q)=P^{\lambda^t}(K)(a,q^{-1}),
\ee
the self-symmetry is a completely new symmetry, that does not categorify any polynomial relation. It uses homological grading in a non-trivial way, and exists only on the homological level.

\item For every pair of non-negative integers $(n,m)$, such  that $n\ge R$ or $m\ge S$,\footnote{For other values of $m$ and $n$ the differentials may exist in the unreduced theory,
which we do not cover in details here.}
 there exists a collection of differentials $d_{n|m}$ on $\CH^{\lambda}(K)$, such that they pairwise anti-commute, and such that the homology of $\CH^{\lambda}(K)$ with respect to $d_{n|m}$ is
isomorphic to the homology $\CH^{sl(n|m),\lambda}(K)$ that categorifies the quantum invariant of a knot $K$
labeled by a representation $\lambda$ of $sl(n|m)$:
\be
H^*(\CH^{\lambda}(K),d_{n|m}) = \CH^{sl(n|m),\lambda}(K).
\ee
Introducing this two-parameter collection of differentials is one of the main goals
of the present paper\footnote{Although some of these operators were implicitly defined in \cite{G} for the uncolored homology, there relation to supergroups was not understood there.}. All previous results on the structural properties of the colored HOMFLY homologies (fundamental \cite{DGR} and symmetric ones  \cite{GS}) predicted only the existence of one-parameter family of anti-commuting differentials, with the parameter $N$ representing the rank of $sl(N)$ specializations.
As we shall explain, various reasons suggest that introduction of a two-parameter family associated with $sl(n|m)$ representations is a much more natural thing to do.

\item There exist {\it colored} differentials that change representations (colors). In other words, they allow passage from one homology theory $\CH^{\lambda}(K)$ to another one $\CH^{\mu}(K)$, for some Young diagram $\mu \subset \lambda$. More precisely, colored differential $d_{\lambda\to\mu}$ on $\CH^{\lambda}(K)$ has the property that the homology of $\CH^{\lambda}(K)$ with respect to $d_{\lambda\to\mu}$ is isomorphic to $\CH^{\mu}(K)$.

The main group of colored differentials, namely the group of $2(R+S-1)$ differentials,
is directly related to particular differentials $d_{n|m}$. More precisely, for every Young diagram $\mu$ that is obtained from $\lambda$ by removing either some of its rows or some of its columns,
there should be two different differentials of the form $d_{\lambda\to\mu}$ such that
\be\label{coliso}
H^*(\CH^{\lambda}(K),d_{\lambda\to\mu})=\CH^{\mu}(K).
\ee
They "lift" the relation  $\CH^{sl(n|m),\lambda}(K)\cong\CH^{sl(n|m),\mu}(K)$ coming from the isomorphism $\lambda \cong \mu$ as $sl(n|m)$ representations for appropriate values of $n$ and $m$, and we also denote such differential by $d^{col}_{n|m}$.  In other words, we have the following commutative diagram:

\begin{tabular}{ccc}
\quad\quad$\CH^{\lambda}(K)$ & $\xrightarrow{d_{\lambda\to \mu}\,\,\,\,(\mathrm{i.e.}\,\,\,\, d^{col}_{n|m})}$ & $\CH^{\mu}(K)$\cr
\quad\quad $\vline$ & & $\vline$ \cr
\quad\quad${\scriptstyle d_{n|m}}$ & & ${\scriptstyle d_{n|m}}$\cr
\quad\quad $\downarrow$ & & $\downarrow$ \cr
\quad\quad$\CH^{sl(n|m),\lambda}(K)$ & $\xrightarrow{\qquad \cong \qquad}$ & $\CH^{sl(n|m),\mu}(K)$
\end{tabular}

\noindent We note that on $\CH^{\lambda}(K)$ the differentials $d_{n|m}$ and $d^{col}_{n|m}$ share many similar properties. In particular, frequently the differential $d_{n|m}$ on $\CH^{\mu}(K)$ acts trivially, and then on $\CH^{\lambda}(K)$ the differentials $d_{n|m}$ and $d^{col}_{n|m}$ coincide.

\item For large classes of knots, including torus knots and two-bridge knots, the size of the homology $\CH^{\lambda}(K)$ is equal to the $|\lambda|$-th power of the dimension of $\CH^{\tableau1}(K)$. For an arbitrary knot $K$, such relation holds asymptotically  when $|\lambda|$ tends to infinity.

\end{itemize}

Our second major conjecture is that the colored HOMFLY homology $\CH^{\lambda}(K)$ of a knot $K$ is {\it quadruply-graded}. We  denote these four gradings by $a,q,t_{r},t_{c}$. Thus, apart from the $a$- and $q$-gradings associated to each generator, we also assign two homological $t$-gradings.
Any of the two $t$-gradings give a categorification of $P^{\lambda}(K)$, i.e. forgetting one of them gives triply-graded categorification of the colored HOMFLY-polynomial. In other words, if we denote the (four-variable) Poincar\'e polynomial of $\CH^{\lambda}(K)$ by $\CP^{\lambda}(K)(a,q,t_r,t_c)$, then we have
\be
\CP^{\lambda}(K) (a,q,t_{r}=-1,t_c=1) = \CP^{\lambda}(K) (a,q,t_r=1,t_{c}=-1) = P^{\lambda}(K)(a,q).
\ee


The appearance of the fourth grading is rather mysterious from the geometric point of view, but it seems to be inevitable. Let us list some of the evidences.
First of all, in the case of symmetric and anti-symmetric representations, the corresponding homology $\CH^{\lambda}(K)$ has been studied in the literature, and although the results in various papers are isomorphic, there are two different conventions for the homological $t$-gradings. The $t_r$-grading that we propose is the $t$-grading assigned to a generator of $\CH^{\lambda}(K)$
in the grading conventions of \cite{FGS2,GS}, whereas one can interpret $t_{c}$ as
the $t$-grading assigned to a generator of $\CH^{\lambda}(K)$ in grading conventions of \cite{AShakirov,DMMSS,DGR,GWalcher}.
As argued in \cite{GS}, the ``mirror symmetry'' --- the symmetry that relates the homology colored by a Young diagram and its transposed --- exchanges the two grading conventions. However, despite being very similar, one cannot pass from one convention to another by an explicit gradings transformation.
The solution that we propose is to keep track of both $t$-grading conventions together, and regard the homology $\CH^{\lambda}$ as
quadruply-graded theory, with four independent gradings $a,q,t_r$ and $t_c$.

The definite advantage of the quadruply-graded theory is that it makes
all of the structural properties and isomorphisms completely explicit. This cannot be achieved by a triply-graded theory: for example in the case of the symmetric representations (see \cite{GS}) only a few of the isomorphisms and symmetries could be made explicit.
In quadruply-graded theory $\CH^{\lambda}(K)$ one has completely explicit mirror symmetry and self-symmetry. Moreover, the
re-grading in all of the isomorphisms of colored differential (\ref{coliso})
can be made explicit, in all gradings. Finally, the exponential growth conjecture for the size of the colored homology becomes a totally refined exponential growth property of an appropriate three-variable Poincar\'e polynomial.
What is even more miraculous, is that all these symmetries, isomorphism re-gradings and exponential growth property have
rather simple and pleasant expressions. Moreover, all of them
become particularly elegant with the introduction of an auxiliary
grading, called the $Q$-grading, defined in a following simple way:
\be\label{Qgrad1}
Q=\frac{q+t_r-t_c}{R}.
\ee
This new grading, which can be considered as certain ``corrected'' $q$-grading,  cannot be seen on the decategorified, polynomial level, and as the formula indicates, both $t$-gradings are
needed for its definition. Therefore we need the quadruply-graded theory in order to have all explicit properties, and only a few of various symmetries and isomorphisms have polynomial counterparts at the decategorified level.

Finally, the symmetries get the nicest form when the homology is written in $(a,Q,t_r,t_c)$-gradings. Because of a linear relation between the gradings (\ref{Qgrad1}), this is just a simple re-grading of $\CH^{\lambda}(K)$. However, due to its importance, we give it a special name, $\tilde{\CH}^{\lambda}(K)$, and refer to it as the {\it tilde-version} of colored HOMFLY homology. Thus, it is given by:
\be
\tilde{\CH}^{\lambda}_{i,j,k,l}(K) \; := \; {\CH}^{\lambda}_{i,Rj-k+l,k,l}(K).
\ee
The only case when the two $t$-gradings coincide is the uncolored case, $\lambda=\tableau1$, and the resulting homology is
triply-graded in agreement with \cite{DGR}. In this case the $q$- and $Q$-gradings coincide, and the mirror symmetry and the self-symmetry are in fact the same symmetry.

Now we start with the detailed description of the structure of the quadruply-graded colored HOMFLY homology.

\subsection{Self-symmetry}
There are two remarkable symmetries in the colored homology when color is given by a rectangular Young diagram $\lambda$. The first one is self-symmetry on the homology ${\CH}^{\lambda}(K)$ and the second one is the mirror symmetry that relates the ${\CH}^{\lambda}(K)$ and ${\CH}^{\lambda^t}(K)$ homologies.

The main advantage of the $Q$-grading is that in the quadruple $(a,Q,t_r,t_c)$-grading both symmetries have a very simple and beautiful form. Thus, we present them in the corresponding tilde-version $\tilde{\CH}^{\lambda}(K)$.

\begin{conjecture}
The homology $\tilde{\CH}^{\lambda}(K)$ enjoys the following self-symmetry:
\be\label{selfsym}
{\boxed{\phantom{\int}\tilde{\CH}^\lambda_{i,j,k,l}(K) \cong  \tilde{\CH}^\lambda_{i,-j,k- Rj,l-Sj}(K).\phantom{\int}}}
\ee
\end{conjecture}

Moreover, we conjecture the existence of an involution $\Phi_{\lambda}:\CH^{\lambda}(K)\to\CH^{\lambda}(K)$
such that
$$(a,Q,t_r,t_c)[\Phi_{\lambda}(x)]=(a(x),-Q(x),t_{r}(x)-RQ(x),t_{c}(x)-SQ(x)).$$

Note that this symmetry does not descend to any relation on the colored HOMFLY polynomial (except in the uncolored case $\lambda=\tableau1$). In fact, due to a nontrivial presence of the $Q$-grading, it can be seen only on the quadruply-graded homology, and as to our knowledge, such symmetry has not been observed before.

For some knots, including torus knots (see below) we also conjecture that this involution comes from some version of the hard Lefshetz isomorphism. Namely, there exists an operator $L_{\lambda}:\CH^{\lambda}(K)\to\CH^{\lambda}(K)$
such that
\begin{equation}
\label{Llambda}
(a,q,Q,t_r,t_c)[L_{\lambda}]=(0,R+S,2,R,S),\quad \Phi_{\lambda}(x)=L_{\lambda}^{-Q(x)}.
\end{equation}
Note, according to the ``dictionary'' \eqref{aqtgradings} this operator is a bound state of $R+S$ D0-branes.

\subsection{Mirror symmetry}

Based on results for symmetric and anti-symmetric representations, it was conjectured \cite{GS} that colored HOMFLY homology enjoys mirror symmetry:
\begin{conjecture}(\cite{GS})\label{mirsymold}
For any knot, there exists an isomorphism (called {\it mirror symmetry}) between $\lambda$- and $\lambda^t$-colored HOMFLY homologies
preserving the $a$-grading and reversing the $q$-grading.
\end{conjecture}

However, to obtain the behavior of the $t$-gradings under this symmetry, one needs to introduce the fourth grading.
We are able to write an explicit formulas for the $t$-degrees change in the case of rectangular $\lambda$.
We start with a {\em different} conjectural symmetry between $\lambda$- and $\lambda^{t}$-colored homology.

\begin{conjecture}
In $(a,Q,t_{r},t_{c})$ gradings one has
\be\label{mirsym}
{\boxed{\phantom{\int}\tilde{\CH}^{\lambda^t}_{i,j,k,l} (K)\cong  \tilde{\CH}^\lambda_{i,j,l,k}(K).\phantom{\int}}}
\ee
\end{conjecture}

In order to obtain the explicit quadruply-graded version of the mirror symmetry from (\ref{mirsymold}), we combine self-symmetry (\ref{selfsym}) and mirror symmetry (\ref{mirsym}):

\begin{corollary}\label{mirsim}
In $(a,Q,t_{r},t_{c})$ gradings one has
\be\label{mirsymnew}
{\boxed{\phantom{\int}\tilde{\CH}^{\lambda^t}_{i,j,k,l}  (K)\cong  \tilde{\CH}^\lambda_{i,-j,l-Sj,k-Rj}(K).\phantom{\int}}}
\ee
\end{corollary}

The last equation enables one to write the explicit change of gradings in the mirror symmetry. Denote by $M_{\lambda}:\CH^{\lambda} \to \CH^{\lambda^t}$ the mirror symmetry isomorphism.
Consider the gradings of a generator $x$ in $\CH^{\lambda}$ theory and the gradings of the generator $M_{\lambda}(x)$ of $\CH^{\lambda^t}$, denoted by $(\hat{a},\hat{q},\hat{Q},\hat{t}_r,\hat{t}_c)$, where $\lambda$ is an $R\times S$ rectangular Young diagram and $\lambda^t$ is $S\times R$ rectangular Young diagram. They are related by the following transformation:
\begin{eqnarray*}
\hat{a} &=& a,~ \qquad \hat{q} = -q,~ \qquad \hat{Q}  =  -Q,\\
\hat{t}_r &=& t_c- S Q = (1+\frac{S}{R})t_c - \frac{S}{R} t_r - \frac{S}{R} q,\\
\hat{t}_c &=& t_r-  R Q = t_c- q.
\end{eqnarray*}

Since the mirror symmetry (\ref{mirsymnew}) inverts the $q$-grading, it ``categorifies'' the following relation between the colored HOMFLY polynomials:
\be
P^{\lambda}(K)(a,q)=P^{\lambda^t}(K)(a,q^{-1})
\ee
which is often attributed to level-rank duality.

\subsection{$sl(n|m)$ differentials}
\label{sec:supermn}

The homology $\CH^{\lambda}(K)$ comes with a large structure of pairwise anti-commuting differentials on it. We conjecture that there exists a differential for every pair of nonnegative integers $(n,m)$, such that $n\ge R$ or $m\ge S$. These differentials should be related to the homologies
corresponding to Lie superalgebras $sl(n|m)$ and representations labeled by Young diagram $\lambda$. The condition that $n\ge R$ or $m\ge S$ gives precisely those pairs $(n,m)$ for which the Lie superalgebra $sl(n|m)$ has an irreducible representation labeled by $\lambda$ (see Section  \ref{repsup}).
This generalizes the structure from \cite{DGR} and \cite{GS}, where the analogous differentials are labelled by a single integer $N$ and were related to the $sl(N)$ representations, i.e. $sl(N)$ specializations of the $\lambda$-colored HOMFLY polynomial. As we shall argue in this paper, the structure corresponding to Lie superalgebras is more natural, and the single parameter $N$ from \cite{DGR} and \cite{GS} is just a super-rank of a corresponding Lie superalgebra:
\be
N = n - m.
\ee

\vskip 0.3cm

This structure has appeared in the physics literature \cite{VafaUNM,WittenK} as a symmetry of brane/anti-brane systems,
which makes it especially easy to implement it in the brane system \eqref{theoryB} or, better yet, in the dual system \eqref{theoryA}
where the rank $N$ enters as the number of five-branes.
In this framework, replacing $sl(N)$ by a Lie superalgebra $sl(n|m)$ can be easily implemented by replacing $N$ five-brane
with a system of $n$ five-branes and $m$ anti-branes.

\begin{figure}[htb]
\centering
\includegraphics[width=3.0in]{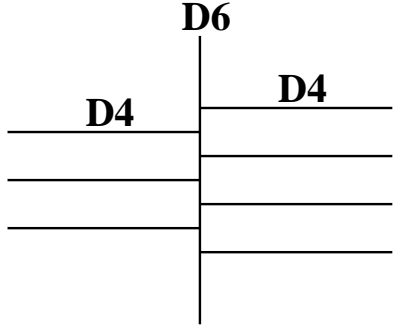}
\caption{In the approach of \cite{fiveknots}, homological invariants of knots labeled by representations of $sl(n|m)$
can be realized by a system of $n$ semi-infinite D4-branes ending on one side of the D6-brane and $m$ semi-infinite D4-branes ending on the other side of the D6-brane.}
\label{fig:D4D6}
\end{figure}

Another, more supersymmetric way to realize quantum and homological invariants of knots labeled by $sl(n|m)$ representations
is to consider $n$ five-branes supported on $\R \times W \times D$ exactly as in \eqref{theoryA} and introduce
$m$ new five-branes supported on $\R \times W \times D'$, such that $\C \cong D' \subset M_4$ and $D \cap D' = \{ 0 \}$,
\bea
\text{space-time} & : & \qquad \R \times T^* W \times M_4 \nonumber \\
n~\text{M5-branes} & : & \qquad \R \times W \times D \label{theoryAmn} \\
m~\text{M5-branes} & : &  \qquad \R \times W \times D' \nonumber
\eea
For instance, if $M_4 \cong \C^2$ is parametrized by coordinates $(z_1,z_2)$, then we can choose $D =\{ z_1 = 0 \}$ and $D' = \{ z_2 = 0 \}$.
Furthermore, identifying $M_4$ with the Taub-NUT space and reducing on the circle fiber takes us from M-theory
to type IIA string theory with a system of D4-branes and D6-branes shown in Figure \ref{fig:D4D6}.

Then, as in \cite{fiveknots}, homological $sl(n|m)$ knot invariants can be obtained by ``counting''
solutions to the following partial differential equations on a 5-manifold $\R \times W \times \R_y$:
\begin{align}
F^+ - \frac{1 }{ 4} B \times B - \frac{1 }{ 2} D_y B & = 0 \notag \\
F_{y i} + D^{j} B_{j i} & = 0 \notag
\end{align}
with the prescribed behavior near the ``interface'' $y=0$:
$$
B \simeq \frac{B_0^{(\pm)} }{ y} \quad \text{as} \quad y \to \pm 0 \,.
$$
Here, $B$ is a section of $\Omega^{2,+} (\R \times W) \otimes {\rm ad} (E)$, so that
$$
(B \times B)_{ij} = \sum_{k} [B_{ik} , B_{jk}]
$$
Relegating further study of this physical system to a future work, we now return to the discussion of mathematical structures that we expect to see in such physical realizations of homological $sl(n|m)$ knot invariants.\\

\begin{conjecture}
For every pair of nonnegative integers $(n,m)$ such that $n\ge R$ or $m\ge S$ there exists a differential $d_{n|m}$ on the homology $\CH^{\lambda}(K)$, satisfying the following properties:

\begin{itemize}
\item {\bf Specialization:}
\be
(\CH^{\lambda}(K),d_{n|m})|_{a=q^{n-m}} \cong \CH^{sl(n|m),\lambda}(K).
\ee

\item  {\bf Anticommutativity:} All $d_{n|m}$ pairwise anti-commute.

\end{itemize}

\end{conjecture}

The homology theory $\CH^{sl(n|m),\lambda}(K)$  categorifies $(n-m)$-specialization of the $\lambda$-colored HOMFLY polynomial. In other words:
\be
\CP^{sl(n|m),\lambda}(K)(q,t=-1)=P^{\lambda}(K)(a=q^{n-m},q).
\ee
Although the right hand side of the above relation equals the quantum polynomial invariant associated with $sl(n|m)$ representations labeled by $\lambda$, it depends only on the super-rank of the Lie superalgebra. However, the homologies $\CH^{sl(n|m),\lambda}(K)$ for various pairs $(n,m)$ with fixed $n-m$ can in general be non-isomorphic. That is one of the reasons why we expect to have differentials parametrized by two integers instead of just one, as it was predicted for the fundamental representation \cite{DGR} and for the
symmetric representations \cite{GS}.

The main reason for the existence of such structure comes from the representation theory of the Lie superalgebras $sl(n|m)$. The condition that $n\ge R$ or $m\ge S$ gives precisely those pairs $(n,m)$ for which the Lie superalgebra $sl(n|m)$ has an irreducible representation labeled by $\lambda$. We denote the subset of the $(x,y)$ plane satisfying $x\ge 0$, $y\ge 0$ and $x\ge R$ or $y\ge S$ by $\Omega$. Then, the set of all admissible pairs $(n,m)$ for which we have a differential $d_{n|m}$ consists exactly of the points from $\Omega$ whose both coordinates are integers.

As for the degrees of the differentials $d_{n|m}$, we have the following:
\be
(a,q)[d_{n|m}]=(-2,2(n-m)),
\ee
in agreement with the specialization $a=q^{n-m}$. The $t_c$-degree is equal to
\be
t_c[d_{n|m}]=-2m-1.
\label{tcdegnm}
\ee
The $t_r$-degree is more subtle, and is given in the next section; in particular, it depends on $R$ and $S$.

\begin{figure}[htb]
\centering
\includegraphics[width=4.0in]{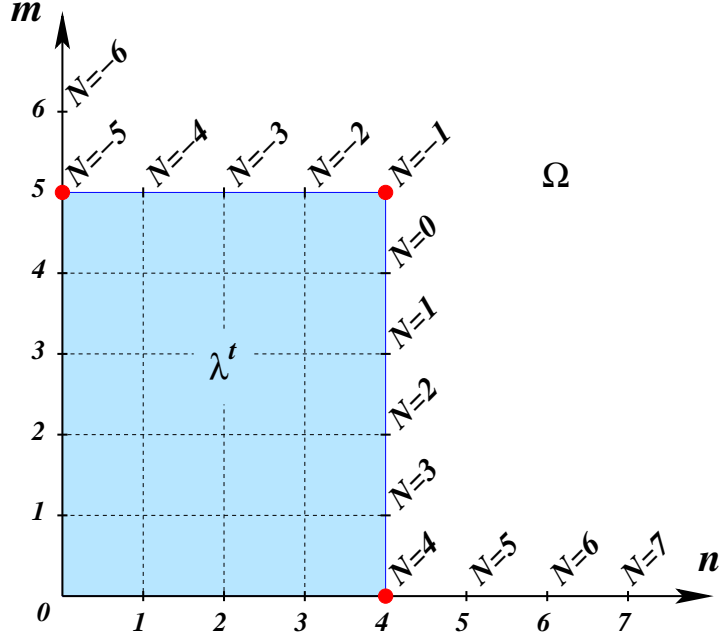}
\caption{Domain $\Omega$ can be visualized as the complement of $\lambda^t$ in a positive quadrant of the $(n,m)$ plane.
Integer points on its boundary correspond to differentials $d_{n|m}$ or, equivalently, $d_N$ with $N=n-m$, among which three are special.
They correspond to the three corners of the domain $\Omega$ and are marked by red bullets.
The differentials for this examples are shown in Figure \ref{sunflower}.}
\label{fig:domain}
\end{figure}

In  the case of the fundamental \cite{DGR} and the symmetric representations \cite{GS}, the conjectured differentials $d_N$ were parametrized by a single integer $N$. It was argued that they correspond to the $sl(N)$ representations for $N>0$, while for the other values of $N$ there was no clear motivation. However, when observed from the $sl(n|m)$ point of view, they become much more natural: $d_N$ becomes in fact precisely the differential $d_{n|m}$, where $n$ and $m$ are uniquely determined by requiring that $N=n-m$ and such that $(n,m)$ lies on the boundary of $\Omega$,
see Figure \ref{fig:domain}.

In the case of the fundamental representation, the $t$-degrees of the differentials $d_N$ were defined in a rather artificial way.
By using 
the conventions from above, we get directly the desired degrees. (Note that in the fundamental representations the two $t$-gradings coincide: $t_r=t_c=t$.)
Indeed, for $N>0$, we have $d_N=d_{N|0}$, and therefore according to \eqref{tcdegnm} their $t$-degrees are equal to $-1$. For $N=0$, we have that $d_0=d_{1|1}$,  so the $t$-degree of $d_0$ equals $-3$, while for $N>0$, we have that $d_{-N}=d_{0|N}$, and so the $t$-degree of  $d_{-N}$ equals $-2N-1$, exactly as predicted in \cite{DGR}.

Although there are infinitely many integer pairs in $\Omega$, for a given knot $K$ only finitely many of the differentials $d_{n|m}$ can act non-trivially. Indeed, since $\CH^{\lambda}(K)$ is finite-dimensional and since $q(d_{n|m})=2(n-m)$ and $t_c(d_{n|m})=-2m-1$,
the differential $d_{n|m}$ becomes void on $\CH^{\lambda}(K)$ for sufficiently large values of $n$ or $m$. Thus we have:

\begin{conjecture}\label{con1}
For a knot $K$ and a rectangular Young diagram $\lambda$, we have
\begin{equation}\label{coneq1}
\CP^{sl(n|m),\lambda}(K)(q,t_r,t_c)=\CP^{\lambda}(K)(a=q^{n-m},q,t_r,t_c),
\end{equation}
when either $n$ or $m$ is sufficiently large.
\end{conjecture}

\subsection{Colored differentials}\label{coldif}

For every knot $K$, the homology theory $\CH^{\lambda}(K)$ comes with a remarkable structure of the so-called {\it colored} differentials that allow passage to homology theories $\CH^{\mu}(K)$ for various partitions $\mu \subset \lambda$. In other words, for certain such partitions $\mu$, there exists a differential $d$ on $\CH^{\lambda}(K)$ such that the homology with respect to this differential is isomorphic to $\CH^{\mu}(K)$, i.e.
\be\label{prob1}
H^*(\CH^{\lambda}(K),d) \cong \CH^{\mu}(K).
\ee
As we shall describe, the possible partitions $\mu$ will be the ones that can be obtained by removing either some rows or some columns from $\lambda$.
All these differentials are in fact ``lifting'' $sl(n|m)$ differentials $d_{n|m}$ for some specific values of $(n,m)$ that are on the boundary of the region $\Omega$, see Figure \ref{fig:domain}. These are precisely the pairs $(n,m)$ such that $\lambda\cong\mu$ as $sl(n|m)$ representations (see section \ref{repsup}). We denote the corresponding colored differential by $d^{col}_{n|m}$, and the ``lifting'' property is given by the following commutative triangle:

\begin{tabular}[c]{ccc}
\quad\quad\quad$\CH^{\lambda}(K)$ & $\xrightarrow{\quad d^{col}_{n|m}\quad}$ & $\CH^{\mu}(K)$\cr
\quad\quad\quad $\vline$ & & $\vline$ \cr
\quad\quad\quad${\scriptstyle d_{n|m}}$ & & ${\scriptstyle d_{n|m}}$\cr
\quad\quad\quad $\downarrow$ & & $\downarrow$ \cr
\quad\quad\quad$\CH^{sl(n|m),\lambda}(K)$ & $\xrightarrow{\qquad \cong \qquad}$ & $\CH^{sl(n|m),\mu}(K)$
\end{tabular}

\medskip
Here the arrow going from $\CH_1$ to $\CH_2$ labeled by $d$ means that a differential $d$ on $\CH_1$ is such that $H^*(\CH_1,d)\cong \CH_2$. In particular, note that the two differentials, both denoted $d_{n|m}$, are essentially different since they are defined on different homology theories. The differentials $d_{n|m}$ and $d^{col}_{n|m}$ are very similar, in particular, the four gradings of $d_{n|m}$ and $d^{col}_{n|m}$ are identical.

We show that for each rectangular diagram $\mu$ that is obtained from $\lambda$ by removing either some rows or some columns, there exist two different differentials on $\CH^{\lambda}$ realizing the isomorphism (\ref{prob1}). As we explained, all these colored differentials are of the form $d^{col}_{n|m}$. However, in order to clearly indicate the color-changing property of the
differential, we also label them by a superscript ``$\pm$'' as in $d^{\pm}_{\lambda\to\mu}$. One of them will be denoted as $d^{+}_{\lambda\to\mu}$: it has positive $Q$-degree (equal to $+2$), and we will be referred to as the {\em positive differential}. The second differential has negative $Q$-degree (equal to $-2$) and will be denoted by $d^{-}_{\lambda\to\mu}$.

We will show that $d^{+}_{\lambda\to\mu}$ and $d^{-}_{\lambda\to\mu}$ are exchanged by the involution $\Phi_{\lambda}$. In a particular case when $\mu=\emptyset$, the homology $\CH^{\mu}$ is one-dimensional and the corresponding colored differentials are called {\it canceling} differentials.
All these differentials are deep generalizations of the differentials $d_{1}$ and $d_{-1}$ in the uncolored theory \cite{DGR}, which can be considered as colored differentials to the trivial representation $ \emptyset$. Positive colored differentials can be seen as analogs of $d_1$ and negative colored differentials can be seen as analogs of $d_{-1}$.

Now, let us give a more precise definition of colored differentials, together with the explicit change of gradings in the isomorphism (\ref{prob1}). We note that the introduction of the fourth grading is again crucial: the $Q$-grading is indispensable for the explicit re-gradings. 

Let $\Sigma$ denote the $S$-invariant of a knot $K$, i.e. it is even integer such that the homology of the uncolored homology $\CH^{\tableau1}(K)$ with respect to the canceling differential $d_1$, which is one-dimensional, has $(\a,q,t_r,t_c)$-degree $(\Sigma,-\Sigma,0,0)$, see \cite{DGR,slice}.

\textbf{Positive row-removing differentials:} For any $0\le k< R$
there exists a ``row-removing differential'' $d^{+}_{R\times S\to k\times S}$ (which is in fact $d^{col}_{R+k|0}$) such that the homology
of $\CH^{R\times S}(K)$ with respect to it is isomorphic to $\CH^{k\times S}(K)$.
In particular, for $k=0$
the differential $d^{+}_{R\times S\to \emptyset}$ is canceling.
These differentials have the following degrees:
\be
(a,q,t_r,t_c)[d^{+}_{R\times S\to k\times S}]=(-2,2R+2k,-2k-1,-1).
\ee
Since $Q=\frac{q+t_r-t_c}{R}$, we have $Q[d^{+}_{R\times S\to k\times S}]=2.$
The explicit grading change in the corresponding isomorphism (\ref{prob1}) is as follows: let $x$ be a generator of $\CH^{k\times S}(K)$ with degree $(\a,Q,t_r,t_c,q)$. Then the degrees of the corresponding generator $\phi(x)$ of $H^*(\CH^{R\times S}(K),d^{+}_{R\times S\to k\times S})$, denoted  by $(\hat{\a},\hat{Q},\hat{t}_r,\hat{t}_c,\hat{q})$, are given by:
\begin{eqnarray*}
\hat{\a}=\a+S(R-k)\Sigma,\qquad\hat{Q}=Q-S(R-k)\Sigma,\\
\hat{t}_r=t_r+(R-k)Q + Sk(R-k) \Sigma,\\
\hat{t}_c=t_c,\qquad\hat{q}=q-S(R-k)(R+k)\Sigma.
\end{eqnarray*}

\textbf{Positive column-removing differentials: } For any $0\le l< S$
there exists a ``column-removing differential'' $d^{+}_{R\times S\to R\times l}$ (which is in fact $d^{col}_{R|l}$) such that the homology
of $\CH^{R\times S}(K)$ with respect to it is isomorphic to $\CH^{R\times l}(K)$.
In particular, for $l=0$
the differential $d^{+}_{R\times S\to \emptyset}$ is canceling and coincides with the canceling differential above.
These differentials have the following degrees:
$$(a,q,t_r,t_c)[d^{+}_{R\times S\to R\times l}]=(-2,2R-2l,-1,-2l-1).$$
One can check that $Q[d^{+}_{R\times S\to R\times l}]=2.$
 The grading change formulas have the form:
\begin{eqnarray*}
\hat{\a}=\a+R(S-l)\Sigma,\qquad \hat{Q}=Q-R(S-l)\Sigma,\\
\hat{t}_r=t_r,\qquad \hat{t}_c=t_c+(S-l)Q+Rl(S-l)\Sigma,\\
\hat{q}=q+(S-l)Q+R(S-l)(l-R)\Sigma.
\end{eqnarray*}

\textbf{Negative row-removing differentials: }
For any $0\le k< R$
there exists a ``row-removing differential'' $d^{-}_{R\times S\to k\times S}$ (which is in fact $d^{col}_{k|S}$) such that the homology
of $\CH^{R\times S}(K)$ with respect to it is isomorphic to $\CH^{k\times S}(K)$.
In particular, for $k=0$
the differential $d^{-}_{R\times S\to \emptyset}$ is canceling.
These differentials have the following degrees:
$$(a,q,t_r,t_c)[d^{-}_{R\times S\to k\times S}]=(-2,2k-2S,-2k-2R-1,-2S-1).$$
One can check that $Q[d^{-}_{R\times S\to k\times S}]=-2.$
The grading change formulas have the form:
\begin{eqnarray*}
\hat{\a}=\a+S(R-k)\Sigma,\quad \hat{Q}=Q+S(R-k)\Sigma,\\
\hat{t}_r=t_r+S(R-k)(R+k)\Sigma,\quad \hat{t}_c=t_c+S^2(R-k)\Sigma,\\
\hat{q}=q+kQ+S(R-k)(S-k)\Sigma.
\end{eqnarray*}

\textbf{Negative column-removing differentials: } For any $0\le l< S$
there exists a ``column-removing differential'' $d^{-}_{R\times S\to R\times l}$ (which is in fact $d^{col}_{0|S+l}$) such that the homology
of $\CH^{R\times S}(K)$ with respect to it is isomorphic to $\CH^{R\times l}(K)$.
In particular, for $l=0$
the differential $d^{-}_{R\times S\to \emptyset}$ is canceling and coincides with the canceling negative differential above.
These differentials have the following degrees:
$$(a,q,t_r,t_c)[d^{+}_{R\times S\to R\times l}]=(-2,-2l-2S,-2R-1,-2l-2S-1).$$
One can check that $Q[d^{+}_{R\times S\to R\times l}]=-2.$
The grading change formulas have the form:
\begin{eqnarray*}
\hat{\a}=\a+R(S-l)\Sigma,\quad \hat{Q}=Q+R(S-l)\Sigma,\\
\hat{t}_r=t_r+R^2(S-l)\Sigma,\quad \hat{t}_c=t_c+R(S-l)(S+l)\Sigma,\\
\hat{q}=q+R(S-l)(S+l)\Sigma.
\end{eqnarray*}

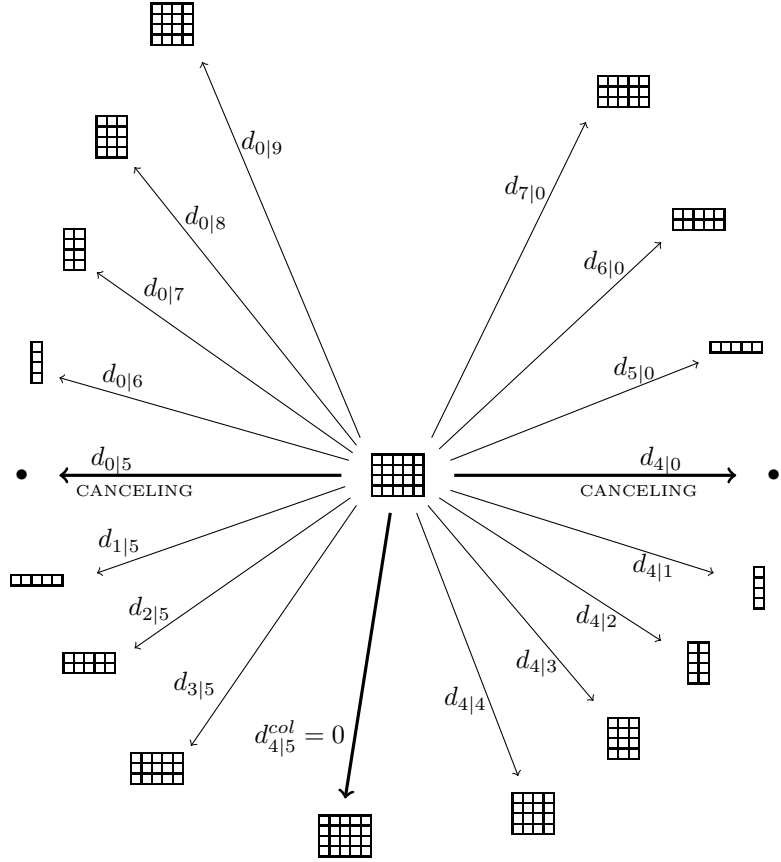
\begin{figure}
\begin{tikzpicture}

\draw (0,0) node {$\tableau{5 5 5 5}$};

\draw[->] (-0.5,0.5)--(-2.6,5.5);
\draw[->] (-0.55,0.4)--(-3.5,4.1);
\draw[->] (-0.6,0.3)--(-4,2.7);
\draw[->] (-0.65,0.2)--(-4.5,1.3);

\draw[->,very thick] (-0.75,0)--(-4.5,0);

\draw[->] (-0.7,-0.15)--(-4,-1.3);
\draw[->] (-0.63,-0.3)--(-3.5,-2.3);
\draw[->] (-0.55,-0.4)--(-2.75,-3.6);

\draw[->,very thick] (-0.1,-0.5)--(-0.7,-4.3);

\draw[->] (0.25,-0.5)--(1.6,-4);
\draw[->] (0.4,-0.4)--(2.6,-3);
\draw[->] (0.55,-0.3)--(3.5,-2.2);
\draw[->] (0.7,-0.2)--(4.2,-1.3);

\draw[->,very thick] (0.75,0)--(4.5,0);

\draw[->] (0.7,0.2)--(4,1.5);
\draw[->] (0.55,0.35)--(3.5,3.1);
\draw[->] (0.45,0.5)--(2.5,4.7);

\draw (-3,6) node {$\tableau{4 4 4 4}$};
\draw (-3.8, 4.5) node {$\tableau{3 3 3 3}$};
\draw (-4.3, 3) node {$\tableau{2 2 2 2}$};
\draw (-4.8,1.5) node {$\tableau{1 1 1 1}$};

\draw (-5,0) node {$\bullet$};

\draw (-4.8,-1.4) node{$\tableau{5}$};
\draw (-4.1,-2.5) node{$\tableau{5 5}$};
\draw (-3.2,-3.9) node{$\tableau{5 5 5}$};

\draw (-0.7,-4.8) node {$\tableau{5 5 5 5}$};

\draw (1.8,-4.5) node{$\tableau{4 4 4 4}$};
\draw (3,-3.5) node{$\tableau{3 3 3 3}$};
\draw (4,-2.5) node{$\tableau{2 2 2 2}$};
\draw (4.8,-1.5) node{$\tableau{1 1 1 1}$};

\draw (5,0) node {$\bullet$};

\draw (4.5,1.7) node {$\tableau{5}$};
\draw (4,3.4) node {$\tableau{5 5}$};
\draw (3,5.1) node {$\tableau{5 5 5}$};

\draw (-1.8,4.4) node {$d_{0|9}$};
\draw (-2.55,3.4) node {$d_{0|8}$};
\draw (-3.1,2.45) node {$d_{0|7}$};
\draw (-3.65,1.3) node {$d_{0|6}$};
\draw (-3.8,0.2) node {$d_{0|5}$};
\draw (-3.5,-0.2) node {{\tiny{CANCELING}}};
\draw (-3.7,-0.9) node {$d_{1|5}$};
\draw (-3.3,-1.8) node {$d_{2|5}$};
\draw (-2.7,-2.8) node {$d_{3|5}$};
\draw (-1.3,-3.5) node {$d^{col}_{4|5}=0$};
\draw (0.9,-3) node {$d_{4|4}$};
\draw (1.85,-2.5) node {$d_{4|3}$};
\draw (2.65,-1.9) node {$d_{4|2}$};
\draw (3.4,-1.23) node {$d_{4|1}$};
\draw (3.5,0.2) node {$d_{4|0}$};
\draw (3.2,-0.2) node {{\tiny{CANCELING}}};
\draw (3.15,1.4) node {$d_{5|0}$};
\draw (2.75,2.8) node {$d_{6|0}$};
\draw (1.7,3.8) node {$d_{7|0}$};

\end{tikzpicture}
\caption{Colored differentials for the $4\times 5$ rectangle.}
\label{sunflower}
\end{figure}

\vskip 0.25 cm

To explain better such a large structure of differentials, we present them in Figure \ref{sunflower}. This diagram represents the case of the representation labeled by the Young diagram
with $R=4$ rows and $S=5$ columns.
An arrow labelled by $d_{n|m}$ going from Young diagram $D_1$ to $D_2$ means that there exists a differential $d^{col}_{n|m}$ on $\CH^{D_1}(K)$ and  $\CH^{D_2}(K)$, such that
\be\label{glav}
H^*(\CH^{D_1}(K),d^{col}_{n|m}) \cong \CH^{D_2}(K).
\ee
There are three ''special" differentials which are represented by thick lines: two of them ($d^{col}_{4|0}$ and $d^{col}_{0|5}$) are canceling, i.e. the homology with respect to them is one-dimensional, and the third one  (the middle one) is  trivial  $d^{col}_{4|5}=0$.

In a general case of an arbitrary $R \times S$ rectangular Young diagram $\lambda$, the picture would be of the same form with $2R+2S-1$ differentials: written in a clockwise direction, they will  be $d_{R+i|0}$, with $i$ ranging $i=R-1,\ldots,0$, $d_{R|i}$, with $i$ ranging $i=0,\ldots,S$, $d_{i|S}$, with $i$ ranging $i=R-1,\ldots,0$, and $d_{0|S+i}$, with $i$ ranging $i=0,\ldots,S-1$, The canceling differentials will be the differentials $d_{R|0}$ and $d_{0|S}$, while the middle one would be $d_{R|S}$.\footnote{Although this middle colored differential $d^{col}_{R|S}$ is trivial, the $sl(R|S)$ colored differential $d_{R|S}$ is in general non-trivial and is related to knot Floer homology, in the same way as the differential $d_0$ in the uncolored case \cite{DGR}. We discuss these HFK-like differentials in Section \ref{sec:symant}.}

In the upper half of the diagram are the positive row- and negative column-removing differentials, with the row-removing ($d_{R+i|0}$ with $0\le i< R$) being in the upper right part and column-removing ($d_{0|S+i}$ with $0\le i< S$) being in the upper left part of the diagram.
The negative row- and positive column-removing differentials are in the lower half of the picture (below the canceling differentials).
The row-removing ones ($d_{i|S}$ with $0\le i<R$) are in the lower left part and column-removing ($d_{R|i}$ with $0\le i<S$) are in the lower right part of the diagram. The positive colored differentials are in the right half of the diagram, while the negative colored differentials are in the left half.\\

We conjecture that the collection of differentials is self-contained under the mirror symmetry and self-symmetry isomorphism. Namely,
Let $\lambda$ be a rectangular Young diagram and let $\mu$ be obtained from $\lambda$ by removal of some of the rows or columns. Then the following identities hold:
\begin{eqnarray}
\Phi_{\lambda}(d^{+}_{\lambda\to\mu}(x))&=&d^{-}_{\lambda\to\mu}(\Phi_{\lambda}(x)),\\
M_{\lambda}(d^{\pm}_{\lambda\to\mu}(x))&=&d^{\pm}_{\lambda^t\to\mu^t}(M_{\lambda}(x))).
\end{eqnarray}
The involution $\Phi$ provides an isomorphism
\be
H^{*}(\CH^{\lambda}(K),d^{+}_{\lambda\to\mu}(x)) \quad \stackrel{\Phi}{\longrightarrow} \quad H^{*}(\CH^{\lambda}(K),d^{-}_{\lambda\to\mu}(x)).
\ee

\subsection{Representation theory interpretation of the colored differentials}\label{repsup}

The main explanation for the existence of all these positive and negative colored differentials with the wanted behavior comes from the representation theory of $sl(n)$ and  $sl(n|m)$. We include the classical Lie algebra $sl(n)$ in the family of Lie superalgebras $sl(m|n)$, by identifying it with $sl(n|0)$.

Indeed, it is well-known that the irreducible polynomial representations of $sl(n)$ are labeled with the Young diagrams with at most $n$ rows. Moreover, some of these representations are isomorphic: namely, for every $k$ with $0\le k <R$, the representations labeled by the rectangular Young diagrams $R \times S$ and $k\times S$ are isomorphic as $sl(R+k)$-representations. This in turn implies that the doubly-graded homologies $\CH^{sl(R+k),R\times S}$ and $\CH^{sl(R+k),k\times S}$ categorifying respectively the quantum polynomial invariants $P^{sl(R+k),R\times S}(q)$ and $P^{sl(R+k),k\times S}(q)$, should be isomorphic. Finally, since $\CH^{sl(R+k),R\times S}$ is isomorphic to the homology of $\CH^{R\times S}$ with respect to the differential $d_{R+k|0}$, we obtain the explanation of the wanted behavior of the positive row-removing differentials.

As for the positive column-removing differentials one should look for the representations of the Lie superalgrebras $sl(n|m)$. There is an analogous theory
for the polynomial representations of $sl(n|m)$ \cite{br1,br2}, extending the classical case of $sl(n)$ algebras. This time, the irreducible representations are labeled by Young diagrams such that the $(n+1)$-st row contains at most $m$ boxes. Again, some of these representations are isomorphic: for every $k$ with $0\le k <S$ we have that the representations labeled by the rectangular Young diagrams $R\times S$ and $R\times k$ are isomorphic as $sl(R|k)$-representations.
Now one gets the homology theory corresponding to the $(sl(R|k), R\times S)$-representations as the homology of $\CH^{R\times S}$ with respect to $d_{R|k}$.
Hence, the above isomorphisms of the $sl(R|k)$-representations give rise to the positive column-removing differentials.

As for the negative differentials, they can be obtained in the same way as the positive differentials by using ``dual'' isomorphisms from above, combined with the mirror symmetry. Indeed, for two Young diagrams $\lambda$ and $\mu$ we have the following isomorphism of representations:
\be
\lambda \cong \mu, \quad\quad {\textrm{ as} }\quad  sl(n|m)\quad  {\textrm{ representations}},
\ee
if and only if
\be
\lambda^t \cong \mu^t, \quad\quad {\textrm{ as }}\quad sl(m|n) \quad {\textrm {representations}}.
\ee
In such way the existence of the negative row-removing differentials $d^{col}_{k|S}$ with $0\le k <R$ is a consequence of the isomorphism of $sl(k|S)$-representations corresponding to Young diagrams $R\times S$ and $k\times S$.
Finally, the existence of negative column-removing differentials $d^{col}_{0|S+k}$ with $0\le k <S$ follows from the isomorphism of $sl(0|S+k)$-representations corresponding to Young diagrams $R\times S$ and $R\times k$.

\subsection{Universal colored differentials}
\label{unicolor}

There exists yet another set of colored differentials on $\CH^{\lambda}(K)$ such that the homology with respect to any of them is isomorphic to $\CH^{\lambda'}(K)$ for some $\lambda'$ obtained from $\lambda$ by removing some of its rows or columns. The difference with respect to the differentials $d^{col}_{n|m}$ from above is that they are {\it universal} in the sense that their $\a$-degree is equal to 0.

Although the insights from the theory of deformations of potentials imply that such universal colored differentials should  exists for any row- or column-removal, we have managed to find explicitly such differentials only in the cases when $\lambda$ has two rows or two columns.

\begin{conjecture}\label{unicol}
Let $K$ be a knot.
Let $\lambda$ be $2 \times S$ rectangular diagram and let $\lambda'$ be $1\times S$ rectangular diagram. Then there exists a differential $d^{\uparrow}$ on  $\CH^{\lambda}(K)$ of $\a$-degree $0$, such that the homology of $\CH^{\lambda}(K)$ with respect to $d^{\uparrow}$ is isomorphic to $\CH^{\lambda'}(K)$.

Let $\mu$ be $R \times 2$ rectangular diagram and let $\mu'$ be $R\times 1$ rectangular diagram. Then their exists a differential $d^{\leftarrow}$ on $\CH^{\mu}(K)$ of $\a$-degree $0$ such that the homology of $\CH^{\mu}(K)$ with respect to $d^{\leftarrow}$ is isomorphic to $\CH^{\mu'}(K)$.
\end{conjecture}


To be more precise, we give  explicit values for the degrees of differentials $d^{\uparrow}$ and $d^{\leftarrow}$ and the degree changes in the isomorphisms of the Conjecture \ref{unicol}. We describe everything in the $(\a,Q,t_r,t_c,q)$-gradings.
The degrees of the differentials are given by:
\begin{eqnarray}
\deg d^{\uparrow}&=&(0,0,-2,0,2),\\
\deg d^{\leftarrow}&=&(0,0,0,2,2).
\end{eqnarray}
The first isomorphism from Conjecture \ref{unicol} gives:
\be\label{isom5}
H^*(\CH^{\lambda}(K),d^{\uparrow}) \cong \CH^{\lambda'}(K).
\ee
Let $x$ be a generator of $\CH^{\lambda'}(K)$ with degree $(\a,Q,t_r,t_c,q)$. Then the degrees of the corresponding generator $\phi(x)$ of $H^*(\CH^{\lambda}(K),d_N)$, denoted  by $(\hat{\a},\hat{Q},\hat{t}_r,\hat{t}_c,\hat{q})$, are given by:
\begin{eqnarray*}
\hat{\a}=2\a,\quad \hat{Q}=2Q,\quad \hat{t}_r=4t_r,\\
\hat{t}_c=2t_c,\quad \hat{q}=4q-2t_c.
\end{eqnarray*}
The second isomorphism from Conjecture \ref{unicol} gives:
\be\label{isom6}
H^*(\CH^{\mu}(K),d^{\leftarrow}) \cong \CH^{\mu'}(K).
\ee
Let $x$ be a generator of $\CH^{\mu'}(K)$ with degree $(\a,Q,t_r,t_c,q)$. Then the degrees of the corresponding generator $\phi(x)$ of $H^*(\CH^{\mu}(K),d_N)$, denoted  by $(\hat{\a},\hat{Q},\hat{t}_r,\hat{t}_c,\hat{q})$, are given by:
\begin{eqnarray*}
\hat{\a}=2\a,\quad \hat{Q}=2Q,\quad \hat{t}_r=2t_r,\\
\hat{t}_c=4t_c,\quad \hat{q}=2q+2t_c.
\end{eqnarray*}

\subsection{Refined exponential growth}

Now we are ready to formulate a quantitative refinement of the ``exponential growth conjecture'' from \cite{GS}.

\begin{conjecture}
Let $K$ be either two-bridge knot or a torus knot.
Let $\tilde{\CP}^{\lambda}(K)(a,Q,t_{r})$ denote the Poincar\'e polynomial (in $(a,Q,t_{r})$ gradings, i.e. after setting $t_c=1$) of $\tilde{\CH}^{\lambda}(K)$. Let $\tilde{\CP}^{\Lambda^{R}}(K)(a,Q,t_{r})$ denote the Poincar\'e polynomial of $\tilde{\CH}^{\Lambda^R}(K)$ ($R$-th antisymmetric colored homology).
Then the following relation holds:
\begin{equation}
\label{reduced refined exp growth}
\tilde \CP^{\lambda}(K)(a,Q,t_{r})=\left(\tilde \CP^{\Lambda^{R}}(K)(a,Q,t_{r})\right)^{S}.
\end{equation}
\end{conjecture}
By using the mirror symmetry and by switching $t_r \longleftrightarrow t_c$, we get the dual statement:

\begin{conjecture}
Let $K$ be either two-bridge knot or a torus knot.
Let $\tilde{\CP}^{\lambda}(K)(a,Q,t_{c})$ denote the Poincar\'e polynomial in $(a,Q,t_{c})$ gradings of $\tilde{\CH}^{\lambda}(K)$. Let $\tilde{\CP}^{S^{S}}(K)(a,Q,t_{c})$ denote the Poincar\'e polynomial of $\tilde{\CH}^{S^S}(K)$ ($S$-th symmetric colored homology).
Then the following relation holds:
\begin{equation}
\label{reduced refined exp growth dual}
\tilde{\CP}^{\lambda}(K)(a,Q,t_{c})=\left(\tilde{\CP}^{S^{S}}(K)(a,Q,t_{c})\right)^{R}.
\end{equation}
\end{conjecture}
\begin{corollary}
The dimensions of the ``rectangular homology groups'' satisfy the equation
\begin{equation}
\label{exp growth}
\dim \CH^{\lambda}(K)=\left(\dim \CH^{\tableau{1}}(K\right))^{|\lambda|}.
\end{equation}
Moreover, the similar statement holds in $(a,Q)$-gradings, i.e. with both $t_r$ and $t_c$ specialized to 1:
\begin{equation}
\label{a q exp growth}
\tilde{\CP}^{\lambda}(K)(a,Q)=\left(\tilde{\CP}^{\tableau{1}}(K)(a,Q)\right)^{|\lambda|}.
\end{equation}
\end{corollary}

\begin{corollary}
The Poincar\'e polynomial of the symmetric homology in $(a,Q,t_{r})$ gradings can be deduced from the uncolored
Poincar\'e polynomial (recall that $t_{r}=t_{c}=t$ for the uncolored homology):
\begin{equation}
\label{reduced refined exp growth for S}
\tilde{\CP}^{S^{r}}(K)(a,Q=q,t_{r}=t,t_c=1)=\left(\CP^{\tableau1}(K)(a,q,t)\right)^{r}.
\end{equation}
The Poincar\'e polynomial of the antisymmetric homology in $(a,Q,t_{c})$ gradings can be deduced from the uncolored
Poincar\'e polynomial:
\begin{equation}
\label{reduced refined exp growth for Lambda}
\tilde{\CP}^{\Lambda^{r}}(K)(a,Q=q,t_{r}=1,t_c=t)=\left(\CP^{\tableau1}(K)(a,q,t)\right)^{r}.
\end{equation}
\end{corollary}

Going back to original $q$-degrees, i.e. the ordinary version of the colored homology $\CH^{\lambda}(K)$, the above relations reduce to two-variable (in variables $a$ and $t$) exponential growths:
\be
{\CP}^{S^{r}}(K)(a,q=1,t_{r}=t,t_c=1)=\left(\CP^{\tableau1}(K)(a,q=1,t)\right)^{r},
\ee
and
\be
{\CP}^{\Lambda^{r}}(K)(a,q=1,t_{r}=1,t_c=t)=\left(\CP^{\tableau1}(K)(a,q=1,t)\right)^{r}.
\ee
Although such nice (and surprising) exponential growth properties should exist for large classes of knots, it is not expected that it should hold in general, i.e. for every knot. For example, the uncolored homology of the knot $9_{42}$ is 9-dimensional (see \cite{DGR}), its (triply-graded) $S^2$-colored homology is computed and the result has 401 generators. Therefore, we don't have any exponential growth property for it.\footnote{This could also be seen on the level of the $S^2$-colored HOMFLY polynomial of $9_{42}$. It has nonzero terms with $a$-degree $-6$, while the uncolored homology is concentrated in $a$-degrees -2,0 and 2. We thank Paul Wedrich for providing us with his computations of the colored HOMFLY polynomials.}
Nevertheless, such exponential growth property should hold in general but only asymptotically, as we briefly explain in the following subsection.

\subsection{Relation to the super-$A$-polynomial}

The exponential growth of the colored HOMFLY homology summarized
in \eqref{reduced refined exp growth}--\eqref{reduced refined exp growth for S}
nicely complements the results of \cite{FGS1,FGS2,Nawata},
where the large-$r$ limit of the $S^r$ homology was studied in detail.
In particular, it was conjectured in \cite{FGS2} that for an arbitrary knot $K$
the $S^r$-colored superpolynomials obey a recursion relation that comes from quantization of
an algebraic curve $A^{\text{super}} (x,y;a,t)=0$.
Specifically, this recursion relation has the following general form:
\be
\CP^{S^{r+n}} (K;a,q,t) + a_{n-1} \, \CP^{S^{r+n-1}} (K;a,q,t)
 + \ldots 
 + a_0 \, \CP^{S^{r}} (K;a,q,t) \; = \; 0
\label{QVCsuper}
\ee
and can be conveniently expressed in the operator form $\hat A^{\text{super}} \; \CP^{S^*} = 0$,
where the coefficients $a_i \equiv a_i (x, a, q, t)$ are rational functions of $a$, $q$, $t$, and $x \equiv q^{r}$.
To make contact with \eqref{a q exp growth}, we need to set $q=1$ which, by definition,
reduces the operator $\hat A^{\text{super}}$ 
to its characteristic polynomial, the super-$A$-polynomial $A^{\text{super}} (x, y; a,t)$.

Furthermore, since setting $q=1$ also implies $x=q^r=1$, we conclude that making contact with \eqref{a q exp growth}
involves comparison with the super-$A$-polynomial evaluated as $x=1$. Something very nice happens in this limit;
namely, in every example studied so far, all coefficients $a_i (x;a,t)$ of the super-$A$-polynomial
vanish, except for $a_{n-1}$. As a result, \eqref{QVCsuper} reduces to a much simpler recursion relation:
\be
\CP^{S^{r+n}} (K;a,q=1,t) \; + \; a_{n-1} (x=1, a, t) \, \CP^{S^{r+n-1}} (K;a,q=1,t) \; = \; 0 \,,
\ee
which indeed takes the form \eqref{a q exp growth}.
Specifically, it can be consistent with \eqref{a q exp growth}
if and only if $a_{n-1} (x=1, a, t) = \CP (K;a,q=1,t)$.
We conjecture that this is always the case, therefore, providing a bridge between our present work and \cite{FGS1,FGS2,Nawata}:

\begin{conjecture}
For any knot $K$, the following relation between the super-$A$-polynomial $A^{\text{super}} (x,y;a,t)$
and the (uncolored) superpolynomial holds:
\be
A^{\text{super}} (x=1,y;a,t) \; = \; y^k + y^{k-1} \, \CP (K;a,q=1,t) \,.
\ee
\end{conjecture}

\subsection{Symmetric and antisymmetric representations}\label{sec:symant}

In \cite{GS} the triply-graded colored HOMFLY homology for the symmetric and antisymmetric representations have been computed for plenty of knots. We obtain a quadruply-graded theory by adding a $t$-degree of the corresponding generator of the colored homology for the dual representation. The $t$-grading we used in \cite{GS} (and called ``old'') corresponds to $t_r$-grading in this paper. The other $t_c$-grading coincides with the $t$-grading used in \cite{AShakirov,DMMSS,DGR,GWalcher}. In the case of  the symmetric representations, we have two additional properties of the quadruply-graded colored homology: $\delta$-grading and HFK-like differential.
To every generator $x$ of $S^r$-colored quadruply-graded theory we can associate a $\delta$-grading by:
\be
\delta(x)=a(x)+\frac{q(x)}{2}-\frac{t_r(x)+t_c(x)}{2}=a(x)+\frac{Q(x)}{2}-t_r(x).
\ee

We call a knot $K$  $S^r$-thin if all generators of $\CH^{S^r}(K)$ have the same $\delta$-grading. In the uncolored case, $r=1$, this coincides with the definition of a thin knot. Also, in the uncolored case, the $\delta$-grading of all generators of a thin knot $K$ is equal to $\Sigma/2$, where $\Sigma$ is defined as in Section \ref{coldif}. We conjecture that every thin knot $K$ is also $S^r$-thin with the $\delta$-grading of all generators of $\CH^{S^r}(K)$ being equal to $r\Sigma/2$. We also conjecture that the refined exponential growth conjecture (\ref{reduced refined exp growth for S}) holds for all thin knots.

The action of the ``special'' differential $d_{1|r}$ on $S^r$-colored homology  is related to the action of the differential $d_{1|1}$ on the uncolored homology. The latter one is denoted $d_0$ in \cite{DGR} and it is conjectured that the homology $H^*(\CH^{\tableau 1}(K),d_0)$ is isomorphic to the knot Floer homology of $K$.
We conjecture the following isomorphism:
\be\label{HFKgrowth}
(\tilde{\CH}^{S^r}(K),d_{1|r})(a,Q=q,t_r=t,t_c=1)=\left(\CH^{\tableau 1}(K),d_{1|1})(a,q,t)\right)^r.
\ee
We note two things: first, the above property does not descend to any relation on the colored HOMFLY polynomial level, due to the non-trivial presence of the $Q$-grading. Secondly, both $d_{1|1}$ on $\CH^{\tableau 1}(K)$ and $d_{1|r}$ on $\CH^{S^r}(K)$
have non-zero $\delta$-grading and therefore act trivially on thin knots.  Hence, the relation (\ref{HFKgrowth}) in the case of thin knots becomes simply (\ref{reduced refined exp growth for S}).
However, once we pass to thick knots, in general the two relations are different, as we shall show in the example below for the (3,4)-torus knot.

\subsection{Rectangular vs. non-rectangular representations}\label{hooks}

What makes rectangular representations special?
Well, these are precisely representations, whose corresponding knot invariants obey Hirota equation.
Specifically, if $\lambda$ is a Young diagram of size $R \times S$, then the corresponding invariants $P^{\lambda} \equiv P_{R,S}$
obey recursion equations of the form
\be
P_{R,S}^2 \; = \; P_{R+1,S} \, P_{R-1,S} + P_{R,S+1} \, P_{R,S-1}
\label{Hirotaeq}
\ee
or
\be
P_{R,S} (u+1) \, P_{R,S} (u-1) \; = \; P_{R+1,S} (u) \, P_{R-1,S} (u) + P_{R,S+1} (u) \, P_{R,S-1} (u)
\label{Hirotaequ}
\ee
where $u$ is the spectral parameter. Both of these equations are variants of the familiar Hirota equation.
For instance, it is easy to verify that
\be
P_{R,S} (\unknot) \; = \; \prod_{i=1}^R \prod_{j=1}^S \frac{a q^{j-i} - a^{-1} q^{-j+i}}{q^{R+S-i-j+1} - q^{-(R+S-i-j+1)}}
\ee
satisfies \eqref{Hirotaeq}.

Another special property of the rectangular Young diagram, is that these are precisely the Young diagrams $\lambda$ such that for every Young diagram $\mu$, the Littlewood-Richardson coefficient $c_{\lambda,\lambda}^{\mu}$ is either 0 or 1. These in turn are precisely those $\lambda$ for which the two-variable HOMFLY polynomial $P^{\lambda}(K)(a,q)$ does not detect mutants (see \cite{MC}).

\section{Examples}
\label{sec:examples}

In this section we give some examples of the quadruply-graded colored HOMFLY homology that exhibit all of the structure and satisfy all of the conjectures from the previous sections. This gives a highly non-trivial check of the existence and consistency of such large and beautiful structure in the colored HOMFLY homology.

\subsection{Trefoil}
\label{sec:trefoil}

We start with the trefoil and its $S^2$ and $\Lambda^2$ colored homologies.
The uncolored HOMFLY homology of the trefoil is (\cite{DGR}):
\be\label{fund31}
\CP^{\tableau1}(3_1)(a,q,t)=a^2q^{-2}+a^2q^2t^2+a^4t^3.
\ee
The uncolored homology is the only one where two $t$-gradings coincide, i.e.
\[
\CP^{\tableau1}(3_1)(a,q,t_r,t_c)=\CP^{\tableau1}(3_1)(a,q,t=t_rt_c).
\]
From (\ref{fund31}) we can get the value of the $S$-invariant
of the trefoil (that we denote $\Sigma$ in this paper):
\be
\Sigma(3_1)=2.
\ee
For easier visualization, triply-graded homologies are usually presented in a diagram form. A generator of $(a,q,t)$-degree $(i,j,k)$
is presented as a dot in the $(q,a)$-plane at the position $(j,i)$ and it
is labeled by its $t$-degree $k$. The triply-graded homology of the
trefoil is presented in a diagram on the left of the Figure \ref{fig:colored3141} below.

\begin{figure}[htb]
\centering
\includegraphics[width=4.5in]{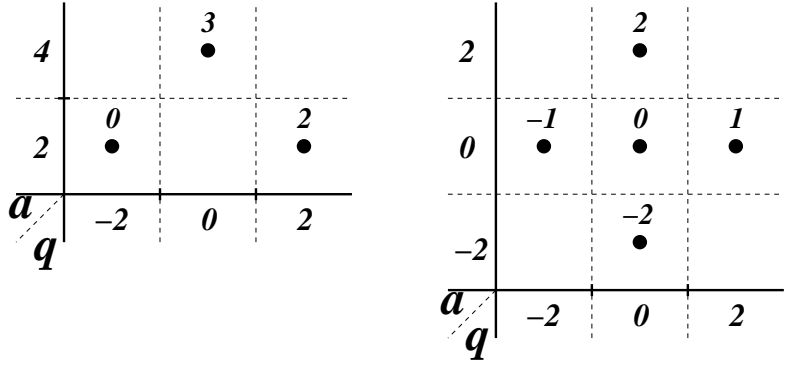}
\caption{The uncolored homologies of the trefoil (left diagram) and figure-eight knot (right diagram).}
\label{fig:colored3141}
\end{figure}

The (reduced) $\tableau2$-colored HOMFLY polynomial of the trefoil knot is equal to:
\[
P^{\tableau2}(3_1)(a,q) = a^{4}q^{-4} + a^4q^2 + a^4q^4 + a^4q^8 - a^6 - a^6q^2 - a^6 q^6 - a^6 q^8 + a^8q^6 \,.
\]
The value of the four-variable $\tableau2$-colored superpolynomial (the Poincar\'e polynomial of the quadruply-graded $\CH^{\tableau2}(3_1)$) of the trefoil is given by
\be\label{quad31}
\CP^{\tableau2}(3_1)(a,q,t_r,t_c)= a^4(q^{-4} +q^2t_r^2t_c^4+q^4t_r^2t_c^6 +q^8t_r^4t_c^8) + $$ $$ a^6 (t_r^3t_c^5+q^2t_r^3t_c^7 + q^6t_r^5t_c^9+q^8t_r^5t_c^{11}) +a^8q^6t_r^6t_c^{12}.
\ee
The triply-graded $S^2$-homology of the trefoil that is computed in \cite{GS} is obtained from (\ref{quad31}) by setting $t_c=1$, whereas the homology obtained in \cite{AShakirov} or \cite{DMMSS} is obtained after setting $t_r=1$.

\begin{figure}[htb]
\centering
\includegraphics[width=4.5in]{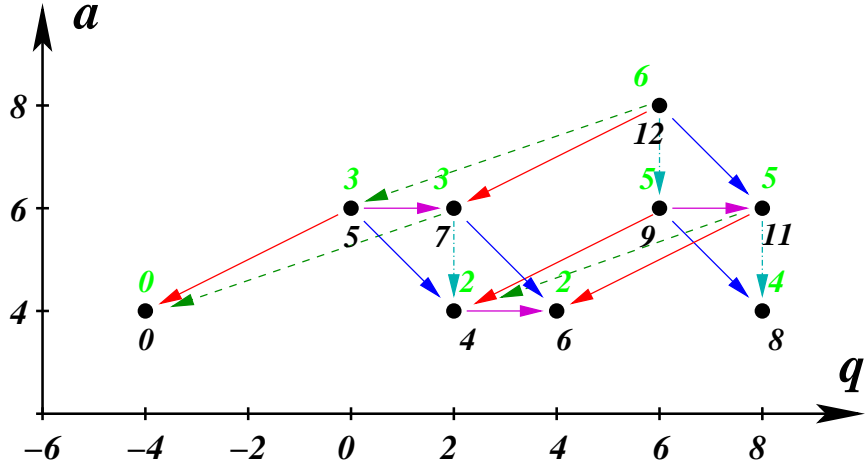}
\caption{The $S^2$-colored homology of the trefoil knot.}
\label{fig:colored31}
\end{figure}

For the second-symmetric representation the $Q$-grading is expressed as $Q=q+t_r-t_c$. Then the tilde-version, i.e. the quadruply-graded $\tableau2$-homology in $(a,Q,t_r,t_c)$-gradings
is given by:
\be\label{tilquad31}
\tilde{\CP}^{\tableau2}(3_1)(a,Q,t_r,t_c)= a^4(Q^{-4} +t_r^2t_c^4+t_r^2t_c^6 +Q^4t_r^4t_c^8) +$$ $$ a^6 (Q^{-2}t_r^3t_c^5+Q^{-2}t_r^3t_c^7 + Q^2t_r^5t_c^9+Q^2t_r^5t_c^{11}) +a^8t_r^6t_c^{12}.
\ee
We note that all generators have the same $\delta$-grading equal to  $+2$. Moreover, the last homology is self-symmetric: in terms of the Poincar\'e polynomial this amounts to:
\[
\tilde{\CP}^{\tableau2}(3_1)(a,Q^{-1}t_r^{-1}t_c^{-2},t_r,t_c)=\tilde{\CP}^{\tableau2}(3_1)(a,Q,t_r,t_c).
\]
The quadruply-graded homologies we represent in the form of a diagram in the following way: for each generator we place a dot in the $(q,a)$-plane. If $(q,a)$-degree of a generator $x$ is $(i,j)$ then
we put a dot at a point with coordinates $(i,j)$. The two $t$-degrees are drawn as labels near the dot: the $t_r$-grading is written above the dot (and colored light green), while the $t_c$-grading of a generator is written below the corresponding dot. In such a way, the homology $\CH^{\tableau2}(3_1)$ is presented in Figure \ref{fig:colored31}.

We can see all five non-trivial colored differentials on this diagram:
There are two canceling ones $d^{col}_{1|0}$ and  $d^{col}_{0|2}$, represented by blue and red arrows, respectively, both leaving a single generator.  Furthermore the two $sl(n|m)$ colored differentials $d^{col}_{1|1}$ (dash-and-dot arrows) and $d^{col}_{0|3}$ (dashed arrows), both leaving the homology isomorphic to $\CH^{\tableau1}(3_1)$. And finally the universal one $d^{\leftarrow}$ (magenta arrows), also leaving the homology isomorphic to $\CH^{\tableau1}(3_1)$.
We note that the re-gradings are exactly as predicted in Sections \ref{coldif} and \ref{unicolor}.

\begin{figure}[htb]
\centering
\includegraphics[width=4.5in]{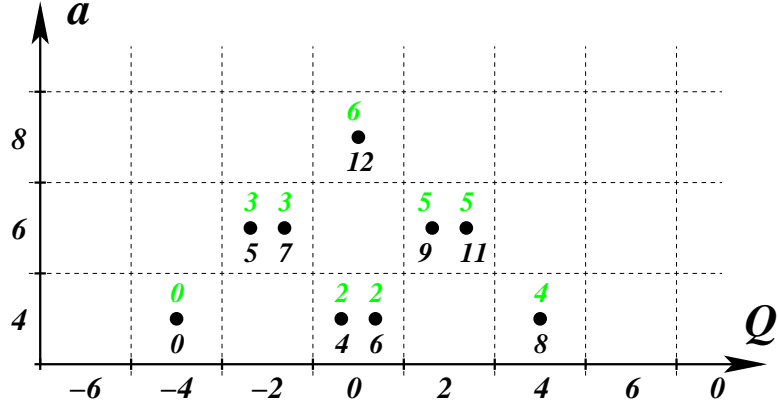}
\caption{The tilde-version of the $S^2$-colored homology of the trefoil knot  ($Q=q+t_r-t_c$).}
\label{fig:colored315}
\end{figure}

\begin{figure}[htb]
\centering
\includegraphics[width=4.5in]{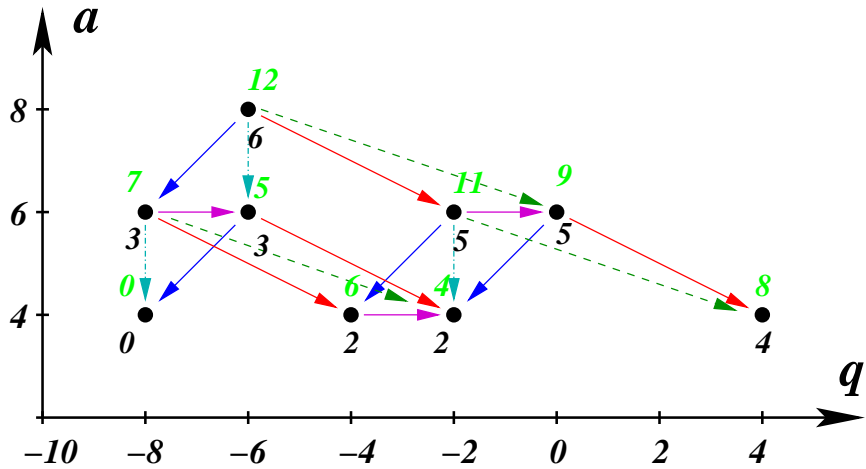}
\caption{The $\Lambda^2$-colored homology of the trefoil knot.}
\label{fig:colored31A2}
\end{figure}

The diagram for the tilde-version of the homology $\CH^{\tableau2}(3_1)$ is shown in Figure \ref{fig:colored315}.
The refined exponential growth is indeed satisfied:
\[
\tilde{\CP}^{\tableau2}(3_1)(a,Q=q,t_r=t,t_c=1)=\left(\CP^{\tableau1}(3_1)(a,q,t)\right)^2.
\]
As for the second antisymmetric representation, we have that
the $\Lambda^2$-colored HOMFLY polynomial of the trefoil is given by
\be
P^{\tableau{1 1}} (3_1)(a,q) \; = \;
a^4 (q^{-8} + q^{-4} + q^{-2} + q^4)
- a^6 (q^{-8} + q^{-6} + q^{-2} + 1)
+ a^8 q^{-6}
\ee
The Poincar\'e polynomial of the quadruply-graded homology $\CH^{\tableau{1 1}}(3_1)$ equals
\be
\CP^{\tableau{1 1}} (3_1) (a,q,t_r,t_c)\; = \;
a^4 (q^{-8} + q^{-4} t_r^6t_c^2 + q^{-2} t_r^4t_c^2 + q^4 t_r^8t_c^4)
+$$ $$ a^6 (q^{-8} t_r^7t_c^3 + q^{-6} t_r^5 t_c^3 + q^{-2} t_r^{11}t_c^5 + t_r^9 t_c^5)
+ a^8 q^{-6} t_r^{12}t_c^6.
\ee
The corresponding diagram, together with the colored differentials, is presented in Figure \ref{fig:colored31A2}.
By comparing Figures \ref{fig:colored31} and \ref{fig:colored31A2}
one can immediately see the mirror symmetry \eqref{mirsymnew}. This symmetry sends $q\mapsto -q$, and commutes with the action of the colored differentials, as can be seen from the diagrams. We remark that all re-gradings from the isomorphisms induced by colored differentials are exactly as predicted in Sections \ref{coldif} and \ref{unicolor}.

\begin{figure}[htb]
\centering
\includegraphics[width=4.5in]{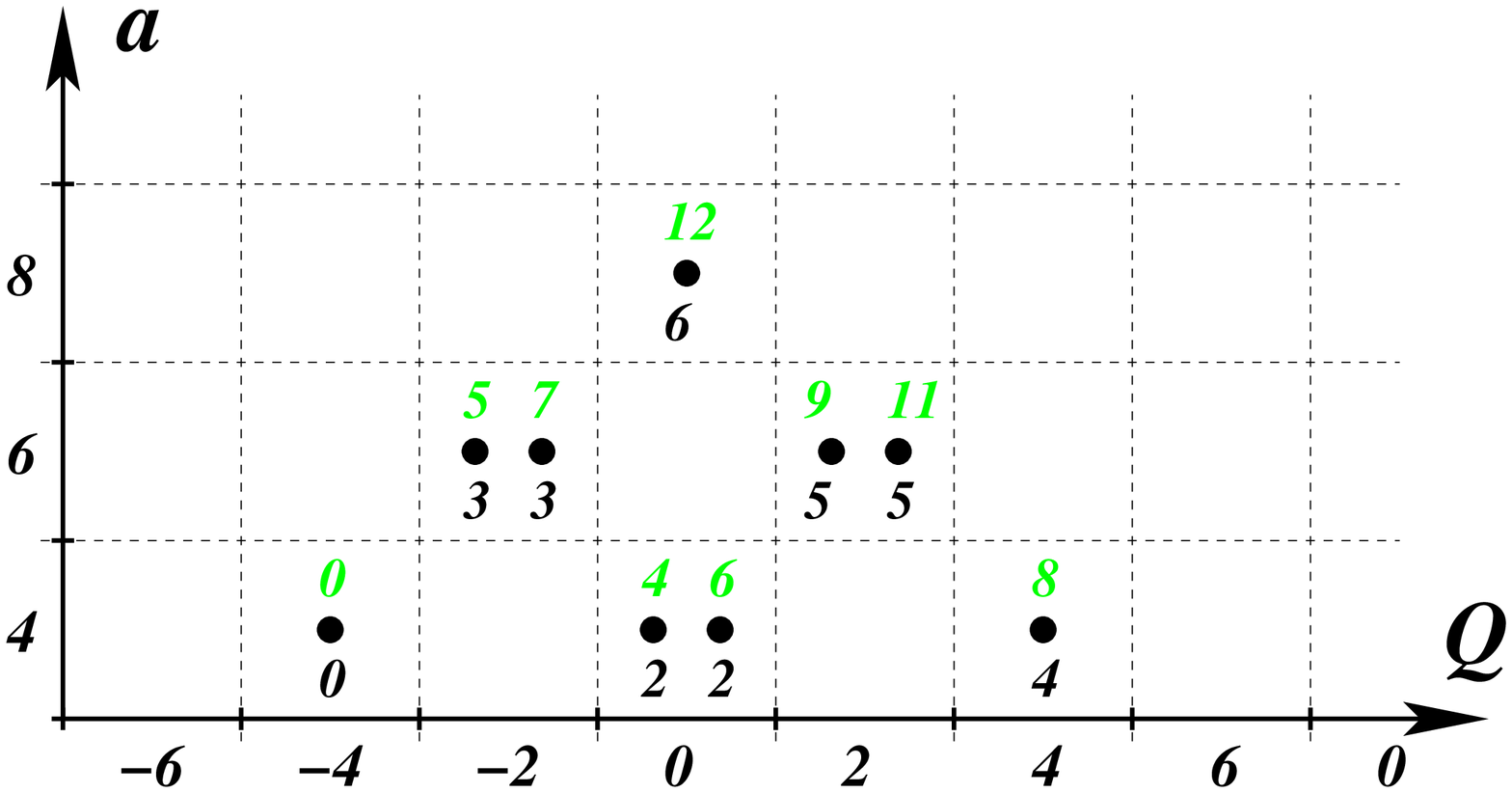}
\caption{The tilde-version of the $\Lambda^2$-colored homology of the trefoil knot  ($Q=\frac{q+t_r-t_c}{2}$).}
\label{fig:colored315A2}
\end{figure}

Mirror symmetry is even better seen on the tilde-version of the homology. For the second-antisymmetric representation $\tableau{1 1}$, the $Q$-grading is obtained by $Q=(q+t_r-t_c)/2$, and so we get:
\be
\tilde{\CP}^{\tableau{1 1}} (3_1) (a,Q,t_r,t_c)\; = \;
a^4 (Q^{-4} +  t_r^6t_c^2 + t_r^4t_c^2 + Q^4 t_r^8t_c^4)
+ $$ $$a^6 (Q^{-2} t_r^7t_c^3 + Q^{-2} t_r^5 t_c^3 + Q^{2} t_r^{11}t_c^5 + Q^2t_r^9 t_c^5)
+ a^8  t_r^{12}t_c^6.
\ee
Its diagram is given in Figure \ref{fig:colored315A2}.
Comparison of Figures \ref{fig:colored315} and \ref{fig:colored315A2} shows that the two are related by a simple exchange of $t_r$ and $t_c$ gradings, as predicted by mirror symmetry (\ref{mirsym}). Also, one can easily check that the explicit
form of the other mirror symmetry (\ref{mirsymnew}) is satisfied, too.


Before passing to further examples, we note that with the explicit mirror symmetries, whenever we know the quadruply-graded $\lambda$-colored homology of a knot $K$, then we automatically have $\lambda^t$-colored homology of $K$, as well:
\[
\tilde{\CP}^{\lambda^t}(K)(a,Q,t_r,t_c)=\tilde{\CP}^{\lambda}(K)(a,Q,t_c,t_r).
\]
Therefore, in the remaining examples we present the results only for one diagram of each pair, $\lambda$ and $\lambda^t$.

\subsection{Figure-eight knot}
We emphasize that our conjectures and predicted structures are for all knots, and not just torus knots. In particular, all homologies computed in \cite{GS} (that includes $S^2$-colored homologies of all prime knots with up to 6 crossings) can be extended to quadruply-graded homologies which have all of the wanted properties. Here we present the explicit expressions for the $S^2$-colored homology of the figure-eight knot.
The uncolored homology of the figure-eight knot is given by:
\[
\CP^{\tableau1}(4_1)(\a,q,t)=\a^2t^2+(q^{-2}t^{-1} + 1 + q^2 t) + \a^{-2}t^{-2},
\]
and its is $S$-invariant is $\Sigma(4_1)=0$.  The triply-graded homology of the figure-eight knot is presented in a diagram on the right of the Figure \ref{fig:colored3141}.
The $\tableau2$-colored HOMFLY polynomial of the figure-eight knot equals to:
$$
P^{\tableau2}(4_1)(\a,q)=\a^4q^4 +\a^2 (-q^{-2} +q^2-  q^4- q^6) +(q^{-6}-q^{-4}+3-q^4+q^6)+$$ $$
\a^{-2}(-q^{-6}-q^{-4}+q^{-2}-q^2) +\a^{-4}q^{-4}.
$$
The quadruply-graded $\tableau2$-colored HOMFLY homology of the figure-eight knot that we computed is given by:
$$
\CP^{\tableau2}(4_1)(\a,q,t_r,t_c)= \a^4q^4t_r^4t_c^8+\a^2 (q^{-2}t_rt_c+t_rt_c^3 +t_r^2t_c^2+ q^2t_r^2t_c^4 +  q^4t_r^3t_c^5+ q^6t_r^3t_c^7) +$$ $$
+(q^{-6}t_r^{-2}t_c^{-4} +q^{-4}t_r^{-1}t_c^{-3}+q^{-2}t_r^{-1}t_c^{-1} + q^{-2}t_c^{-2}+3+q^2t_c^2 + q^2t_rt_c+q^4t_rt_c^3+q^6t_r^2t_c^4)+$$ $$
+\a^{-2}(q^{-6}t_r^{-3}t_c^{-7}+q^{-4}t_r^{-3}t_c^{-5}+q^{-2}t_r^{-2}t_c^{-4}+t_r^{-2}t_c^{-2}+t_r^{-1}t_c^{-3}+q^2t_r^{-1}t_c^{-1}) + \a^{-4}q^{-4}t_r^{-4}t_c^{-8}.
$$

\begin{figure}[htb]
\centering
\includegraphics[width=4.5in]{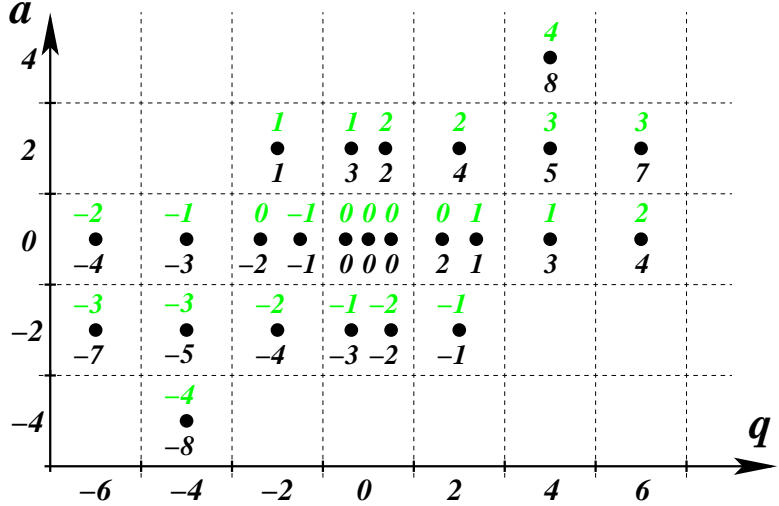}
\caption{The $S^2$-colored homology of the figure-eight knot.}
\label{fig:colored41S2}
\end{figure}

It is presented in a diagram form in Figure \ref{fig:colored41S2}.
All generators have the same $\delta$-degree 0, in accordance with the fact that $4_1$ is a thin knot.
The homology $\CH^{\tableau2}(4_1)$ has five colored differentials: two canceling ones $d^{col}_{1|0}$ and $d^{col}_{0|2}$, further two colored $sl(m|n)$ differentials $d^{col}_{1|1}$ and $d^{col}_{0|3}$ both leaving homology isomorphic to $\CH^{\tableau1}(4_1)$ and the universal one $d^{\leftarrow}$ also leaving  homology isomorphic to $\CH^{\tableau1}(4_1)$. All re-gradings are as given in Sections \ref{coldif} and \ref{unicolor}.

\begin{figure}[htb]
\centering
\includegraphics[width=4.5in]{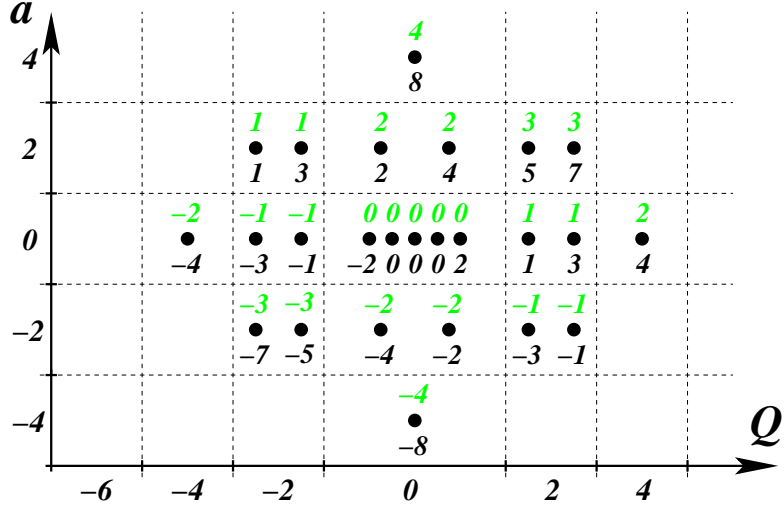}
\caption{The tilde-version of $S^2$-colored homology of the figure-eight knot.}
\label{fig:colored41t2}
\end{figure}

The $Q$-grading is given by $Q=q+t_r-t_c$.
So, the tilde-version of the $\tableau2$-colored HOMFLY homology of the figure-eight knot becomes:
$$
\tilde{\CP}^{\tableau2}(4_1)(\a,Q,t_r,t_c)=\a^4t_r^4t_c^8+\a^2 (Q^{-2}t_rt_c+Q^{-2}t_rt_c^3 +t_r^2t_c^2+ t_r^2t_c^4 +  Q^2t_r^3t_c^5+ Q^2t_r^3t_c^7) +$$ $$
(Q^{-4}t_r^{-2}t_c^{-4} +Q^{-2}t_r^{-1}t_c^{-3}+Q^{-2}t_r^{-1}t_c^{-1} + t_c^{-2}+3+t_c^2 + Q^2t_rt_c+Q^2t_rt_c^3+Q^4t_r^2t_c^4)+$$ $$
\a^{-2}(Q^{-2}t_r^{-3}t_c^{-7}+Q^{-2}t_r^{-3}t_c^{-5}+t_r^{-2}t_c^{-4}+t_r^{-2}t_c^{-2}+Q^2t_r^{-1}t_c^{-3}+Q^2t_r^{-1}t_c^{-1}) +\a^{-4}t_r^{-4}t_c^{-8}.
$$
The tilde-version is presented in Figure \ref{fig:colored41t2} in $(a,Q,t_r,t_c)$-gradings.
It indeed satisfies self-symmetry (\ref{selfsym}) and refined exponential growth (\ref{reduced refined exp growth for S}).

\subsection{(3,4)-torus knot}
In this section we describe the quadruply-graded $S^2$ colored homology of the (3,4)-torus knot $T_{3,4}$ (i.e. $8_{19}$ knot) and show that it satisfies all of the desired properties. This knot is a thick knot, and has nontrivial differential $d_{1|1}$ (i.e. $d_0$) on the uncolored homology, which makes this a highly-nontrivial consistency check of all predicted properties.
The uncolored homology of $T_{3,4}$ is given by \cite{DGR}:
$$
\CP^{\tableau1}(T_{3,4})(\a,q,t)=\a^6q^{-6}+\a^8q^{-4}t^3+\a^6q^{-2}t^2+\a^8t^5+\a^6q^2t^4+\a^8q^{4}t^7+\a^6q^{6}t^6+ $$ $$
+\a^{10}t^8+\a^8q^{-2}t^5+\a^8q^{2}t^7+\a^6t^4.
$$
Its $S$-invariant is $\Sigma(T_{3,4})=6$.
There is a non-trivial differential $d_0$ (or $d_{1|1}$) on $\CH^{\tableau1}(T_{3,4})$ of $(a,q,t)$-degree $(-2,0,-3)$ leaving 5-dimensional homology isomorphic to the knot Floer homology of $T_{3,4}$:

\be
(\CH^{\tableau1}(T_{3,4}),d_{1|1})=\a^6q^{-6}+\a^8q^{-4}t^3+\a^6t^4+\a^8q^{4}t^7+\a^6q^{6}t^6.
\ee
The quadruply-graded $\tableau2$-colored homology of $T_{3,4}$ that we have computed has 121 generators. Its Poincar\'e polynomial is given by:
\begin{multline*}
\tilde{\CP}^{\tableau2}(T_{3,4})(\a,Q,t_r,t_c)=a^{12}(Q^{-12}+Q^{-8}t_r^2(t_c^4+t_c^6)+Q^{-6}t_r^4(t_c^8+t_c^{10})+\\ Q^{-4}t_r^4(t_c^8+t_c^{10}+t_c^{12})+
Q^{-2}t_r^6(t_c^{12}+t_c^{14})+t_r^6(t_c^{12}+t_c^{14}+t_c^{16}+t_c^{18})+t_r^8t_c^{16}+Q^2t_r^8(t_c^{16}+t_c^{18})+\\ Q^4t_r^8(t_c^{16}+t_c^{18}+t_c^{20})+Q^6t_r^{10}(t_c^{20}+t_c^{22})+Q^8t_r^{10}(t_c^{20}+t_c^{22})+Q^{12}t_r^{12}t_c^{24})+\\
\a^{14}(Q^{-10}t_r^{3}(t_c^{5}+t_c^{7})+Q^{-8}t_r^{5}(t_c^{9}+t_c^{11})+Q^{-6}t_r^{5}(t_c^{9}+2t_c^{11}+t_c^{13})+\\ Q^{-4}t_r^{7}(2t_c^{13}+3t_c^{15}+t_c^{17})
+Q^{-2}t_r^{7}(t_c^{13}+2t_c^{15}+2t_c^{17}+t_c^{19})+\\
Q^{-2}t_r^{9}(t_c^{17}+t_c^{19})+t_r^9(2t_c^{17}+3t_c^{19}+t_c^{21})
+ Q^2t_r^{9}(t_c^{17}+2t_c^{19}+2t_c^{21}+t_c^{23})+
 Q^2t_r^{11}(t_c^{21}+t_c^{23})+\\
 Q^4t_r^{11}(2t_c^{21}+3t_c^{23}+t_c^{25})+
Q^6t_r^{11}(t_c^{21}+2t_c^{23}+t_c^{25})+Q^8t_r^{13}(t_c^{25}+t_c^{27})+Q^{10}t_r^{13}(t_c^{25}+t_c^{27}))+\\
+\a^{16}(Q^{-8}t_r^{6}t_c^{12}+Q^{-6}t_r^{8}(t_c^{14}+2t_c^{16}+t_c^{18})+Q^{-4}t_r^{8}(t_c^{16}+t_c^{18})+Q^{-4}t_r^{10}t_c^{20}+\\
Q^{-2}t_r^{10}(t_c^{18}+3t_c^{20}+2t_c^{22})+t_r^{10}(t_c^{20}+t_c^{22}+t_c^{24})+ t_r^{12}(t_c^{22}+2t_c^{24}+t_c^{26}) + \\
 Q^2t_r^{12}(t_c^{22}+3t_c^{24}+2t_c^{26}) + Q^4t_r^{12}(t_c^{24}+t_c^{26}) +Q^4t_r^{14}t_c^{28}+q^6t_r^{14}(t_c^{26}+2t_c^{28}+t_c^{30})+q^8t_r^{14}t_c^{28})+\\
+\a^{18}(Q^{-4}t_r^{11}(t_c^{21}+t_c^{23})\!+\!Q^{-2}t_r^{13}(t_c^{25}+t_c^{27})+t_r^{13}(t_c^{25}+t_c^{27})+Q^{2}t_r^{15}(t_c^{29}+t_c^{31})+\\
+Q^{4}t_r^{15}(t_c^{29}+t_c^{31}))
+\a^{20}t_r^{16}t_c^{32}.
\end{multline*}
The Poincar\'e polynomial in $(a,q,t_r,t_c)$-gradings can be obtained from the above expression by:
\[
\CP^{\tableau2}(T_{3,4})(a,q,t_r,t_c)=\tilde{\CP}^{\tableau2}(T_{3,4})(a,Q=q,t_r=t_r q^{-1},t_c=t_c q).
\]
This homology satisfies all of the expected properties: there are four colored $sl(m|n)$ differentials $d^{col}_{1|0}$, $d^{col}_{0|2}$, $d^{col}_{1|1}$ and $d^{col}_{0|3}$, with the properties and re-gradings exactly matching the predicted behaviour from Section \ref{coldif}. Moreover, there exists a universal colored differential $d^{\leftarrow}$ as predicted in Section \ref{unicolor}.
Furthermore, the homology $\tilde{\CH}^{\tableau2}(T_{3,4})$ satisfies the self-symmetry property:
\[
\tilde{\CP}^{\tableau2}(T_{3,4})(a,Q^{-1}t_r^{-1}t_c^{-2},t_r,t_c)=\tilde{\CP}^{\tableau2}(T_{3,4})(a,Q,t_r,t_c).
\]
It also satisfies the refined exponential growth property:
\[
\tilde{\CP}^{\tableau2}(T_{3,4})(a,Q=q,t_r=t,t_c=1)=\left(\CP^{\tableau1}(T_{3,4})(a,q,t)\right)^2.
\]
Moreover, since HFK differential $d_0$ on $\CH^{\tableau1}(T_{3,4})$ is nontrivial, this time the HFK-like exponential growth property (\ref{HFKgrowth}) is giving a new relation on the homologies. Indeed, there exists a differential $d_{1|2}$ on $\tilde{\CH}^{\tableau2}(T_{3,4})$ of $(a,Q,t_r,t_c)$-degree $(-2,0,-3,-5)$, so that the homology $H^*(\tilde{\CH}^{\tableau2}(T_{3,4}),d_{1|2})$ is 25-dimensional with the following
Poincar\'e polynomial:
\begin{multline*}
(\tilde{\CH}^{\tableau2}(T_{3,4}),d_{1|2})(\a,Q,t_r,t_c)=\a^{12}(Q^{-12}+Q^{-6}t_r^4(t_c^8+t_c^{10})+t_r^6(t_c^{16}+ t_c^{18})+t_r^8 t_c^{16}+\\ Q^{6} t_r^{10} (t_c^{20}+t_c^{22}) + Q^{12} t_r^{12} t_c^{24})
+ \a^{14} (Q^{-10} t_r^3 (t_c^5+t_c^7) + Q^{-4} t_r^7 (t_c^{13}+t_c^{15}) + \\ Q^{-2} t_r^7 (t_c^{17}+t_c^{19})+Q^{2} t_r^9 (t_c^{21}+t_c^{23})  +Q^{4} t_r^{11} (t_c^{21}+t_c^{23}) + Q^{10} t_r^{13} (t_c^{25}+t_c^{27}) ) +\\
+ \a^{16} (Q^{-8} t_r^6 t_c^{12} + t_r^{10} (t_c^{22} + t_c^{24} ) + Q^8 t_r^{14} t_c^{28}).
\end{multline*}
Therefore
\be\label{HFKgrowtht34}
(\tilde{\CH}^{\tableau2}(T_{3,4}),d_{1|2})(a,Q=q,t_r=t,t_c=1)=\left(\CH^{\tableau 1}(T_{3,4}),d_{1|1})(a,q,t)\right)^2.
\ee

\subsection{$(2,2)$ representation}
\label{sec:22}
In this section we describe the $\tableau{2 2}$-colored HOMFLY homology of the trefoil knot.
It should categorify the reduced $\tableau{2 2}$-colored HOMFLY polynomial of the trefoil:
\begin{multline*}
P^{\tableau{2 2}} (3_1)(\a,q) =
\a^{16}
 + \a^{14} (-q^{-6} - 2q^{-4} - q^{-2} - q^2 - 2 q^4 - q^6)\\
+ \a^{12} (4 + 2q^{-10} + 2q^{-8} + 2 q^{-6} + q^{-4} + 3 q^{-2} + 3 q^2 + q^4 + 2 q^6 + 2 q^8 + 2 q^{10})\\
+ \a^{10} (-2 - q^{-14} - 2 q^{-12} - q^{-10} - q^{-8} - 3 q^{-6}\\ - 4 q^{-4} - 3 q^{-2} - 3 q^2 - 4 q^4 - 3 q^6 - q^8 - q^{10} - 2 q^{12} - q^{14})\\
+ \a^{8} (2 + q^{-16} + q^{-10} + 2 q^{-8} + q^{-6} + q^{-4} + q^{-2} + q^2 + q^4 + q^6 + 2 q^8 + q^{10} + q^{16})
\end{multline*}
 The  $\tableau{2 2}$-colored HOMFLY homology of the trefoil knot that we computed, has 81 generators. Its four-variable Poincar\'e polynomial is given by:
\begin{multline*}
\CP^{\tableau{2 2}} (3_1) (\a,q,t_r,t_c) =
\a^{16} t_r^{24}  t_c^{24}+ \a^{14} (q^{-6}{t_r^{19}t_c^{17}} + q^{-4}({t_r^{17}t_c^{17}} +{t_r^{19}t_c^{19}})+\\ q^{-2}{t_r^{17}t_c^{19}}+ q^2 t_r^{23}t_c^{21} + q^4 (t_r^{21}t_c^{21} + t_r^{23}t_c^{23}) + q^6 t_r^{21}t_c^{23})+  \\
 + \a^{12} (q^{-10}({t_r^{12}t_c^{10}} + {t_r^{14}t_c^{12}}) + 2q^{-8}{t_r^{12} t_c^{12}} + q^{-6}(t_r^{10}{t_c^{12}} +  t_r^{12}{t_c^{14}}) +q^{-4}t_r^{18}{t_c^{14}} +\\
+ q^{-2}({t_r^{16}t_c^{14}}  + 2t_r^{14}{t_c^{16}})+ (3t_r^{16} t_c^{16}+t_r^{18}t_c^{18}) + q^2 (t_r^{14}t_c^{16} + 2  t_r^{16}t_c^{18}) + \\
+ q^4 t_r^{14}t_c^{18} + q^6 (t_r^{20}t_c^{18} + t_r^{22}t_c^{20}) + 2 q^8 t_r^{20}t_c^{20} +  q^{10} (t_r^{18}t_c^{20} + t_r^{20} t_c^{22}))+ \\
 + \a^{10} (q^{-14}t_r^{7}{t_c^5} + q^{-12}({t_r^{5}t_c^5} + t_r^{7}{t_c^7}) + q^{-10}t_r^{5}{t_c^7} + q^{-8}t_r^{13}{t_c^9} +  q^{-6} (2t_r^{11}{t_c^9}  +t_r^{13}{t_c^{11}})+\\
 + q^{-4}(t_r^{9}{t_c^9} + 3t_r^{11}{t_c^{11}}) + q^{-2}(2t_r^{9}{t_c^{11}}+t_r^{11}{t_c^{13}})+ (t_r^{9} t_c^{13} +t_r^{17}t_r^{13})+\\
 +  q^2 (2t_r^{15}t_c^{13} +t_r^{17}t_c^{15}) + q^4 (t_r^{13}t_c^{13} +3t_r^{15}t_c^{15})
+ q^6 (2 t_r^{13}t_c^{15} +  t_r^{15}t_c^{17}) + \\
+ q^8 t_r^{13} t_c^{17} + q^{10} t_r^{19}t_c^{17} + q^{12} (t_r^{17}t_c^{17} +t_r^{19}t_c^{19}) + q^{14} t_r^{17}t_c^{19}) +\\
 + \a^{8} (q^{-16} + q^{-10}t_r^{6}{t_c^4} + q^{-8}(t_r^{4}{t_c^4} + t_r^{6}{t_c^6}) + q^{-6}t_r^{4}{t_c^6}
+ q^{-4}t_r^{12}{t_c^8} + q^{-2}t_r^{10}{t_c^8}+\\
 (t_r^{8}t_c^8 + t_r^{10}t_c^{10}) + q^2 t_r^{8}t_c^{10} + q^4 t_r^{8}t_c^{12}  + q^6 t_r^{14}t_c^{12} + q^8(t_r^{12}t_c^{12} + t_r^{14} t_c^{14}) + q^{10} t_r^{12}t_c^{14} + q^{16} t_r^{16}t_c^{16}).
\end{multline*}
Let us list various differentials and their $(\a,q,t_r,t_c)$-degrees. Everywhere in this section $\lambda=\tableau{2 2}$.

\noindent
1. {Categorification:}
$$
\CP^{\lambda} (3_1) (\a,q,t_r=-1,t_c=1)= P^{\lambda} (3_1)(\a,q) = \CP^{\lambda} (3_1) (\a,q,t_r=1,t_c=-1).
$$
2. { Canceling differentials:} There exist differentials $d^{col}_{2|0}$ and $d^{col}_{0|2}$, with  degrees $(-2,4,-1,-1)$ and $(-2,-4,-5,-5)$, respectively. Remaining generators have degrees $(8,-16,0,0)$ and $(8,16,16,16)$, respectively.

\noindent
3. { Row-removing differentials:} There exist differentials $d^{col}_{3|0}$ and $d^{col}_{1|2}$ of degrees $(-2,6,-3,-1)$ and $(-2,-2,-7,-5)$, respectively, such that
\[
H^*\left( \CH^{\lambda}(3_1), d^{col}_{3|0} \right) \; \cong \; \CH^{\tableau{2}}(3_1),
\]
\[
H^*\left( \CH^{\lambda}(3_1), d^{col}_{1|2} \right) \; \cong \; \CH^{\tableau{2}}(3_1).
\]
4. { Column-removing differentials:} There exist differentials $d^{col}_{0|3}$ and $d^{col}_{2|1}$ of degrees  $(-2,-6,-5,-7)$ and $(-2,2,-1,-3)$, respectively, such that
\[
H^*\left( \CH^{\lambda}(3_1), d^{col}_{0|3} \right) \; \cong \; \CH^{\tableau{1 1}}(3_1),
\]
\[
H^*\left( \CH^{\lambda}(3_1), d^{col}_{2|1} \right) \; \cong \; \CH^{\tableau{1 1}}(3_1).
\]
5. { Universal colored differentials:} There exist differentials $d^{\uparrow}$ and $d^{\leftarrow}$ of  degrees $(0,2,-2,0)$ and $(0,2,0,2)$, respectively, such that


\[
H^*\left( \CH^{\lambda}(3_1), d^{\uparrow} \right) \; \cong \; \CH^{\tableau{2}}(3_1),
\]
\[
H^*\left( \CH^{\lambda}(3_1), d^{\leftarrow} \right) \; \cong \; \CH^{\tableau{1 1}}(3_1).
\]
The remaining properties are best seen in the $(\a,Q,t_r,t_c)$-gradings. Recall that in this case the $Q$-grading can be expressed via the other ones by the equation $$Q=(q+t_r-t_c)/2.$$ We denote the $\lambda$-colored homology in these four gradings by $\tilde{\CH}^{\lambda}(3_1)$ and the corresponding four-variable Poincar\'e polynomial by $\tilde{\CP}^{\lambda}(3_1)(\a,Q,t_r,t_c)$. In terms of $\CP^{\lambda}(3_1)$, we have
\[
\tilde{\CP}^{\lambda}(3_1)(\a,Q,t_r,t_c)={\CP}^{\lambda}(3_1)(\a,Q^{1/2},t_rQ^{1/2},t_cQ^{-1/2}).
\]
6. { Self-symmetry:}
\[
 \tilde{\CH}^{\lambda}_{i,j,k,l}(3_1) \cong \tilde{\CH}^{\lambda}_{i,-j,k-2j,l-2j}(3_1).
\]
7. { Mirror symmetry:}
Since $\tableau{2 2}^t = \tableau{2 2}$, the mirror symmetry in this case is an another, different symmetry on $\tilde{\CH}^{\lambda}(3_1)$:
\[
 \tilde{\CH}^{\lambda}_{i,j,k,l}(3_1) \cong \tilde{\CH}^{\lambda}_{i,j,l,k}(3_1).
\]
8. { Refined exponential growth:}
\[
\tilde{\CP}^{\lambda}(3_1)(\a,Q,t_r,t_c=1)=\left(\tilde{\CP}^{\tableau{1 1}}(3_1)(\a,Q,t_r,t_c=1)\right)^2.
\]


\subsection{$(2,2,2)$ representation}
\label{sec:222}

In this section we give the result for the $\tableau{2 2 2}$-colored homology of the trefoil. Although this homology is very large (729-dimensional) it
contains all colored differentials predicted in Section \ref{coldif}.
Due to its size, and in order to simplify the presentation, we give the homology and the differentials in $(\a,q,t_c)$-gradings, and simply denote $t_c$ by $t$.
In $(\a,q,t_c)$-gradings the Poincar\'e polynomial $\CP^{\lambda}(3_1)$  is given by:

$$
{\scriptstyle{
\CP^{\tableau{2 2 2}}(3_1)=\a^{12}(q^{-36}+q^{-30}t^4+q^{-28}(t^4+t^6)+q^{-26}(t^4+t^6)+q^{-24}(t^6+t^8)+q^{-22}t^8 }} $$ $$
{\scriptstyle{+q^{-20}(2t^8+t^{10})+q^{-18}(t^8+2t^{10})+q^{-16}(t^8+2t^{10}+t^{12})+q^{-14}(t^{10}+2t^{12}) }} $$ $$
{\scriptstyle{+q^{-12}(3t^{12}+t^{14})+q^{-10}(2t^{12}+2t^{14})+q^{-8}(t^{12}+3t^{14})+q^{-6}(t^{12}+2t^{14}+t^{16})+q^{-4}(t^{14}+2t^{16}) }} $$ $$
{\scriptstyle{+2q^{-2}t^{16}+q^{0}(2t^{16}+2t^{18})+q^{2}(t^{16}+2t^{18})+}} $$ $$ {\scriptstyle{q^{4}(t^{16}+2t^{18}+t^{20})+q^{6}(t^{18}+t^{20})+}}
{\scriptstyle{+q^{8}t^{20}+q^{10}t^{20}+q^{12}(t^{20}+t^{22})+q^{14}(t^{20}+t^{22})+q^{16}t^{22}+q^{24}t^{24}) }} $$ $$
{\scriptstyle{+\a^{14}(q^{-36}t^5+q^{-34}(t^5+t^7)+q^{-32}(t^5+t^7)
+q^{-30}(t^7+t^9)+q^{-28}(2t^9+t^{11})+q^{-26}(3t^9+3t^{11}) }} $$ $$
{\scriptstyle{+q^{-24}(2t^9+5t^{11}+t^{13})+q^{-22}(t^9+4t^{11}+3t^{13})+q^{-20}(2t^{11}+5t^{13}+t^{15})+q^{-18}(5t^{13}+4t^{15}) }} $$ $$
{\scriptstyle{+q^{-16}(4t^{13}+7t^{15}+t^{17})+q^{-14}(2t^{13}+8t^{15}+3t^{17})+ q^{-12}(t^{13}+5t^{15}+6t^{17})+q^{-10}(2t^{15}+7t^{17}+2t^{19}) }} $$ $$
{\scriptstyle{+q^{-8}(6t^{17}+5t^{19})+q^{-6}(4t^{17}+8t^{19}+t^{21})+q^{-4}(2t^{17}+8t^{19}+3t^{21})+ q^{-2}(t^{17}+5t^{19}+5t^{21}) }} $$ $$
{\scriptstyle{+q^{0}(2t^{19}+5t^{21}+t^{23})+q^{2}(4t^{21}+2t^{23})+q^{4}(3t^{21}+4t^{23})+q^{6}(2t^{21}+5t^{23}+t^{25}) }} $$ $$
{\scriptstyle{+q^{8}(t^{21}+4t^{23}+2t^{25})+q^{10}(2t^{23}+2t^{25})+q^{12}t^{25}+q^{14}t^{25}+q^{16}(t^{25}+t^{27})+q^{18}(t^{25}+t^{27})+q^{20}t^{27})+}} $$ $$
{\scriptstyle{+\a^{16}(q^{-34}(t^{10}+t^{12})+q^{-32}(t^{10}+2t^{12})
+q^{-30}(t^{10}+3t^{12}+t^{14})+q^{-28}(2t^{12}+2t^{14})}}  $$ $${\scriptstyle{+q^{-26}(t^{12}+3t^{14}+2t^{16})+q^{-24}(3t^{14}+5t^{16}+t^{18})+q^{-22}(2t^{14}+8t^{16}+3t^{18})}}$$ $${\scriptstyle{q^{-20}(t^{14}+7t^{16}+6t^{18})+q^{-18}(4t^{16}+8t^{18}+2t^{20}) }} $$ $$
{\scriptstyle{+q^{-16}(t^{16}+7t^{18}+5t^{20})+q^{-14}(5t^{18}+10t^{20}+2t^{22})+ q^{-12}(2t^{18}+11t^{20}+5t^{22})+}} $$ $${\scriptstyle{q^{-10}(t^{18}+9t^{20}+9t^{22}+t^{24})+
q^{-8}(4t^{20}+9t^{22}+2t^{24})+q^{-6}(t^{20}+7t^{22}+5t^{24})}}  $$ $${\scriptstyle{q^{-4}(4t^{22}+7t^{24}+t^{26})+ q^{-2}(2t^{22}+9t^{24}+3t^{26})+}} $$ $$
{\scriptstyle{\quad+q^{0}(t^{22}+7t^{24}+6t^{26})+q^{2}(4t^{24}+6t^{26}+t^{28})+q^{4}(t^{24}+4t^{26}+t^{28})+q^{6}(2t^{26}+2t^{28})+}} $$ $$
{\scriptstyle{\quad+q^{8}(t^{26}+2t^{28})+q^{10}(t^{26}+3t^{28}+t^{30})+q^{12}(2t^{28}+t^{30})+q^{14}(t^{28}+t^{30})+}} $$ $$
{\scriptstyle{\!\!\!\!+\a^{18}(q^{-32}t^{17}
+q^{-30}(t^{15}+2t^{17}+t^{19})+q^{-28}(3t^{17}+2t^{19})+q^{-26}(2t^{17}+3t^{19})+q^{-24}(t^{17}+3t^{19}+2t^{21})+}} $$ $$
{\scriptstyle{\quad+q^{-22}(2t^{19}+4t^{21}+t^{23})+q^{-20}(t^{19}+6t^{21}+4t^{23})+q^{-18}(6t^{21}+7t^{23}+t^{25})+q^{-16}(3t^{21}+8t^{23}+2t^{25})+}} $$ $$
{\scriptstyle{\quad+q^{-14}(t^{21}+6t^{23}+4t^{25})+ q^{-12}(3t^{23}+6t^{25}+t^{27})+q^{-10}(t^{23}+7t^{25}+4t^{27})+q^{-8}(6t^{25}+7t^{27}+t^{29})+}} $$ $$
{\scriptstyle{\quad+q^{-6}(3t^{25}+8t^{27}+2t^{29})+q^{-4}(t^{25}+5t^{27}+3t^{29})+ q^{-2}(2t^{27}+3t^{29})+q^{0}(t^{27}+3t^{29}+t^{31})+}} $$ $$
{\scriptstyle{\quad+q^{2}(3t^{29}+2t^{31})+q^{4}(2t^{29}+3t^{31})+q^{6}(t^{29}+2t^{31}+t^{33})+q^{8}t^{31})+}} $$ $$
{\scriptstyle{
\!\!\!\!+\a^{20}(q^{-28}(t^{22}+t^{24})+q^{-26}(t^{22}+2t^{24})+q^{-24}(t^{22}+3t^{24}+t^{26})+q^{-22}(2t^{24}+t^{26})+q^{-20}
(t^{24}+2t^{26})}} $$ $$
{\scriptstyle{\quad+q^{-18}(2t^{26}+2t^{28})+q^{-16}(2t^{26}+4t^{28}+t^{30})+q^{-14}(t^{26}+5t^{28}+2t^
{30})+ q^{-12}(3t^{28}+3t^{30})+}} $$ $$
{\scriptstyle{\quad+q^{-10}(t^{28}+2t^{30})+q^{-8}(2t^{30}+t^{32})+q^{-6}(t^{30}+2t^{32})+q^{-4}(t^{30}+3t^{32}+t^{34})+
q^{-2}(2t^{32}+t^{34})+}} $$ $$
{\scriptstyle{q^{0}(t^{32}+t^{34}))+\a^{22}(q^{-24}t^{29}+q^{-22}(t^{29}+t^{31})+q^{-20}
(t^{29}+t^{31})+q^{-18}t^{31} +q^{-14}t^{33}+ q^{-12}(t^{33}+t^{35})+}} $$ $$
{\scriptstyle{q^{-10}(t^{33}+t^{35})+q^{-8}t^{35}) +}}
{\scriptstyle{\!\!\!\!+ \a^{24}q^{-18}t^{36}.
}}
$$
On such a big homology $\CH^{\lambda}(3_1)$ categorifying the reduced $\lambda$-colored HOMFLY polynomial of the trefoil knot, there exist all of the predicted colored differentials:

\noindent
1. Canceling differential $d^{col}_{3|0}$ of $(\a,q,t_c)$-degree $(-2,6,-1)$ such that the homology $H^*(\CH^{\lambda}(3_1),d^{col}_{3|0})$ is one-dimensional, and the remaining generator has $(\a,q,t_c)$-degree $(12,-36,0)$.

\noindent
2. Canceling differential $d^{col}_{0|2}$ of degree $(-2,-4,-5)$, such that the homology $H^*(\CH^{\lambda}(3_1),d^{col}_{0|2})$ is one-dimensional, and the remaining generator has $(\a,q,t_c)$-degree $(12,24,24)$.

\noindent
3. Positive row-removing differential $d^{col}_{4|0}$ of  degree $(-2,8,-1)$, such that
\[
H^*(\CH^{\lambda}(3_1),d^{col}_{4|0}) \cong \CH^{\tableau2}(3_1).
\]
4. Negative row-removing differential $d^{col}_{1|2}$ of  degree $(-2,-2,-5)$, such that
\[
H^*(\CH^{\lambda}(3_1),d^{col}_{1|2}) \cong \CH^{\tableau2}(3_1).
\]
5. Positive row-removing differential $d^{col}_{5|0}$ of  degree $(-2,10,-1)$, such that
\[
H^*(\CH^{\lambda}(3_1),d^{col}_{5|0}) \cong \CH^{\tableau{2 2}}(3_1).
\]
6.  Negative row-removing differential $d^{col}_{2|2}$ of  degree $(-2,0,-5)$, such that
\[
H^*(\CH^{\lambda}(3_1),d^{col}_{2|2}) \cong \CH^{\tableau{2 2}}(3_1).
\]
7.   Negative column-removing differential $d^{col}_{0|3}$ of  degree $(-2,-6,-7)$, such that
\[
H^*(\CH^{\lambda}(3_1),d^{col}_{0|3}) \cong \CH^{\tableau{1 1 1}}(3_1).
\]
8. Positive column-removing differential $d^{col}_{3|1}$ of  degree $(-2,-4,-3)$, such that
\[
H^*(\CH^{\lambda}(3_1),d^{col}_{3|1}) \cong \CH^{\tableau{1 1 1}}(3_1).
\]
9.  Universal colored differential $d^{\leftarrow}$ of  degree $(0,2,2)$, such that
\[
H^*(\CH^{\lambda}(3_1),d^{\leftarrow}) \cong \CH^{\tableau{1 1 1}}(3_1).
\]

\subsection{$(2,1)$ ``hook'' representation}

The $\tableau{2 1}$-colored HOMFLY polynomial of the trefoil is given by
\begin{align}
P^{\tableau{2 1}} (3_1)(a,q) \; = \;
& -a^{12} + a^{10} (q^{-6} + q^{-4} + q^{-2} + q^2 + q^4 + q^6) \notag \\
& + a^{8} (-q^{-10} - 2 q^{-6} - 3 q^{-2} - 3 q^2 - 2 q^6 - q^{10}) \notag \\
& + a^{6} (q^{-10} + 2 q^{-6} - q^{-4} + 2 q^{-2} + 2 q^2 - q^4 + 2 q^6 + q^{10}). \notag
\end{align}
This polynomial has 31 terms, and therefore any homology theory categorifying it, should have at least 31 generators. This already indicates distinction between rectangular and hook diagrams: this time, already for the trefoil, the exponential growth property  cannot be satisfied (the uncolored homology of the trefoil is 3-dimensional and $31>3^3$). However, some of the structural properties still exist in the colored HOMFLY homology theories for the hook-shaped Young diagrams.

As a $sl(n|m)$ representation, $\lambda = \tableau{2 1}$ has three nice isomorphisms for particular values of $n$ and $m$. Namely we have that $\tableau{2 1}\cong \tableau1$ as $sl(2|0)$, as $sl(1|1)$ and as $sl(0|2)$ representations. Therefore, we expect three colored differentials, $d^{col}_{2|0}$, $d^{col}_{1|1}$ and $d^{col}_{0|2}$ on $\CH^{\tableau{2 1}}(K)$, such that

\begin{eqnarray}
H^*(\CH^{\tableau{2 1}}(K),d^{col}_{2|0})&\cong& \CH^{\tableau{1}}(K)\label{ht1},\\
H^*(\CH^{\tableau{2 1}}(K),d^{col}_{1|1})&\cong&\CH^{\tableau{1}}(K) \label{ht2},\\
H^*(\CH^{\tableau{2 1}}(K),d^{col}_{0|2})&\cong& \CH^{\tableau{1}}(K).\label{ht3}
\end{eqnarray}
The $(a,q,t)$-degrees of these three differentials are:
\be
\deg d^{col}_{2|0} = (-2,4,-1), \quad\quad \deg d^{col}_{1|1} = (-2,0,-3), \quad\quad \deg d^{col}_{0|2} = (-2,-4,-5).
\ee
The requirement that $\CH^{\tableau{2 1}}(3_1)$ should categorify $P^{\tableau{2 1}} (3_1)(a,q)$, together with the existence of the colored differentials $d^{col}_{2|0}$, $d^{col}_{1|1}$ and $d^{col}_{0|2}$
essentially determine the $\tableau{2 1}$-colored homology of the trefoil knot. It has 41 generators and its Poincar\'e polynomial is given by:
\begin{multline}
\label{barebones}
\CP^{\tableau{2 1}} (3_1)(a,q,t)=
a^{12} t^{13} + a^{10} (q^{-6} t^8 + q^{-4} t^8 + q^{-2} t^{10} + t^{10}+t^{11}\\ + q^2 t^{12} + q^4 t^{12} + q^6 t^{14})
+ a^{8} (q^{-10} t^3 + 2 q^{-6} t^5 + q^{-4} t^5 + q^{-4} t^6+3 q^{-2} t^7 + t^7 + t^8 + 3 q^2 t^9 \\ + q^4 t^9 + q^4 t^{10} + 2 q^6 t^{11} + q^{10} t^{13})
+ a^{6} (q^{-10} + 2 q^{-6} t^2 + q^{-4} t^3 + 2 q^{-2} t^4 + t^4+t^5\\+ 2 q^2 t^6 + q^4 t^7 + 2 q^6 t^8 + q^{10} t^{10}).
\end{multline}
Additional nice properties that this homology satisfies, and which also exist in the rectangular case, is the mirror symmetry. Since ${\tableau{2 1}}^t={\tableau{2 1}}$ this becomes a symmetry on $\CH^{{\tableau{2 1}}}(3_1)$:

\be
M_{\tableau{2 1}} : \CH^{\tableau{2 1}}_{i,j,k}(3_1) \to \CH^{\tableau{2 1}}_{i,-j,k-j}(3_1),
\ee
with the property that
\be
M_{\tableau{2 1}} d_{n|m} \; = \; d_{m|n} M_{\tableau{2 1}}.
\ee

There are some properties that we predict in the case of rectangular representations and which do not exist for the hook diagram. First of all, already for the trefoil the exponential growth property is not satisfied. The homology is 41-dimensional, and the ${\tableau{2 1}}$-colored HOMFLY polynomial already has 31 terms.

Second, this time the isomorphisms (\ref{ht1})-(\ref{ht3}) are valid only after
appropriate collapse of tri-grading, i.e. after setting $a=q^{n-m}$ in the case of the differential $d_{n|m}$.
Additionally, there is no explicit re-grading in these isomorphisms, and we could not find the second $t$-grading. Therefore the theory that we propose is triply-graded.

Finally, we note that there exists the $SL(2)$ Lefschetz-like action on homology \eqref{barebones} with raising operator of degree $(0,4,2)$.
According to the ``dictionary'' \eqref{aqtgradings}, this raising operator corresponds to forming a bound state with a BPS state of D0-brane charge $n=4$.


\section{Colored HOMFLY homology for torus knots}
\label{sec:algmodel}

In this part we describe various algebraic and geometric approaches to the colored homology of
torus knots.

\subsection{Stable HOMFLY homology of torus knots}
\label{sec:torus}

We will be interested in the stable limit of the colored HOMFLY polynomials of $(n,m)$ torus knots at $m\to\infty$.
Using e.g. \cite[eq. (5.4)]{stevan}, one can prove the following result.

\begin{proposition}
Let $\kappa(\mu)=\sum_{(i,j)\in \mu}(i-j).$
Define
$$P_{st}^{\lambda}(T({m,n})) := q^{-mn\kappa(\lambda)+\frac{m}{n}\kappa(n\lambda)}a^{-\frac{m(n-1)|\lambda|}{2}}P^{\lambda}(T({m,n})).$$
Then
\begin{equation}
\label{stable poly}
\lim_{m\to\infty}P_{st}^{\lambda}(T({m,n}))=P^{n\lambda}(a,q).
\end{equation}
\end{proposition}
As a corollary, one can compute the stable colored invariants of torus knots using the hook formula (\ref{plambda}).
One can easily check that the unreduced stable polynomial $P^{\lambda}(T(n,\infty))=P^{n\lambda}(\unknot ; a,q)$ is
a series in $(-a)$ and $q$ with nonnegative coefficients.

Let us compute the stable polynomial of the $(2,\infty)$ torus knot colored by the Young diagram $\lambda=(2,1)$.
By (\ref{stable poly}) it is the same as the polynomial of the unknot colored by the Young diagram $2\lambda=(4,2)$.
The hook formula (\ref{plambda}) gives the following answer (up to an overall scaling):
$$P^{2\lambda}(\unknot)=\frac{(1-aq^{-1})(1-a)^2(1-aq)(1-aq^2)(1-aq^{3})}{(1-q)^2(1-q^2)^2(1-q^4)(1-q^5)}.$$
Note that
$$P^{\lambda}(\unknot)=\frac{(1-aq^{-1})(1-a)(1-aq)}{(1-q)^2(1-q^3)}.$$
The reduced polynomial
$P^{red}_{\lambda}(T(2,\infty))=\frac{P^{\lambda}(T(2,\infty))}{P^{\lambda}(a,q)}$ has positive and negative terms
in the $q$-expansion of the coefficient at $a^{0}$:
$$\frac{(1-q^3)}{(1-q^2)^2(1-q^4)(1-q^5)}=1 + 2 q^2 - q^3 + 4 q^4 - q^5 + 6 q^6 - 2 q^7 + 8 q^8 - 2 q^9 +
 11 q^{10}+\dots$$


We would like to use the equation (\ref{stable poly}) as a definition of the stable colored homology of $(n,\infty)$ torus knot.
More precisely, we conjecture that as a vector space $\CH^{\lambda}(n,\infty)$ is isomorphic to $\CH^{n\lambda}(\unknot)$, moreover, the $a,q,t_{c}$ gradings on these spaces match. This suggests that $\CH^{\lambda}(n,\infty)$
is a free supercommutative algebra with bosonic and fermionic generators labelled by the boxes of $n\lambda$.

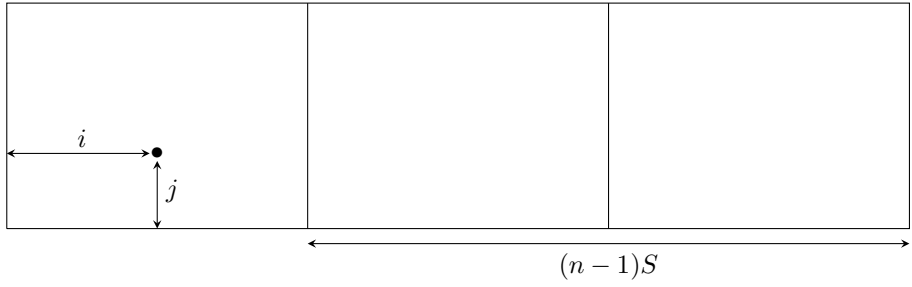
\begin{figure}
\begin{tikzpicture}
\draw (0,0)--(0,3)--(12,3)--(12,0)--(0,0);
\draw (4,0)--(4,3);
\draw (8,0)--(8,3);
\draw (2,1) node {$\bullet$};
\draw [<->,>=stealth] (0,1)--(1.9,1);
\draw [<->,>=stealth] (2,0)--(2,0.9);
\draw [<->,>=stealth] (4,-0.2)--(12,-0.2);
\draw (1,1.2) node {$i$};
\draw (2.2,0.5) node {$j$};
\draw (8,-0.5) node {$(n-1)S$};
\end{tikzpicture}
\caption{Generators $u_{ij}^{(n)}$}
\label{nrect1}
\end{figure}

\begin{figure}
\begin{tikzpicture}
\draw (0,0)--(0,3)--(12,3)--(12,0)--(0,0);
\draw (4,0)--(4,3);
\draw (8,0)--(8,3);
\draw (10,1) node {$\bullet$};
\draw [<->,>=stealth] (8,1)--(9.9,1);
\draw [<->,>=stealth] (10,0)--(10,0.9);
\draw [<->,>=stealth] (0,-0.2)--(8,-0.2);
\draw (9,1.2) node {$i$};
\draw (10.2,0.5) node {$j$};
\draw (4,-0.5) node {$(n-1)S$};
\end{tikzpicture}
\caption{Generators $\xi_{ij}^{(n)}$}
\label{nrect2}
\end{figure}
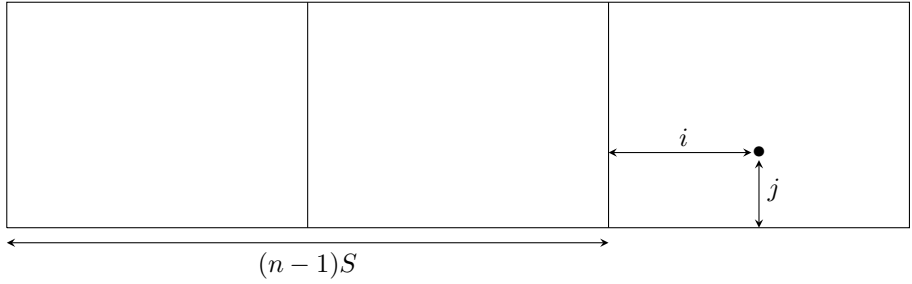

Let $\lambda$ be an $R \times S$ rectangular Young diagram, and let $K=T(p,q)$ be a $(p,q)$-torus knot.
We assume that  $\CH^{\lambda}(K)$ is a quotient of $\CH^{\lambda}(T(p,\infty))$, the latter space being
$(a,q,t_{c})$-graded isomorphic to the $p\lambda$-colored homology of the unknot. Under this assumption,
the multiplicative generators of $\CH^{\lambda}(K)$ coincide with the ones of $\CH^{\lambda}(T(p,\infty))$,
while some additional relations are imposed.  Let us determine the gradings of these generators.
The diagram $p\lambda$ is a $R\times pS$ rectangle, which can be naturally divided into $p$ different $R\times S$ rectangles.
Let us denote by $u_{ij}^{(n)}$ a bosonic generator in $n$-th rectangle {\em from the right} at box $(i,j)$,
and by   $\xi_{ij}^{(n)}$ a fermionic generator in $n$-th rectangle {\em from the left} at box $(i,j)$, as in Figures \ref{nrect1} and \ref{nrect2} .
Their gradings are given by the following equations:
\begin{eqnarray*}
(a,q)\left[u_{ij}^{(n)}\right]=(0,2(nS-i+R-j+1)),\\
(t_c,t_r)\left[u_{ij}^{(n)}\right]=2(nS-i),2((n-1)R+j-1)),\\
(a,q)\left[\xi_{ij}^{(n)}\right]=(2,2((n-1)S+i-j)),\\
(t_c,t_r)\left[\xi_{ij}^{(n)}\right]=(2((n-1)S+i)-1,2((n-1)R+j)-1).
\end{eqnarray*}
Note that $t_{r}$ is constant in rows and $t_{c}$ is constant in columns. This explains the notations
$t_{r}$ and $t_{c}$. Using the equation $Q=\frac{q+t_r-t_c}{R}$, we get
$$Q\left[u_{ij}^{(n)}\right]=2n,\quad Q\left[\xi_{ij}^{(n)}\right]=2n-2.$$

It is important to note that these elements generate {\em unreduced} stable homology. To work with the reduced homology,
one should consider $u_{ij}^{(n)},\xi_{ij}^{(n)}$ for $n\ge 2$ since the subalgebra generated by  $u_{ij}^{(1)},\xi_{ij}^{(1)}$
coincides with the unreduced homology of the unknot.

There exists a map $M_{\lambda} :\CH^{\lambda}(T(p,\infty))\to \CH^{\lambda^t}(T(p,\infty))$
such that for any $w$ we have
\begin{equation}
\label{gradings of M}
(Q,t_{c},t_{r})[M_{\lambda} (w)]=(Q,t_{r},t_{c})[w].
\end{equation}
The  map $M_{\lambda}$ is defined by the equations
$$M_{\lambda} \left(u_{ij}^{(n)}\right):=u_{R-j+1,S-i+1}^{(n)},\quad M_{\lambda} \left(\xi_{ij}^{(n)}\right):=\xi_{ji}^{(n)}.$$
One can check that it satisfies (\ref{gradings of M}).
The refined exponential growth conjecture also holds for the stable homology $\CH^{\lambda}(T(p,\infty)):$
\be
\label{stable exp growth}
\CP^{R\times S}(T(p,\infty),a,Q,t_{r})=\left(\CP^{\Lambda^R}(T(p,\infty),a,Q,t_{r})\right)^{S}.
\ee
The exponential growth conjecture follows from the fact that $t_{r}$ is constant in rows:
if we consider generators in $(a,Q,t_{r})$ gradings, then we get $S$ copies of $\CH^{\Lambda^R}(T(p,\infty))$ in the same gradings, which in turn implies (\ref{stable exp growth}).

Let us recall that in (\ref{Llambda}) we conjectured that the self-duality map $\Phi$ can be obtained as follows:
there exists a ``D0-brane'' operator $L_{\lambda}:\CH^{\lambda}(K)\to\CH^{\lambda}(K)$
such that
\begin{equation}
(a,q,Q,t_r,t_c)[L_{\lambda}]=(0,R+S,2,R,S),\quad \Phi_{\lambda}(x)=L_{\lambda}^{-Q(x)}.
\end{equation}
We conjecture that for torus knots colored by $\lambda=R\times S$ rectangular diagrams
such operator $L_{\lambda}$ coincides with the multiplication operator by $u_{S,1}^{(2)}.$ One can check that it has
the prescribed degrees. Moreover, the mirror map $M_{\lambda}$ is chosen in such a way that $M_{\lambda} (u_{S,1}^{(2)})=u_{R,1}^{(2)}$,
and so for the rectangular diagrams we have an equation
$$M_{\lambda} \circ L_{\lambda}=L_{\lambda^t}\circ M_{\lambda}, $$
and therefore
$$M_{\lambda} \circ \Phi_{\lambda} = \Phi_{\lambda^t}\circ M_{\lambda} .$$
This is indeed the required compatibility relation between $M$ and $\Phi$.

Now let us describe the action of colored differentials
in this algebraic model.
The positive row-removing differentials are given by the equation
$$d^{+}_{R\times S\to k\times S}\left(\xi_{i,j}^{(n)}\right)=u_{S+1-i,j-k}^{(n)}, j>k.$$
The positive column-removing differentials are given by the equation
$$d^{+}_{R\times S\to R\times l}\left(\xi_{i,j}^{(n)}\right)=u_{S+l+1-i,j}^{(n)}, i>l.$$
The negative row-removing differentials are given by the equation
$$d^{-}_{R\times S\to k\times S}\left(\xi_{1,k+1}^{(2)}\right)=1.$$
The negative column-removing differentials are given by the equation
$$d^{-}_{R\times S\to R\times l}\left(\xi_{l+1,1}^{(2)}\right)=1.$$
The differentials vanish on all other generators of the colored HOMFLY homology.
One can check that the degrees of these differentials agree with the expected ones.

The duality between positive and negative differentials is party explained by the identities
$$d^{+}_{R\times S\to k\times S}(\xi_{1,k+1}^{(2)})=u_{S,1}^{(2)},\ d^{-}_{R\times S\to k\times S}(\xi_{1,k+1}^{(2)})=1,$$
$$d^{+}_{R\times S\to R\times l}(\xi_{l+1,1}^{(2)})=u_{S,1}^{(2)},\ d^{-}_{R\times S\to R\times l}(\xi_{l+1,1}^{(2)})=1.$$
Finally, we describe the action of the universal differentials. As a model example, consider the unreduced $S^2$ homology of the unknot. It has two even generators $u_1,u_2$ and two odd generators $\xi_1,\xi_2$.
We define the universal differential on the bottom row by the formula
$$d(f(u_1,u_2))=u_2\frac{f(u_1,u_2)-f(-u_1,u_2)}{u_1}.$$
In other words,
$$d(u_1^{a}u_2^{b})=\begin{cases}
2u_1^{a-1}u_2^{b+1},& {\rm if}\ a\ {\rm is\ odd},\\
0,& {\rm if}\ a\ {\rm is\ even}.\\
\end{cases}$$
It is clear that $d^2=0$ and the homology of $d$ is isomorphic to $\mathbb{C}[u_1^2]$, as predicted by the grading restrictions in Section \ref{sec:gradings}. Similarly, on higher levels we define the differential by the equation
$d(\xi_1)=\xi_2$, so $d(\xi_1\wedge \xi_2)=0$. The total homology is a free algebra generated by $u_1^2$ and $\xi_1\xi_2$.

In general, we can always split variables into pairs, define a universal differential in each pair as above and extend it to the whole algebra by the Leibnitz rule. For example, in stable $S^2$ homology of the $(2,\infty)$ torus knot we have two even generators $u_3,u_4$ and two odd generators $\xi_3,\xi_4$. As in the above example, one can define a universal differential by the equations
$$d(u_3^{a}u_4^{b})=\begin{cases}
2u_3^{a-1}u_4^{b+1},& {\rm if}\ a\ {\rm is\ odd},\\
0,& {\rm if}\ a\ {\rm is\ even},\\
\end{cases}
\quad d(\xi_3)=\xi_4.
$$
It clearly agrees with the description of the universal differential for the trefoil in Section \ref{sec:trefoil}.

\subsection{Differential forms}
\label{sec:symtorusforms}

\begin{definition}
The {\em unreduced} moduli space $\overline{\CM}_{p,q}(r)$ is defined in the affine space
with the coordinates $u_1, \ldots, u_{rp}; v_1,\ldots,v_{rq}$ by the equation
\begin{equation}
\label{equv}
(1+u_1z+u_2z^2+\ldots+u_{rp}z^{rp})^{q}=(1+v_1z+v_2z^2+\ldots+v_{rq}z^{rq})^{p}.
\end{equation}
which should hold at every coefficient of its expansion in powers of  $z$.
\end{definition}

The space $\overline{\CM}_{p,q}(r)$ can be also defined in the affine space
with the coordinates $u_1,\ldots,u_{rp}$ by the coefficients in the $z$-expansion of the equation
$$(1+u_1z+u_2z^2+\ldots+u_{rp}z^{rp})^{q/p},$$
starting from $(qr+1)$-st. Indeed, one can rewrite the equation (\ref{equv}) as
$$(1+u_1z+u_2z^2+\ldots+u_{pr}z^{pr})^{q/p}=(1+v_1z+v_2z^2+\ldots+v_{qr}z^{qr}).$$
This will express $v_i$ through $u_j$. The dimension of $\overline{\CM}_{p,q}(r)$ equals to $r$.
Indeed, let us forget the scheme structure and study the underlying subset of the affine space.
One has $U(z)^{q}=V(z)^{p}$ for coprime $p$ and $q$ iff there exists a polynomial $F(z)$ such that
$U(z)=F(z)^{p}, V(z)=F(z)^{q}$. Since $F(z)$ starts from 1 and has degree $r$, we have $r$ parameters at our disposal.

\begin{conjecture}
\label{symmetric torus unreduced}
The {\em unreduced} $S^r$ homology of the $(p,q)$ torus knot is the space of differential forms on $\overline{\CM}_{p,q}(r)$:
\begin{equation}
\overline{\CH}^{S^r} (T(p,q)) \; = \; \Omega^{\bullet}(\overline{\CM}_{p,q}(r)) \,.
\label{HMunreduced}
\end{equation}
\end{conjecture}

As an illustration of this conjecture, let us consider an example of the unknot. We can choose $p=q=1$, so (\ref{equv}) takes the form
$$1+u_1z+u_2z^2+\ldots+u_{r}z^r=1+v_1z+v_2z^2+\ldots+v_{r}z^r$$
Therefore $\overline{\CM}_{p,q}(r)=\Spec\C[u_1,\ldots,u_r],$
and $$\Omega^{\bullet}(\overline{\CM}_{p,q}(r))=\C[u_1,\ldots,u_r,du_1,\ldots,du_r].$$
This agrees with the above description of the unreduced $S^r$-colored triply graded homology of the unknot.

For a more interesting example, consider the $S^2$-colored trefoil knot.
The  generators are $u_2$, $u_3$, $u_4$, and the defining equations are
\be
\label{s2 colored trefoil}
-\frac{3}{16} u_3 (u_2^2 - 4 u_4)=\frac{3}{128} (u_2^4 - 8 u_2 u_3^2 - 8 u_2^2 u_4 + 16 u_4^2)=0,
\ee
$$-\frac{1}{32} u_3 (-3 u_2^3 + 2 u_3^2 + 12 u_2 u_4)=0.$$

\begin{definition}
The {\em reduced} moduli space $\CM_{p,q}(r)$ is defined in the affine space
with the coordinates $u_{r+1},\ldots,u_{pr}; v_{r+1},\ldots,v_{qr}$ by all coefficients in the $z$-expansion of the equation
\begin{equation}
\label{equvred}
(1+u_{r+1}z^{r+1}+u_{r+2}z^{r+2}+\ldots+u_{pr}z^{pr})^{q}=(1+v_{r+1}z+v_2z^2+\ldots+v_{qr}z^{qr})^{p}
\end{equation}
\end{definition}

One can check that $\CM_{p,q}(r)$ defines a single point 0, so it is zero-dimensional.
We propose the following

\begin{conjecture}
\label{symmetric torus reduced}
The {\em reduced} $S^r$ homology of the $(p,q)$ torus knot is the space of differential forms on $\CM_{p,q}(r)$:
\begin{equation}
\CH^{S^r} (T(p,q)) \; = \; \Omega^{\bullet}(\CM_{p,q}(r)) \,.
\label{HMreduced}
\end{equation}
\end{conjecture}

For $r=1$ this conjecture is equivalent to the main conjecture of \cite{GORS}, since by the results of \cite{Gdaha} the right hand side of (\ref{HMreduced}) can be identified with a certain isotypic component of the finite-dimensional representation of rational Cherednik algebra with parameter $p/q$.
Moreover, it follows from \cite{GORS} for $r=1$ and from \cite{EGL} for general $r$ that the $(a,q)$-character
of $\Omega^{*}(\overline{\CM}_{p,q}(r))$ is equal to the unreduced $S^r$-colored HOMFLY polynomial of the $(p,q)$ torus knot.
One can also check (\cite{EGL}) that the $(a,q)$-character of $\Omega^{*}(\CM_{p,q}(r))$ is equal to the reduced $S^r$-colored HOMFLY polynomial
of the same knot. 

Let us describe the gradings on the space $\Omega^{*}(\CM_{p,q}(r))$.
The $a$-grading is defined by the degree of a differential form, so that $a(u_i)=0,\ a(du_i)=2.$ Furthermore, we define the $(q,t_r,t_c)$ gradings by the formula
\be
\label{deg ui for symmetric}
(q,t_c,t_r)[u_i]=\left(2i,2i-2,2\left\lfloor\frac{i-1}{r}\right\rfloor\right),
\ee
$$(q,t_r,t_c)[du_i]=\left(2i-2,2i-1,2\left\lfloor\frac{i-1}{r}\right\rfloor+1\right).$$
It is easy to check that the defining equations of $\CM_{p,q}(r)$ are homogeneous in $q$-grading and not homogeneous in $t$-gradings. Therefore, strictly speaking, on $\Omega^{*}(\CM_{p,q}(r))$ we get  $(a,q)$ bigrading and a pair of filtrations $(t_r,t_c)$.

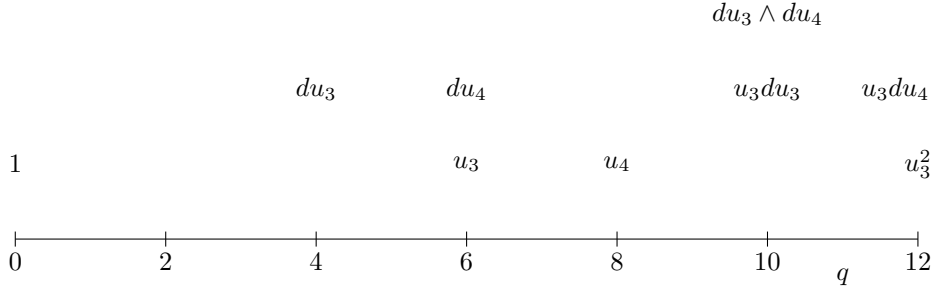
\begin{figure}
\begin{tikzpicture}
\draw (0,0) node {$1$};
\draw (6,0) node {$u_3$};
\draw (8,0) node {$u_4$};
\draw (12,0) node {$u_3^2$};
\draw (4,1) node {$du_3$};
\draw (6,1) node {$du_4$};
\draw (10,1) node {$u_3du_3$};
\draw (11.7,1) node {$u_3du_4$};
\draw (10,2) node {$du_3\wedge du_4$};
\draw (0,-1)--(12,-1);
\draw (0,-1.1)--(0,-0.9);
\draw (2,-1.1)--(2,-0.9);
\draw (4,-1.1)--(4,-0.9);
\draw (6,-1.1)--(6,-0.9);
\draw (8,-1.1)--(8,-0.9);
\draw (10,-1.1)--(10,-0.9);
\draw (12,-1.1)--(12,-0.9);
\draw (0,-1.3) node {0};
\draw (2,-1.3) node {2};
\draw (4,-1.3) node {4};
\draw (6,-1.3) node {6};
\draw (8,-1.3) node {8};
\draw (10,-1.3) node {10};
\draw (12,-1.3) node {12};
\draw (11,-1.5) node {$q$};
\end{tikzpicture}
\caption{Algebraic model for the $S^2$ homology of the trefoil.}
\label{s2tref}
\end{figure}

For example, in the uncolored case the reduced homology of the trefoil knot has only one bosonic generator $u_2$,
and the defining equation for $\CM_{2,3} (1)$ reads as
\be
u_2^2=0
\label{tref3876}
\ee
Similarly, for the $(3,4)$ torus knot the reduced uncolored homology has two bosonic generators $u_2$ and $u_3$,
with the following defining equations for $\CM_{3,4} (1)$:
\be
u_2u_3=0,\ \frac{2 u_2^3}{9} = u_3^2.
\label{knot34829054}
\ee
Finally, the reduced $S^2$-colored homology of the trefoil knot
has two even generators $u_3$ and $u_4$ and the defining equations can be obtained from (\ref{s2 colored trefoil}) by setting $u_2=0$:
$$u_3u_4=u_4^2=u_3^3=0.$$
Their differentials have the form
$$u_3du_4+u_4du_3=2u_4du_4=3u_3^2du_3=0.$$
Therefore one can check that the monomial basis in $\Omega^{*}(\CM_{2,3} (2))$ is given by
\be
1,~ u_3,~ u_4,~ u_3^2,~ du_3,~ du_4,~ u_3du_3,~ u_3du_4,~ du_3\wedge du_4.
\label{tref456}
\ee
We illustrate this homology with the Figure \ref{s2tref} (compare with the Figures in Section \ref{sec:trefoil}).

For example, by computing the gradings of the basic elements \eqref{tref456} one finds
the Poincar\'e polynomial of the $S^2$-colored trefoil knot (in the $(a,q,t_r)$ grading):
\be
\CP^{S^2}(3_1)=1+q^6t^2+q^8t^2+q^{12}t^4+a^2q^4t^3+a^2q^6t^3+a^2q^{10}t^5+a^2q^{12}t^5+a^4q^{10}t^6.
\ee
More generally, for the $S^r$-colored reduced homology of the trefoil knot
the $(a,q,t_r)$ gradings of the generators are given by the equations:
$$(a,q,t_r)[u_i]=(0,2i,2),\  (a,q,t_r)[du_i]=(2,2i-2,3).$$
In particular, the monomial $u_{r+1}^{i}\cdot du_{r+1}\wedge du_{r+2}\wedge du_{r+j}$
has grading
$$i(0,2r+2,2)+\sum_{l=1}^{j}(2,2r+2l-2,3)=(2j,2(r+1)i+2jr+j(j-1),2i+3j)$$
and contributes the term
$a^{2j} q^{2(r+1)i+2rj+j(j-1)}t_r^{2i+3j}$ in the Poincar\'e polynomial
\begin{equation}\label{f1}
\CP^{S^r}(3_1)= a^{2r}q^{-r} \sum_{i+j\le r} { a^{2j} q^{2(r+1)i+2rj+j(j-1)}t_r^{2i+3j} \frac{ [r]!}{[i]![j]![r-i-j]!}, }
\end{equation}
obtained from the categorification of the Volume Conjecture and the super-$A$-polynomial \cite{FGS1,FGS2,FGSS}.
Here $[N]!$ is the unbalanced quantum factorial:
\[
[N]!=[N]\cdot [N-1]\cdot \ldots[2]\cdot[1]',
\]
\[
[N]'=\frac{q^N-1}{q-1}=1+q+\ldots +q^{N-1}.
\]
and one can easily verify that (\ref{f1}) satisfies the refined exponential growth conjecture
as well as the rest of the structural properties of section \ref{sec:gradings}.

Finally, let us match the gradings \eqref{deg ui for symmetric} with those in the construction of Section \ref{sec:alggrad}.
We have $R=1, S=r$, so there are  generators $u_{i,1}^{(n)},\xi_{i,1}^{(n)},\ 1\le i\le r$ with gradings
$$(a,q,t_c,t_r)[u_{i,1}^{(n)}]=(0,2+2rn-2i,2rn-2i,2n-2),\ $$
$$(a,q,t_c,t_r)[\xi_{i,1}^{(n)}]=(2,2r(n-1)+2i-2,2r(n-1)+2i-1,2n-1).$$
Indeed, these gradings agree with (\ref{deg ui for symmetric}) if one identifies
$$u_{i,1}^{(n)}=u_{rn+1-i},\ \xi_{i,1}^{(n)}=du_{r(n-1)+i}.$$



\subsection{Differentials for supergroups}

It has been conjectured in \cite[Conj. 3.22]{G} that the differentials $d_{0},d_{\pm n}$ belong to a bigger algebra.
In particular, for the uncolored homology of $(4,5)$ torus knot the following operators were introduced:
$\alpha_2$ of $(a,q,t)$ degree $(-2,-2,-5)$, $\alpha_1$ of $(a,q,t)$ degree $(-2,2,-3)$.
Their introduction was motivated by the observation that the ``usual'' differentials $d_N$ applied to the (unique)
generator with top $a$-degree 6 do not generate uncolored triply graded homology: there are 7 nontrivial differentials
$d_{0},d_{\pm n}$, while there are 9 generators with $a$-degree 4. The extra differentials solve this problem.
We can now identify these operators with the supergroup differentials:
$$\alpha_2=d_{1|2}, \quad \alpha_1=d_{2|1}.$$
In total, for the uncolored (4,5) knot we have 9 nontrivial differentials:
$$d_1=d_{1|0},\ \quad d_{-1}=d_{0|1},$$
$$d_2=d_{2|0},\ \quad d_0=d_{1|1},\ \quad d_{-2}=d_{0|2},$$
$$d_{3}=d_{3|0},\ \quad \alpha_1 =d_{2|1},\ \quad \alpha_2=d_{1|2},\ \quad d_{-3}=d_{0|3}.$$
For a general $(n,m)$ torus knot $(m>n)$ we will have all differentials $d_{a|b}$ with $a+b<n$,  of total number
$$2+3+4+...+n=\frac{n(n+1)}{2}-1=\frac{(n+2)(n-1)}{2}.$$
It has been explained in \cite{GORS} that the differentials are tightly related to the
(conjectural) action of the rational Cherednik algebra $H_{c}$ on the triply graded homology of torus knots.
More precisely, \cite[Definitions 7.2-7.3]{GORS} assign a differential to every $S_n$-equivariant copy
of the standard $(n-1)$-dimensional irreducible representation $V_{n}$ inside $H_c$. These differentials satisfy some nice
properties: for example, if the two copies $\alpha$ and $\beta$ are {\em pure}, then by \cite[Lemma 7.7]{GORS} the corresponding differentials
$d_{\alpha}$ and $d_{\beta}$ anticommute.

The new supergroup differentials can be naturally embedded in this framework.
Namely, recall that $H_{c}$ has generators $x_i, y_j$.
We define $\alpha(m|k)$ to be a copy of  $V_{n}$ spanned by $x_{i}^{m}y_{i}^{k}$,
and
$$d_{m|k}=d_{\alpha(m|k)}.$$
Since this copy is pure in the sense of \cite{GORS}, the corresponding differentials anticommute.
Moreover, one can check that $\alpha(m|k)$ span $\Hom_{S_{n}}(V_{n},H_{c})$, so the supergroup differentials
form a complete collection of differentials for the uncolored torus knot homology. The grading conventions of \cite{GORS}
assure that $d_{m|n}$ have the prescribed gradings.

\subsection{Examples}

We have compared the results of \ref{sec:trefoil} with the algebraic model in (\ref{tref456}). Let us make the similar comparisons for other torus knots in Section \ref{sec:examples}.

The stable reduced $\tableau{2 2}$-colored homology
of $(2,\infty)$ torus knot have 4 even and 4 odd generators
$u_{ij}^{(2)}$ and $\xi_{ij}^{(2)}$ (see Section \ref{sec:torus}). They have the following $(a,q,t_c,t_r)$-degrees:
$$\deg[u_{11}^{(2)}]=(0,10,6,4),\ \deg[u_{12}^{(2)}]=(0,8,6,6),$$
$$\deg[u_{21}^{(2)}]=(0,8,4,4),\ \deg[u_{22}^{(2)}]=(0,6,4,6),$$
$$\deg[\xi_{11}^{(2)}]=(2,4,5,5),\ \deg[\xi_{12}^{(2)}]=(2,2,5,7),$$
$$\deg[\xi_{21}^{(2)}]=(2,6,7,5),\ \deg[\xi_{22}^{(2)}]=(2,4,7,7).$$
As explained in Section \ref{sec:torus},  $Q\left[u_{ij}^{(2)}\right]=4$ and $Q\left[\xi_{ij}^{(2)}\right]=2$.
For example, the ``bottom row'' of the homology (with appropriate shift of gradings) is 16-dimensional and has monomial basis
$$1,~u_{21}^{(2)},~\left(u_{21}^{(2)}\right)^2,~\left(u_{21}^{(2)}\right)^3,~\left(u_{21}^{(2)}\right)^4,$$
$$u_{11}^{(2)},~u_{11}^{(2)}u_{21}^{(2)},~u_{11}^{(2)}~\left(u_{21}^{(2)}\right)^2,$$
$$u_{12}^{(2)},~u_{12}^{(2)}u_{21}^{(2)},~u_{12}^{(2)}~\left(u_{21}^{(2)}\right)^2,$$
$$u_{22}^{(2)},~u_{22}^{(2)}u_{21}^{(2)},~u_{22}^{(2)}~\left(u_{21}^{(2)}\right)^2,$$
$$\left(u_{22}^{(2)}\right)^2,\left(u_{11}^{(2)}\right)^2.$$
The self-symmetry map $\Phi$ reflects first four sets of monomials in the vertical axis, and preserved two remaining monomials.

The stable reduced $\tableau{2 2 2}$-colored homology
of $(2,\infty)$ torus knot have 6 even and 6 odd generators. They have the following $(a,q,t_c,t_r)$-degrees:
$$\deg[u_{11}^{(2)}]=(0,12,6,6),\ \deg[u_{12}^{(2)}]=(0,10,6,8),$$
$$\deg[u_{13}^{(2)}]=(0,8,6,10),\ \deg[u_{21}^{(2)}]=(0,10,4,6),$$
$$\deg[u_{22}^{(2)}]=(0,8,4,8),\ \deg[u_{23}^{(2)}]=(0,6,4,10),$$
$$\deg[\xi_{11}^{(2)}]=(2,4,5,7),\ \deg[\xi_{12}^{(2)}]=(2,2,5,9),$$
$$\deg[\xi_{13}^{(2)}]=(2,0,5,11),\ \deg[\xi_{21}^{(2)}]=(2,6,7,7).$$
$$\deg[\xi_{22}^{(2)}]=(2,4,7,9),\ \deg[\xi_{23}^{(2)}]=(2,2,7,11).$$
Indeed,  $Q\left[u_{ij}^{(2)}\right]=4$ and $Q\left[\xi_{ij}^{(2)}\right]=2$.

Finally, consider the $S^2$-colored $(3,4)$ torus knot.
In the notations of Section \ref{sec:symtorusforms}, we have
4 generators $u_3,u_4,u_5,u_6$ of $(a,q,t_c,t_r)$ gradings
$$\deg u_3=(0,6,4,2),\ \deg u_4=(0,8,6,2),$$
$$\deg u_5=(0,10,8,4),\ \deg u_6=(0,12,10,4).$$
Note that in the notations of Section \ref{sec:torus},
one would have
$$u_3=u_{21}^{(2)},\ u_4=u_{11}^{(2)},\ u_5=u_{21}^{(3)},\ u_6=u_{11}^{(3)},$$
and the grading conventions agree.
The defining relations have the form
$$u_3^3 - 9 u_4 u_5 - 9 u_3 u_6=2 u_3^2 u_4 - 3 (u_5^2 + 2 u_4 u_6)=u_3 u_4^2 + u_3^2 u_5 - 3 u_5 u_6=0,$$ $$5 u_3^4 - 12 u_4^3 - 72 u_3 u_4 u_5 - 36 u_3^2 u_6 + 54 u_6^2=5 u_3^3 u_4 - 9 u_4^2 u_5 - 9 u_3 (u_5^2 + 2 u_4 u_6) =0.$$
One can check that the monomial basis on the ``bottom row'' is given by the following 25 monomials:
$$1,\ u_3,\ u_3^2,\ u_3^3,\ u_3^4,\ u_3^5,\ u_3^6,$$
$$u_4,\ u_4u_3,\ u_4u_3^2,\ u_4u_3^3,\ u_4u_3^4,$$
$$u_5,\ u_5u_3,\ u_5u_3^2,\ u_5u_3^3,$$
$$u_6,\ u_6u_3,\ u_6u_3^2,\ u_6u_3^3,$$
$$u_4^2\, u_4^2u_3,\ u_4^2u_3^2,$$
$$u_5^2, u_4^3.$$
The self-symmetry map $\Phi$ reflects first five sets of monomials in the vertical axis, and preserved two remaining monomials.


\section{HOMFLY homology from A-model and B-model}
\label{sec:mirror}

In this section, we wish to look for geometric and physical models for colored HOMFLY homology,
eventually formulating it in the language of symplectic geometry or the so-called A-model.

Indeed, in mirror symmetry, the so-called B-model is famous for its formulation in the language of
complex (algebraic) geometry, while the A-model side involves symplectic invariants, Fukaya category, and quantum cohomology \cite{KontsevichHMS}.
Clearly, the algebraic model described in the previous section is closer to the B-model, and therefore
it will be our natural starting point here.
Then, we will relate it to various problems --- that are interesting in their own right --- in mirror symmetry and in symplectic geometry.
We will give several reformulations, in terms of A-model and B-model, open and closed ({\it i.e.} with and without D-branes),
in all of which colored HOMFLY homology will be realized as cohomology of a suitable BRST operator $Q$, such that $Q^2 = 0$:
\be
\CH^{\lambda} (K) \quad = \quad  Q\text{-cohomology} \,.
\ee
In physics literature, $Q$ is often called ``supercharge'' and
its cohomology is often called the {\it space of BPS states} since its elements
are the so-called BPS states (= supersymmetric configurations).
Even though precisely this interpretation of knot homology was proposed in \cite{GSV} and studied from various
vantage points in \cite{AShakirov,DGH,GukovRTN,IK,fiveknots}, the A-model and B-model reformulations discussed here
appear to be completely new and do not make direct contact with any of the existent geometric / physical models of knot homologies.
Therefore, by pursuing some of these new geometric descriptions further one is likely to learn a lot about colored HOMFLY homology and, possibly, even about homological mirror symmetry.


\subsection{B-model on a supermanifold}
\label{sec:Bmodel}

With this goal in mind, we start by reformulating the Conjectures \ref{symmetric torus unreduced} and \ref{symmetric torus reduced}
in terms of B-model on a super-manifold $\C^{d|d}$ with a certain Landau-Ginzburg potential $W_{super} (K; \lambda)$, such that
\be
\CH^{\lambda} (K) \; \cong \; Jac (W_{super} (K;\lambda))
\label{HasJac}
\ee
Note, that $\C^{d|d}$ is a super-Calabi-Yau manifold (see \cite{AVsuper,Schwarz,Wtwistor}).

Specifically, consider the potential $W (T(p,q); S^r)$ on the $u$-space given by the formula
\be
W (T(p,q); S^r) =\Coef_{(p+q)r+1}(1+u_1z+\ldots+u_{pr}z^{pr})^{\frac{q+p}{p}}.
\label{WTpqr}
\ee
Let us rename $du_i$ by $\xi_i$ and introduce a potential $W_{super} (T(p,q); S^r)$ on the superspace $\C^{pr|pr}$
with coordinates $u_i,\xi_i$ by the formula
\be
W_{super} (T(p,q); S^r) :=\sum_{i}\frac{\partial W}{\partial u_{i}}\xi_{i} \,.
\label{WsuperTpqr}
\ee
The space of differential forms $\Omega^{\bullet}\left(\overline{\CM}_{p,q}(r)\right)$ is isomorphic to the Jacobi ring of $W_{super} (T(p,q); S^r)$:
\begin{equation}
\label{W super}
\Omega\left[\overline{\CM}_{p,q}(r)\right]=\C[u_1,\ldots,u_{pr},\xi_1,\ldots,\xi_{pr}]/\left(\frac{\partial W_{super}}{\partial u_{i}},\frac{\partial W_{super}}{\partial \xi_{i}}\right).
\end{equation}
According to the Conjecture \ref{symmetric torus unreduced}, both spaces represent the unreduced $S^r$-colored homology of the $(p,q)$ torus knot. Indeed, the partial derivatives of $W (T(p,q); S^r)$ are equal to
$$\frac{\partial W (T(p,q); S^r)}{\partial u_{i}}
=\Coef_{(p+q)r+1}\frac{\partial}{\partial u_{i}}(1+u_1z+\ldots+u_{pr}z^{pr})^{\frac{q+p}{p}}=$$
$$\frac{q+p}{p}\Coef_{(p+q)r}\left[z^{i}(1+u_1z+\ldots+u_{pr}z^{pr})^{\frac{q}{p}}\right]=$$
$$\frac{q+p}{p}\Coef_{(p+q)r+1-i}(1+u_1z+\ldots+u_{pr}z^{pr})^{\frac{q}{p}},\ i=1\ldots pr.$$
Therefore the ring of functions on $\overline{\CM}_{p,q}(r)$ coincides with the Jacobi ring of $W$:
\begin{equation}
\C\left[\overline{\CM}_{p,q}(r)\right]=\C[u_1,\ldots,u_{pr}]/\left(\frac{\partial W}{\partial u_{i}}\right).
\end{equation}

The full space of differential forms $\Omega^{\bullet}\left(\overline{\CM}_{p,q}(r)\right)$ is a quotient of the space 
of all differential forms on $\C^{pr}$ by the equations
\begin{equation}
\label{forms eq}
\frac{\partial W}{\partial u_{i}}=0,\quad \sum_{j}\frac{\partial^2 W}{\partial u_{i}\partial u_j}\cdot du_{j}=0.
\end{equation}
The Jacobi ring of $W_{super}$ is defined by the following equations:
$$\frac{\partial W_{super}}{\partial \xi_{i}}=\frac{\partial W}{\partial u_{i}}=0$$
$$\frac{\partial W_{super}}{\partial u_{i}}=\sum_{j}\frac{\partial^2 W}{\partial u_{j}\partial u_{i}}\xi_{j}=0.$$
Since these equations coincide with (\ref{forms eq}), we get the relation (\ref{W super}).

In the mirror symmetry literature \cite{Schwarz,Sethi,AVsuper,Katzarkov},
the setup we just described is called a Landau-Ginzburg B-model with target space $\C^{d|d}$
and superpotential $W_{super} (K;\lambda)$. The function $W_{super} (K;\lambda)$ is called
{\it superpotential} even when the target space is an ordinary (bosonic) manifold,
{\it e.g.} $W(u_i)$ in our discussion might be called superpotential of a Landau-Ginzburg model on $\C^{pr}$.
As a final clarification of terminology, we should mention that a model is called {\it Landau-Ginzburg model}
when the superpotential is non-zero, and is referred to as the {\it sigma-model} (on a certain target manifold) when $W=0$.

In our case, the colored HOMFLY homology (of torus knots) is realized as the Jacobi ring\footnote{also known as
the ``chiral ring'' or ``ring of observables,'' see below} \eqref{HasJac} - \eqref{W super}
of a Landau-Ginzburg model on $\C^{pr|pr}$ with the superpotential $W_{super} (T(p,q); S^r)$.
There are several aspects of this B-model that make it rather special and interesting,
the main of which is undoubtedly the fact that the target space $\C^{pr|pr}$ is a supermanifold.

It is instructive, though, to start with a more familiar and conventional B-model with a target space $X$.
Then, we will extend it to a Landau-Ginzburg model with a superpotential $W$ and will see
why \eqref{HasJac} - \eqref{W super} can be interpreted as a ring of observables with and without D-branes
(called B-branes in the context of the B-model).
In general, we shall use $X$ for the target space of the B-model and $Y$ for the target space of the $A$-model,
so that $X$ and $Y$ are mirror manifolds in most of our applications.

First, as a warm-up and a lightning review of the relevant facts about B-model, let us consider
a B-model with target space $X$.
The space of observables --- or, to be more precise, the space of {\it closed (string) observables} ---
in such model is then identified with the $\bar \p$-cohomology of $X$ with values in $\bigwedge^p T_{X}$:
\be
\label{hom of polyvectors}
H^q \left( X, \bigwedge^p T^{1,0}_{X} \right)
\ee
Indeed, every $(0,q)$-form on $X$ with values in $\bigwedge^p T_{X}$
\be
A = d \bar z^{\bar k_1} \ldots d \bar z^{\bar k_q}
A^{j_1 \ldots j_p}_{\bar k_1 \ldots \bar k_q} \frac{\p}{\p z_{j_1}} \ldots \frac{\p}{\p z_{j_p}}
\ee
in the topological B-model defines an observable:
\be
\CO_A = \eta^{\bar k_1} \ldots \eta^{\bar k_q} A^{j_1 \ldots j_p}_{\bar k_1 \ldots \bar k_q} \theta_{j_1} \ldots \theta_{j_p}
\ee
Contracting this ``$(-p,q)$-form'' with the holomorphic $(n,0)$-form $\Omega$ on $X$ we obtain an isomorphism
between this space and $(n-p,q)$-forms, where $n$ is the total (complex) dimension of $X$.
Moreover, $\{ Q, \CO_A \} = - \CO_{\bar \p A}$, so that observables in the B-model are Dolbeault cohomology
of forms valued in exterior powers of the holomorphic tangent bundle, {\it i.e.} elements of \eqref{hom of polyvectors}.

Incorporating a superpotential $W \ne 0$ leads to a Landau-Ginzburg model where
the $Q$-operator acts as a differential $Q = \bar \p + Q_{\text{bdry}}$, such that
\be
Q^2 = W \cdot \text{id}
\ee
As a result, the space of closed string states in the Landau-Ginzburg model can be described
as a Hochschild homology or, equivalently, as a hypercohomology:
\be
HH_* (MF (X, W)) \cong \mathbb{H}^* (\Lambda^* \Omega_{X}, dW \wedge)
\ee
where $\Lambda^* \Omega_{X}$ denotes the exterior powers of the sheaf of differential forms.
The dual cohomology theory is
\be
HH^* (MF (X, W)) \cong \mathbb{H}^* (\Lambda^* T_{X}, [W, - ])
\ee
where $\Lambda^* T_{X}$ denotes exterior powers of the tangent sheaf and $[W, - ]$ denotes the Lie bracket with $W$.
In particular, if $f_1, \ldots, f_n$ are sections of line bundles $\CE_i$ on a smooth variety $V$,
then they define a function $W: \text{Tot} (\oplus \CE_i^*) \to \C$ and we have the following equivalence
of categories\footnote{The analogous statement for the topological A-model is:
$$
Fuk (Y,W) \cong Fuk (Sing (W)) \,.
$$}
\be
MF (\text{Tot} (\oplus \CE_i^*), W) \cong QCoh ( \{ f_i = 0 \} )
\ee
that allows to describe B-branes in a Landau-Ginzburg theory with the superpotential $W$
({\it i.e.} objects in the category of matrix factorizations of $W$)
in terms of the category of coherent sheaves on the critical variety of $W$.

\begin{theorem}[\cite{Orlov}]
\be
H^0 (MF (W)) \cong \CD^b_{sing} (W^{-1} (0))
\ee
\end{theorem}

Note, this brings us very close to our applications \eqref{HMunreduced} and \eqref{HMreduced},
where the moduli spaces $\overline{\CM}_{p,q}(r)$ and $\CM_{p,q}(r)$
are defined by polynomial equations $f_i (u_j) = 0$ in the space $X = \C^{pr}$ of $u$'s.
In situations like this, when the equations $f_i (u_j) = 0$ integrate to a superpotential $W (u_j)$
the chiral ring (also known as the $(c,c)$ ring) of the Landau-Ginzburg model on $X = \C^n$
is isomorphic to the quotient of the ring of polynomial functions $\C[u_1, ... , u_n]$ by $dW$,
\be
\CR_{\text{closed}} (\C^n, W) \; = \; \C [u_i] / (dW)
\label{Rclosed}
\ee
When $W$ is quasihomogeneous, this ring is graded by the R-charge symmetry $U(1)_R$.

In order to interpret the colored HOMFLY homology, though,
in the Landau-Ginzburg model we need not only to find a natural home for the moduli spaces
$\overline{\CM}_{p,q}(r)$ and $\CM_{p,q}(r)$
but also to describe the corresponding spaces of differential forms \eqref{HMunreduced} and \eqref{HMreduced} in the B-model language.
One way to do this is to replace the target space $\C^n$ by the supermanifold $X = \C^{n|n}$
and to replace the superpotential $W(u_i)$ by $W_{super} (u_i, \xi_i)$ defined in \eqref{WsuperTpqr}.
This gives then
\be
\CR_{\text{closed}} (\C^{n|n}, W_{super}) \; = \; \Omega^* (\CM)
\ee
where $\CM = \{ u \in \C \vert f_i (u) \equiv \p_i W (u) = 0 \}$.
Specifically, for the superpotential $W_{super}$ defined by \eqref{WTpqr} and \eqref{WsuperTpqr}
we get a B-model realization of the $S^r$-colored HOMFLY homology of torus knots.
For instance, from \eqref{s2 colored trefoil}, \eqref{tref3876}, and \eqref{knot34829054}
we find the following potentials in simple examples considered earlier:
$$W (3_1; \tableau1)=u_2^3,$$
$$W_{super} (3_1; \tableau1)=3u_2^2\xi_2,$$
$$W (8_{19}; \tableau1)=-\frac{7 u_2^4}{243} + \frac{14 u_2 u_3^2}{27},$$
$$W_{super} (8_{19}; \tableau1) = \left(-\frac{28 u_2^3}{243} + \frac{14 u_3^2}{27}\right) \xi_2+\frac{28 u_2 u_3}{27} \xi_3,$$
$$W (3_1; \tableau2)=\frac{5}{256}u_3 (3 u_2^4 - 8 u_2 u_3^2 - 24 u_2^2 u_4 + 48 u_4^2).$$

\vskip 0.25cm

Although conceptually this Landau-Ginzburg model should be the effective two-dimensional $\CN=2$ theory
on $D$ in the brane construction \eqref{theoryB} studied in \cite{DGH,FGS2,FGSS},
the setup appears to be surprisingly different.
Yet, in both cases, the Poincar\'e polynomial of the colored HOMFLY homology is realized as
a supersymmetric index (character) of the two-dimensional $\CN=2$ theory on $D$.
One important difference is that here we have a different 2d $\CN=2$ theory for each particular value of the color,
whereas in loc. cit. a single 2d $\CN=2$ theory determines $S^r$-colored superpolynomials for all values of $r$.
It would be interesting to understand better how the new Landau-Ginzburg realization of the colored HOMFLY homology
discussed here relates to the 2d $\CN=2$ theory on $D$ studied in \cite{DGH,FGS2,FGSS}.

\subsection{B-branes and open B-model}

There is another way, however, to describe the space of differential forms that is familiar to practitioners of mirror symmetry.
This will require introducing one more ingredient in our story, namely D-branes, and considering the space of open string states,
{\it i.e.} the states of open strings stretched between branes.
In the context of topological B-model branes are usually called B-branes.
Mathematically, B-branes are objects of the derived category of coherent sheaves $\CD^b (X)$ in the case of the sigma-model (when $W=0$),
or the category of matrix factorizations $MF (X, W)$ in the case of Landau-Ginzburg model with the target space $X$ and superpotential $W$.
We already saw both of these categories in our previous discussion and now we will spend more time discussing their objects, or B-branes.

Note that both the derived category of coherent sheaves and the category of matrix factorizations
are also familiar to practitioners of knot homology: these categories play a key role
in the constructions of {\it e.g.} \cite{CKamnitzer,KRa,KRb,Webster,Wu,Yonezawa}.

Now, let us describe the ``open'' analogue of \eqref{hom of polyvectors} and \eqref{Rclosed},
{\it i.e.} the space of states of open strings ending on a brane $\CB$ in the B-model,
starting with the sigma-model on $X$ (with $W=0$), as we did in our previous discussion,
and then extend it to branes in more general Landau-Ginzburg models.
In both cases, the conclusion will be that open strings ending on $\CB$ form an algebra, which is the Ext algebra\footnote{The algebra structure comes from the Yoneda product on the self-Ext groups.}
\be
\CR_{\text{open}} (\CB) \; = \; \Ext^* (\CB,\CB) \,.
\label{Ropen}
\ee
In order to see how it comes about, let us first consider two different branes, $\CB_1$ and $\CB_2$,
represented by shaves $\CE$ and $\CF$ supported on $S \subset X$.
Then, states of open strings stretched between $\CB_1$ and $\CB_2$ are elements of
\be
\label{homs of tensor prod}
H^q \left( S , \CE^{\vee} \otimes \CF \otimes \bigwedge^p \CN_{S/ X} \right)
\ee
There is a spectral sequence with \eqref{homs of tensor prod} as the second page that converges to \cite{Sharpe}:
\be
\Ext^{p+q}_{X} (\CE, \CF)
\ee
and in many cases one has
\be
\Ext^{n}_{X} (\CE, \CF) \; \cong \;
\bigoplus_{p+q=n} H^q \left( S , \CE^{\vee} \otimes \CF \otimes \bigwedge^p \CN_{S/ X} \right)
\ee
because the spectral sequence degenerates at $E_2$.


Now let us see how the boundary chiral ring \eqref{Ropen} can help us to produce
yet another interpretation of \eqref{HMunreduced} and \eqref{HMreduced} in {\it open} B-model.
The simplest example is the unknot colored by $\lambda = S^r$, whose HOMFLY homology
$\overline{\CH}^{S^r} (\unknot)$ is realized as the space of differential forms on
$V = \Spec\C[u_1,\ldots,u_r]$.
The corresponding choice of B-brane is closely related to SYZ brane \cite{SYZ} that plays an important role in mirror symmetry.
Namely, it is well known that for a ``zero-brane'' ({\it i.e.} a skyscraper sheaf $\CO_p$)
supported at a smooth point on a manifold $X$ of complex dimension $r$ the open string algebra \eqref{Ropen}
is the exterior algebra
\be
\Ext^* (\CO_p , \CO_p) \; \cong \; \Lambda^* (V)
\ee
where $V = T_p X  \cong  \Lambda^* \C^r$.
In order to upgrade $\Lambda^* (V)$ to the space of differential forms
$\Omega^* (V) = \C [V] \otimes \Lambda (V^{\vee})$,
one simply needs to replace $X \cong \C^r$ by its complexification $X_{\C} \cong X \times X$,
and consider a B-brane $\CB = \CO_{\Delta}$, where $\Delta \cong \C^r$
is the diagonal in $X \times X \cong \C^r \times \C^r$.
Note, the new target space $X \times X \cong \C^{2r}$
can be viewed as a SYZ torus fibration with a singular fiber at the origin.
To facilitate a discussion of the mirror A-model that will follow next,
it is convenient to replace $X = \C^r$ by $X = (\C^*)^r$.
Then, the target space $X \times X \cong (\C^*)^{2r}$ is also a $T^{2r}$
fibration, but without singular fibers and one has
\be
\CR_{\text{open}} (\CB) \; = \; \Ext^* (\CO_\Delta,\CO_\Delta) \; \cong \; \C [u_i^{\pm 1} , du_i]
\label{openDD}
\ee
Before we proceed to the mirror A-model, let us briefly discuss another interesting feature of this example.


Note that \eqref{openDD} can be identified with $HH^* (\CD^b (X)) \cong \C [u_i^{\pm 1} , du_i]$
which, according to our previous discussion, describes {\it closed} string states, not {\it open}.
This is an illustration of a more general phenomenon:
open string states \eqref{Ropen} of a B-brane supported on the diagonal $\Delta \subset X \times X$
can be identified with the closed string states \eqref{hom of polyvectors} of a B-model on $X$.
The mathematical content of this statement is summarized in the following:

\begin{theorem}[Hochschild-Kostant-Rosenberg isomorphism \cite{HKR}]
\be
\Ext^{n}_{X \times X} (\CO_{\Delta}, \CO_{\Delta}) \; \cong \;
\bigoplus_{p+q=n} H^q \left( X , \bigwedge^p T_{X} \right)
\ee
where $\Delta$ is the diagonal in $X \times X$.
\end{theorem}

Moreover, Yekutieli \cite{Yekutieli} shows that there is an isomorphism from the hypercohomology to self-Ext:
\be
\mathbb{H}^n (X, D_{poly} (X)) = \Ext^n_{X \times X} (\CO_X, \CO_X)
\ee
which is compatible with the cup product on the left and the Yoneda product on the right.
And using the isomorphism between the Hochschild cohomology and the polyvector field cohomology \cite{GSchack}, one has
\be
HH^* (\CD^b (X)) \; \cong \; \Ext^*_{X \times X} (\CO_\Delta,\CO_\Delta) \; \cong \; H^* (X, \Lambda^* T_X )
\label{HHExtDD}
\ee


\subsection{HOMFLY homology from symplectic geometry}
\label{sec:Amodel}

Now, let us see what kind of A-models one gets by applying mirror symmetry to the open and closed B-models in the above discussion.

In general, we denote by $Y$ the mirror of $X$.
The simplest example is $X = (\C^*)^r$ viewed as a (trivial) $T^r$ fibration over $\R^r$.
Then, mirror symmetry is simply a T-duality (a ``Fourier transform'') along the fibers,
and the mirror manifold is $Y \cong T^* T^r$.
The mirror of the statement $HH^* (\CD^b (X)) \cong \C [u_i^{\pm 1} , du_i]$
is the statement about the Hochschild cohomology of the Fukaya category $\CF (Y)$,
\be
HH^* \left( \CF (T^* T^r) \right) \; \cong \; \C [u_i^{\pm 1} , du_i]
\label{HHTTex}
\ee
which, by analogy with \eqref{HHExtDD}, can be also realized as the boundary OPE algebra
in the A-model of $Y \times \bar Y$, where $\bar Y$ denotes the space $Y$ with the symplectic form $- \omega$.
Indeed, in general, when a symplectic manifold $Y$ has ``enough'' Lagrangians,
the Hochschild cohomology its Fukaya category $\CF (Y)$ is expected \cite{Abouzaid,KontsevichHMS,Costello}
to be isomorphic to the Lagrangian Floer cohomology of the (anti-)diagonal $\Delta_Y \hookrightarrow Y \times \bar Y$,
\be
HH^* \left( \CF (Y) \right) \; \cong \; HF^* (\Delta_Y , \Delta_Y)
\label{HHHFdiag}
\ee
Thus, in our basic example of \eqref{openDD} and \eqref{HHTTex} mirror symmetry maps a B-brane $\CB = \CO_{\Delta_X}$
to a Lagrangian A-brane supported on the conormal bundle $T_{\Delta}^{\vee}$ to $T_{\Delta} \subset T^n \times T^n$
in $Y \times Y = T^* T^r \times T^* T^r$, so that
\be
HF^* (T_{\Delta}^{\vee},T_{\Delta}^{\vee}) \; \cong \; \C [u_i^{\pm 1} , du_i]
\ee
This is an open A-model version of \eqref{HHTTex} and a concrete illustration of \eqref{HHHFdiag}.

In general, this leads us to the following picture, where the (colored) HOMFLY homology is realized
either in terms of closed A-model of a symplectic manifold $Y$
\be
\CH^{S^r} (K) \; = \; HH^* \left( \CF (Y) \right)
\ee
or in terms of open A-model (Floer homology) of a certain Lagrangian submanifold $L$,
\be
\CH^{S^r} (K) \; = \; HF^* (L,L)
\label{HKLL}
\ee
(In our previous discussion, $L$ was the (anti-)diagonal $\Delta_Y \hookrightarrow Y \times \bar Y$,
but we allow more general Lagrangian submanifolds and more general symplectic manifolds here.)
For torus knots, the closed A-model version can be obtained by applying
the supermanifold version of mirror symmetry \cite{AVsuper,Katzarkov,Sethi}
to the B-model in section \ref{sec:Bmodel}.
We plan to return to a systematic study of this mirror symmetry elsewhere.

Here, let us simply point out that even for non-torus knots one might hope to find a realization of the form \eqref{HKLL}.
For example, for the figure-eight knot $K=4_1$ and $r=1$ a natural guess is
\be
L = {\bf S}^2 \text{~with 3 punctures}
\ee
embedded in some symplectic manifold $Y$, such that $HF^* (L,L) \cong H^* (L)$.



\section{``Bottom row'' of the colored HOMFLY homology}
\label{sec:bottom}

In this part we wish to study the ``bottom row'' of the colored superpolynomial, defined as
\be
\CP^{\lambda}_{\text{bottom}} (q,t) \; : = \; \lim_{\a \to 0} \, \a^{\#} \, \CP^{\lambda} (\a,q,t) \,,
\label{bottomrow}
\ee
where $\a^{\#}$ denotes the appropriate power of $a$, such that the product $\a^{\#} \, \CP^{\lambda} (\a,q,t)$
contains only non-negative $a$-degrees starting from zero.

\subsection{Combinatorial interpretation for torus knots}
\label{sec:combbot}

First, let us recall the conjectural combinatorial description of the uncolored triply graded homology
for a $(p,q)$ torus knot \cite{G,GM,ORS}. Consider a $p\times q$ rectangle and a $NW-SE$ diagonal in it.

\begin{proposition}(\cite{G})
\label{prop:bottom1}
The dimension of the ``bottom row'' of the reduced uncolored homology equals to the number of lattice paths
in this rectangle strictly below the diagonal:
$$\dim \CH_{\text{bottom}} (T(p,q)) \; = \; c_{p,q}=\frac{(p+q-1)!}{p!q!}.$$
\end{proposition}

For example, if $q=p+1$, this dimension equals to the $p$-th Catalan number.

\begin{proposition}(\cite{G})
The dimension of the ``$k$-th row'' of the reduced uncolored homology equals to the number of lattice paths
below the diagonal with $k$ marked corners:
$$\dim \CH_{a=k} (T(p,q)) \; = \; \frac{(p+q-k-1)!}{ p\cdot q\cdot k!\cdot (p-k-1)!(q-k-1)!}.$$
\end{proposition}

A $q$-analogue of this formula was discussed in \cite[eq. (3.46)]{BEM}.
Both of these propositions motivate the following combinatorial construction,
completed in \cite{ORS} where it was matched to the algebro-geometric constructions of \cite{os}.

The basis in $\CH_{p,q}$ is enumerated by the lattice paths in the $p\times q$ rectangle with marked corners. To every such path we assign three gradings, where $a$-grading corresponds to the number of marked corners. Let us describe $q-$ and $t-$ gradings.
For simplicity, we will focus on the bottom row, considered in details in \cite{GM}. Corrections for the higher rows can be found in \cite{ORS}.

Let us interpret a lattice path as a Young diagram $D$. Consider the $\C^{*}\times \C^{*}$ action on the Hilbert scheme of
$N$ points on $\C^2$, and let us restrict it to a one-parameter subgroup $(t^p,t^q)$. One can check (see {\it e.g.} \cite{Nakajima})
that the action of this subgroup has isolated fixed points parametrized by the Young diagrams with $N$ boxes. Therefore the corresponding moment map is a Morse function on $\Hilb^{N}(\C^2)$, and its gradient flow induces the Bialynicki-Birula cell decomposition of the Hilbert scheme \cite{Nakajima,Elstro}. We are interested in the dimension of the cell $C_{D}$ labelled by the fixed point $D$, which is given by the following formula:
$$\dim C_{D}=|D|+h^{+}_{p/q}(D),$$
where
\begin{equation}
\label{hplus}
h^{+}_{p/q}(D)=\sharp\left\{c\in D~~ \vline~~ \frac{a(c)}{l(c)+1}\le \frac{p}{q}< \frac{a(c)+1}{l(c)}\right\}.
\end{equation}
Here $a(c)$ and $l(c)$ denote the arm and the leg of a cell $c$ in a diagram $D$:

\begin{tikzpicture}
\draw  (0,0)--(0,4)--(5,0)--(0,0);
\draw  (1.5,0.5)--(1.5,0.7)--(1.7,0.7)--(1.7,0.5)--(1.5,0.5);
\draw [thick] (0,3)--(0.5,3) -- (0.5,2) -- (2.3,2)--(2.3,1.3)--(3,1.3)--(3,0.5)--(4,0.5)--(4,0)--(0,0)--(0,3);
\draw [<->,>=stealth] (1.7,0.6) -- (3,0.6);
\draw [<->,>=stealth] (1.6,0.7) -- (1.6,2);
\draw (2.3,0.4) node {${a(c)}$};
\draw (1.3,1.3) node {${l(c)}$};
\draw (2.5,-0.2) node {$p$};
\draw (-0.2,2) node {$q$};
\draw [->,>=stealth] (2.6,-0.15) -- (5,-0.15);
\draw [->,>=stealth] (2.4,-0.15) -- (0,-0.15);
\draw [->,>=stealth] (-0.15,1.8) -- (-0.15,0);
\draw [->,>=stealth] (-0.15,2.2) -- (-0.15,4);
\end{tikzpicture}

The statistics $h^{+}_{p/q}(D)$ was introduced in \cite{LW} in connection with some combinatorial models of the bigraded Catalan numbers.

\begin{proposition}
\label{prop:bottom3}
The gradings of Dyck paths are given by the following equations:
\begin{enumerate}

\item{The $q$-grading assigned to $D$ equals to $2\left[|D|+2h^{+}_{p/q}(D)\right]$.}

\item{The $t$-grading assigned to $D$ equals to $2|D|$.}
\end{enumerate}
\end{proposition}

Motivated the refined exponential growth conjecture, we assume that the dimension of the $S^r$-colored homology is equal to the $r$-th power of the uncolored homology, {\it i.e.} the $r$-th power of the generalized Catalan number $c_{p,q}$.
Therefore it is natural to assume that one can find a basis in the $S^r$-colored homology labelled by all possible $r$-tuples of Young diagrams in $p\times q$ rectangle below the diagonal.

\begin{conjecture}
\label{prop:bottom4}
The bottom row of the $S^r$ colored homology of the $(p,q)$ torus knot has a basis labeled by the $r$-tuples $(D_1,\ldots,D_r)$ of the Dyck paths in $p\times q$ rectangle. The $(Q,t_r)$ gradings of Dyck paths are given by the following equations:
\begin{enumerate}

\item{The $Q$-grading assigned to $D$ equals to $2\sum_{i=1}^{r}\left[|D_i|+2h^{+}_{p/q}(D_i)\right]$.}

\item{The $t_r$-grading assigned to $D$ equals to $2\sum_{i=1}^{r}|D_i|$.}
\end{enumerate}
\end{conjecture}

The $q$-grading (and, therefore, the $t_c$-grading) has yet to be constructed. We expect it to be related to the combinatorial constructions of \cite{BDZ} and \cite{HMZ}.


\subsection{Coupled instanton-vortex counting}
\label{sec:inst-surf}

The ``bottom row'' \eqref{bottomrow} has a simple and beautiful interpretation in terms of instanton / vortex counting.
As a warm-up, let us start with a simple vortex counting problem that is easy to do ``by hand''
and that has all the essential features. Then, we shall upgrade it to the equivariant instanton counting relevant to arbitrary knots and links.

Let $\CV_m$ be the moduli space of $m$ abelian vortices on a two-dimensional plane, $D = \R^2$.
In other words, $\CV_m$ is the moduli space of solutions to the PDEs
\be
\begin{array}{rl}
* F_A & =~ i |\phi|^2 - it \\
\bar \partial_A \phi & =~ 0
\end{array}
\label{abelianvorteq}
\ee
where $A$ is a connection on a $U(1)$ bundle of first Chern class $m$, and $t$ is a parameter.
These equations describe supersymmetric configurations (BPS states) in a two-dimensional theory
on $D$ with $\CN=(2,2)$ supersymmetry:
\be
\CT_{\text{vortex}}~:~ U(1)~\text{gauge theory with a single charged field}~\phi
\ee
It is well known that $\CV_m$ is a K\"ahler manifold of (real) dimension $2m$. Namely,
\be
\CV_m \; = \; \Sym^m (\C) \; \equiv \; \C^m / S_m \,,
\label{vortsymm}
\ee
where one can think of coordinates on $\CV_m$ as vortex positions on $D = \R^2 \cong \C$.

Since $D$ admits a circle action, one can consider the equivariant character, $Ch_q (\CV_m)$,
with respect to the rotation group $U(1)_q$ acting on $D$.
Indeed, $U(1)_q$ acts on $\CV_m$ in a natural way, by equal phase rotations on all factors in the symmetric product \eqref{vortsymm}.
Moreover, identifying the vortex moduli space $\CV_m = \Sym^m (\C)$
with the space of monic polynomials of degree $m$,
\be
f (x) \; = \; \prod_{j=1}^m (x - x_j) = x^m + a_1 x^{m-1} + \ldots + a_m \,,
\ee
we immediately deduce the isomorphism $\CV_m \cong \C^m = (\{ a_1, \ldots, a_m \})$
and also learn that $U(1)_q$ acts on the space $\CV_m$ with weights $(1, 2, \ldots, m)$.
Therefore, the $U(1)_q$-equivariant character of $\CV_m$ is
\be
Ch_q (\CV_m) \; = \; \frac{1}{(1-q)(1-q^2) \ldots (1-q^m)} \,,
\label{unknotvortchar}
\ee
in which a careful reader will recognize the ``bottom row'' of the \textsl{unnormalized} $S^m$-colored HOMFLY polynomial of the unknot.

This is not an accident, of course, and the physics setup \eqref{theoryB} predicts that the ``bottom row''
of any $\lambda$-colored HOMFLY homology of any knot $K$ has a similar interpretation in terms
of equivariant vortex counting on $D = \R^2$.
In order to explain this relation, we need to generalize our warm-up example in two important ways:
first, we need to introduce the homological $t$-grading, and secondly we need a generalization to arbitrary knots.
Both of these problems can be achieved by embedding $D = \R^2$ into a larger space $M_4 = \R^4$,
and realizing our vortex counting problem on $D$ as a special case of counting solutions
to coupled \textsl{instanton-vortex equations} or, equivalently, via equivariant instanton
counting on the 4-manifold $M_4$ in the presence of ramification along $D \subset M_4$.

However, before we proceed to generalizations, let us stay a bit longer in the world of vortex equations
and explore the picture which seems to emerge:
\be
\text{knot}~K \quad \leadsto \quad \text{2d vortex theory}~\CT_{\text{vortex}} (K) \quad \leadsto \quad \CV_m (K)
\ee
so that a suitable version $Ch_{q,t} (\CV_m)$ of the equivariant character computes the bottom row of
the colored superpolynomial \eqref{bottomrow} with $\lambda = S^m$,
\be
Ch_{q,t} (\CV_m) \; = \; \CP^{S^m}_{\text{bottom}} (q,t) \,.
\label{naiveVmP}
\ee
Already at this stage the reader might suspect that incorporating the $t$-grading and passing to knot homologies
can be addressed by considering equivariant cohomology or K-theory of the vortex moduli space,
\be
H^*_{U(1)_q} (\CV_m (K)) \,,
\label{naiveHvort}
\ee
where $t$-grading is identified with the homological grading.
Although this idea is a little too naive, it is actually on the right track, and it is instructive to pursue it a little further.

Therefore, as the natural next step let us consider a non-abelian generalization of \eqref{abelianvorteq}
that describes supersymmetric solutions in the 2d theory:
\be
\CT_{\text{vortex}}~:~ U(p+1)~\text{gauge theory with}~p+1~\text{fundamental fields}~\phi_i
\ee
The corresponding vortex equations
\be
\begin{array}{rl}
* F_A & =~ i\sum_{i=1}^{p+1} \phi_i \phi_i^{\dagger} - i t \\
\bar \partial_A \phi_i & =~ 0
\end{array}
\label{vortexeq}
\ee
involve a $U(p+1)$ gauge connection and $p+1$ Higgs fields $\phi_i$ in the fundamental
representation of the gauge group. Again, the $m$-vortex moduli space is a K\"ahler manifold
of real dimension
\be
\dim~ \CV_m^{U(p+1)} = 2 m (p+1) \,,
\label{mmndimension}
\ee
In particular, the single-vortex moduli space is well known to be a $(2p+2)$-dimensional space
\be
\CV_1^{U(p+1)} \; \cong \; \C \times \cp^p
\ee
where the factor $\C$ parametrized by the ``center-of-mass'' position,
while $\cp^p$ encodes the ``internal degrees of freedom'' of a single non-abelian vortex.
The rotation symmetry $U(1)_q$ acts on $\C$ with weight 1,
so we expect that it contributes to the equivariant character a factor $\frac{1}{1-q}$.
However, unlike our first example (related to the unknot), now the moduli space $\CV_1^{U(p+1)}$
has non-trivial topology and, therefore, if $t$-grading is simply the homological grading on \eqref{naiveHvort}
the suitable generalization of the equivariant character should be a product of $\frac{1}{1-q}$
and the Poincar\'e polynomial of $\cp^p$,
\be
\frac{1}{1-q} \sum_{i=0}^p t^{2i}
\ee
This turns out to be the correct answer for the bottom row of the HOMFLY homology of the $(2,2p+1)$ torus knot.

Now let us explain why all these examples are ``working'' and where the connection with vortex equations come from.
The reason, in fact, is already contained in section~\ref{sec:geomphys}
where we reviewed the interpretation~\eqref{HHBPS} of homological knot invariants
in terms of supersymmetric configurations (BPS states) in the string theory setup~\eqref{theoryB}.
What this interpretation tells us is that specializing to the bottom row \eqref{bottomrow},
{\it i.e.} taking the limit $\a \to 0$, means the ``large volume limit.''
Indeed, according to the identification of the parameters \eqref{NvolCP1},
in this limit the interesting geometry (and topology) of the Calabi-Yau space $X$
is replaced by the simplest Calabi-Yau 3-fold, namely a flat space:
\be
X \;\overset{\a \to 0}{\leadsto}\; \C^3 \,.
\label{Xatozero}
\ee
In other words, the bottom row of the colored HOMFLY homology is described by
a much simpler ``toy model'' of \eqref{theoryB} in which the pair $(X,L_K)$ is replaced by $(\C^3, L_K)$.
(Not much happens to the Lagrangian submanifold $L_K$ in this limit.)

The next step is to look at this system from the vantage point of the 4-manifold $M_4 \cong \R^4$
and the defect (called the ``surface operator'' \cite{Ramified}) supported on $D \cong \R^2$.
Before we took the limit $\a \to 0$, the compactification on $X$ produces abelian gauge theory
on $M_4$ with gauge group $U(1)$ \cite{KKV}.
Incorporating an extra brane in this setup means including the so-called surface operator supported on $D$
or, mathematically speaking, ramification in the gauge theory on $M_4$.

Let $\CM_{k,m}$ be the moduli space of abelian instantons on $M_4 = \R^4$ with ramification along $D \subset M_4$.
In other words, $\CM_{k,m}$ is the moduli space of solutions to the self-duality equation on $M_4 \setminus D$,
\be
F_A^+ \; = \; 0
\label{sdeq}
\ee
with the second Chern class $k$, with the monopole number $m = \frac{1}{2\pi} \int_D F_A$,
and with the prescribed behavior along $D$.
The latter is what we call the ramification data;
{\it e.g.} a simple example of the so-called \textsl{tame} ramification
can be obtained by introducing a $\delta$-function source
in the self-duality equation \eqref{sdeq} on $M_4$:
\be
F_A^+ \; = \; 2\pi \alpha (\delta_D)^+ \,.
\ee
The gauge connections which solve this equation have a first-order pole at $D$.
More generally, one can study solutions to self-duality equations with singularities of arbitrarily high order,
\be
A = dz  \left( \frac{\alpha_1}{z} + \frac{\alpha_2}{z^2} + \ldots \right) + d \bar z \left( \ldots \right) \,.
\ee
This is called \textsl{wild} ramification and is precisely how the dependence on the knot $K$ will enter,
via the choice of the ramification data.\footnote{For instance, it was conjectured~\cite{DGH} that $(2,2p+1)$ torus knots
correspond to wild ramification of order $p$.}
Let us denote the corresponding moduli space by $\CM_{k,m} (K)$.

Much like the vortex moduli space $\CV_m$ enjoyed an action of the rotation group $U(1)_q$,
the moduli space of ``ramified instantons'' on $M_4 \setminus D \cong \C \times \C^*$
has an action of the rotation group $U(1)_q \times U(1)_t$,
where $U(1)_q$ acts by rotations of $D = \C$ and $U(1)_t$ acts by rotations of its normal bundle, $\C^*$.

Then, in the limit \eqref{Xatozero} the physical interpretation \eqref{HHBPS} of the colored HOMFLY homology predicts
\be
Ch_{q,t} (\CM_{0,m} (K)) \; = \; \bar \CP^{S^m}_{\text{bottom}} (K; q,t)
\ee
which is very close to \eqref{naiveVmP} and, in fact, is exactly the sought-after proper generalization of it.
The connection with vortex equations is now easy to see. At least for some knots
the moduli space $\CM_{k,m} (K)$ can be equivalently described as the moduli space of solutions
to the coupled \textsl{instanton-vortex equations}:
\begin{eqnarray}
F_A^+ & = & 2\pi \phi \otimes \phi^{\dagger} \; (\delta_D)^+  \nonumber \\
*_{D} F_B & = & i \phi^{\dagger} \otimes \phi - i t {\bf 1}_E \label{coupledvort} \\
\bar \partial_{A,B} \phi & = & 0 \nonumber
\end{eqnarray}
where $A$ is a unitary connection on the line bundle $\CL$ over $M_4$,
$B$ is a unitary connection on the bundle $E$ over $D \subset M_4$,
$\phi$ is an element $\phi \in H^0 (\Hom (E,\CL \vert_D ))$, and $\phi^{\dagger}$ is its adjoint.

As far as we know the instanton-vortex equations \eqref{coupledvort} are new and have not appeared in the literature previously.
Their closest cousin is a set of the so-called coupled vortex equations studied {\it e.g.} in~\cite{BGP1,BGP2},
where both $E$ and $\CL$ are bundles over $D$. In other words, our equations \eqref{coupledvort} can be viewed
as analogs of the coupled vortex equations where one of the bundles is extended over $M_4$.
For instance, for a $(2,2p+1)$ torus knot $T^{2,2p+1}$ one takes $E$ to be a unitary bundle of rank $p+1$.

Then, as before, to get the ``bottom row'' of the \textsl{unreduced} $S^m$-colored HOMFLY homology we need to consider
the equivariant character of $\CM_{0,m} (K)$ with instanton number $k=0$ and vortex number $m$.
In situations where $\CM_{k,m} (K)$ can be identified with the moduli space of solutions to \eqref{coupledvort}
this boils down to studying non-abelian vortex equations on $D$.
Thus, for $2p+1$ torus knots we end up with $U(p+1)$ vortex equations,
so that the equivariant character of $\CM_{0,m} (T^{2,2p+1}) \cong \CV_m^{U(p+1)}$
essentially consists of $p+1$ copies of the abelian vortex character that we discussed earlier, {\it cf.} \eqref{unknotvortchar}.

Our prediction is that the equivariant character of the $m$-vortex moduli space $\CV_m^{U(p+1)}$
is the ``bottom row'' of the $S^m$-colored HOMFLY homology of the $(2,2p+1)$ torus knot, {\it i.e.}
\begin{eqnarray}\nonumber
\!\!\!\!\!\!\!\!Ch_{q,t} (\CV_m^{U(p+1)})  &=& \frac{q^{-pm}}{(q;q)_m}  \sum_{0\le k_p \le \ldots \le k_2 \le k_1 \le m}
\left[\!\begin{array}{c} m \\k_1 \end{array}\!\right]\left[\!\begin{array}{c} k_1\\k_2 \end{array}
\!\right]\cdots\left[\!\begin{array}{c} k_{p-1}\\k_p \end{array}\!\right]  \times\\
& & \nonumber \!\!\!\!\!\!\!\!\!\!\!\!\!\!\!\!\!\!\!\! \times \,\,\, q^{(2m+1)(k_1+k_2+\ldots+k_p)-\sum_{i=1}^p k_{i-1}k_i}
t^{2(k_1+k_2+\ldots+k_p)} ,
\end{eqnarray}
with the convention $k_0:=m$.
For example, when $m=2$ and $p=1$, {\it i.e.} for the $(2,3)$ torus knot (= trefoil knot $3_1$) colored by $\lambda = \tableau{2}$ we have
\be
Ch_{q,t} (\CV_2^{U(2)}) \; = \; q^{-2} \frac{1 + q^3 t^2 + q^4 t^2 + q^6 t^4}{(1-q)(1-q^2)}
\label{S2chtref}
\ee
which can be verified using equivariant $U(2)$ action on $\CV_2^{U(2)}$, see {\it e.g.} \cite{DGH} and \cite{Taipale}.

It would be interesting to understand better the relation between vortex moduli spaces
and the combinatorics of the Bialynicki-Birula cell decomposition discussed in Section \ref{sec:combbot}
that also determines the bottom row of the colored HOMFLY homology.
Thus, in the above example \eqref{S2chtref} of the $S^2$-colored trefoil, the corresponding
moduli space of two $U(2)$ vortices $\CV_2^{U(2)}$ has been studied in the literature \cite{Tong,Shifman,Baptista}
and is known to have two strata, which correspond to separated vortices and coincident vortices.
The first stratum --- of complex dimension 4 --- is simply the product of two single-vortex moduli spaces $\CV_1^{U(2)}$,
modulo the permutation of the two vortices,
\be
\CU_1 ( \CV_2^{U(2)} ) \; \cong \; \C \times \frac{\C \times \cp^1 \times \cp^1}{\Z_2}
\ee
with the ``diagonal'' removed.
The second stratum --- of complex dimension 3 ---  is known to be the space of Hecke modifications
(times the center of mass position),
\be
\CU_2 ( \CV_2^{U(2)} ) \; \cong \;  \C \times \mathbb{W} \cp^2_{(1,1,2)} \,.
\ee
See \cite{GS} for color graphics illustrating the bottom row of the $S^2$-colored HOMFLY homology
and the role of $\mathbb{W} \cp^2_{(1,1,2)}$ in its construction.
In this, and more general examples, it would be interesting to match the geometry of these strata
with the Bialynicki-Birula cell decomposition discussed in Section \ref{sec:combbot}.

Another direction for future work, already mentioned in Section \ref{sec:Bmodel},
is the relation between theory $\CT_{\text{vortex}} (K)$ and the surface operator theory considered in \cite{FGS2,FGSS}.
Indeed, here, in Section \ref{sec:Bmodel}, and in \cite{FGS2,FGSS} the Poincar\'e polynomial
of the colored HOMFLY homology (a.k.a. the superpolynomial) is realized as a certain index (or, character)
of a two-dimensional $\CN=2$ theory that describes the physics on $D$ in the brane construction \eqref{theoryB}.
Yet, the key difference is that here and in Section \ref{sec:Bmodel} there is a different $\CN=2$ theory
for each choice of ``color'' $\lambda$; the index of such theories determines colored superpolynomials
and {\it a priori} all such theories are unrelated.
In \cite{FGS2,FGSS}, on the other hand, a single $\CN=2$ theory $\CT_K$ (that depends only on the knot $K$)
determines the entire tower of $S^m$-colored homological invariants, via specializations of its supersymmetric index:
\be
\bar \CP^{S^m}_K (a,q,t) \; = \; \CI_{\CT_K} (x = q^m; a,q,t) \,.
\ee
Moreover, recursion relations found in \cite{FGS2,Nawata,FGSS}, which have a natural meaning\footnote{{\it e.g.} the characteristic
variety of the corresponding $q$-difference equations has a nice interpretation as the moduli space of SUSY parameters, {\it etc.}}
in the $\CN=2$ theory $\CT_K$, relate equivariant characters $Ch_{q,t} (\CM_{0,m} (K))$
with different values of the ``vortex number'' $m = \frac{1}{2 \pi} \int_D F$.
Similarly, from the vantage point of Section \ref{sec:Bmodel}, the recursion relations of \cite{FGS2,Nawata,FGSS}
relate B-models with target manifolds of different (super-)dimension
in a way reminiscent of the notions of ``endoscopy'' and ``transfer'' in the geometric Langlands correspondence \cite{FW}.
Thus, in the case of the trefoil knot we have:
\begin{align}
\CP^{S^{m+2}}_{\text{bottom}}
- (q^{-1} - t^2 q^m + (1+q) q^{2m} t^2) \CP^{S^{m+1}}_{\text{bottom}}
& \\
+ t^4 q^{3m-1} (q^m - 1) \CP^{S^{m}}_{\text{bottom}} & = 0 \notag
\end{align}


\begin{appendix}
\section{Rectangle-colored invariants of $(2,m)$ torus knots}

As we use the quantum invariants of $(2,m)$ torus knots colored
by the representations labeled by the rectangular Young diagrams in many examples throughout the paper, we would like
to provide a reader with a handy formula for them.

Let us first recall the general formula for the colored invariants of torus knots \cite{LZ,RJ,stevan}. Define the coefficients $c_{\lambda,n}^{\mu}$ by the equation
\begin{equation}
s_{\lambda}(x_1^n,x_2^n,\ldots)=\sum_{|\mu|=n|\lambda|}c_{\lambda,n}^{\mu}s_{\mu}.
\end{equation}
The $\lambda$-colored invariant of the $(n,m)$ torus knot can be expressed via the colored invariants of the unknot by the following formula (we omit an overall scaling factor):
\begin{equation}
\label{plethysm}
P^{\lambda}(T(n,m))=\sum_{|\mu|=n|\lambda|}q^{-\frac{m}{n}\kappa(\mu)}c_{\lambda,n}^{\mu}P^{\mu}(\unknot),
\end{equation}
where $\kappa(\mu)=\sum_{(i,j)\in \mu}(i-j)=\frac{1}{2}\sum_{j}\mu_j(\mu_j-2j+1)$ is the content of $\mu$.

Therefore to compute $P^{\lambda}(T(m,n))$, one just has to determine the
``plethysm'' coefficients $c_{\lambda,n}^{\mu}$.There is no closed formula for them  for general $\lambda,\mu$ and $n$. However, for $n=2$ and $\lambda=(S^{R})$, a simple formula for the coefficients $c_{\lambda,n}^{\mu}$ was found in \cite{CRemmel}.

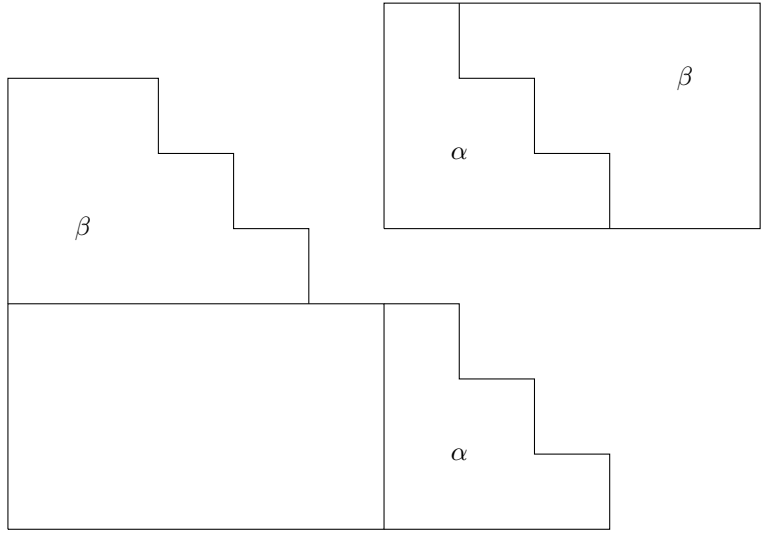
\begin{figure}[ht]
\begin{tikzpicture}
\draw (0,0)--(0,3)--(5,3)--(5,0)--(0,0);
\draw (0,3)--(0,6)--(2,6)--(2,5)--(3,5)--(3,4)--(4,4)--(4,3);
\draw (5,3)--(6,3)--(6,2)--(7,2)--(7,1)--(8,1)--(8,0)--(5,0);
\draw (6,1) node {$\alpha$};
\draw (1,4) node {$\beta$};
\draw (5,4)--(5,7)--(10,7)--(10,4)--(5,4);
\draw (6,7)--(6,6)--(7,6)--(7,5)--(8,5)--(8,4);
\draw (6,5) node {$\alpha$};
\draw (9,6) node {$\beta$};
\end{tikzpicture}
\caption{Diagram $\mu$ is balanced if diagrams $\alpha$ and $\beta$ fill the rectangle}
\label{CRemmelfig}
\end{figure}

A Young diagram $\mu=(\mu_1,\ldots,\mu_{2R})$ is called $(S,R)$--balanced, if
$$\mu_i+\mu_{2R+1-i}=2S\ \mbox{for all}\ 1\le i\le R$$
Balanced diagrams are in 1-to-1 correspondence with splittings of the $S\times R$ rectangle into pairs of complimentary diagrams, as shown
on Figure \ref{CRemmelfig}.

For an $(S,R)$--balanced diagram $\mu$, let
${\rm sgn}(\mu)=(-1)^{|\beta|}=(-1)^{\sum_{i=1}^{R}\mu_i}$. Let $B(S,R)$ denote the set of  $(S,R)$--balanced diagrams. By \cite[Theorem 1]{CRemmel}, one has
$$s_{(S^R)}(x_1^2,x_2^2,\ldots)=\sum_{\mu\in B(S,R)}(-1)^{{\rm sgn}(\mu)}s_{\mu},$$
hence
\begin{equation}
\label{CRemmeleq}
c_{(S^R),2}^{\mu}=\begin{cases}(-1)^{{\rm sgn}(\mu)},\ \mbox{if}\ \mu\in B(S,R),\\ 0,\ \mbox{otherwise}\end{cases}.
\end{equation}
The equations (\ref{plethysm}) and (\ref{CRemmeleq}) determine the invariant $P^{(S^R)}(T(2,m))$ completely.

\section{Unreduced colored HOMFLY homology}

\subsection{General structure}
As explained in Section \ref{sec:geomphys}, the unreduced $\lambda$-colored homology $\bar{\CH}^{\lambda}(K)$ for any rectangular diagram $\lambda$ and any knot $K$ is a free module over the unreduced $\lambda$-colored homology of the unknot $\bar{\CH}^{\lambda}(U)$. Therefore, for every knot $K$ the unreduced homology $\bar{\CH}^{\lambda}(K)$ is determined by the reduced homology ${\CH}^{\lambda}(K)$
and the unreduced homology of the unknot $\bar{\CH}^{\lambda}(U)$. In terms of the corresponding Poincare polynomials we have:
\be\label{simprel}
\bar{\CP}^{\lambda}(K)={\CP}^{\lambda}(K)\bar{\CP}^{\lambda}(U).
\ee
Hence, the unreduced colored homology  is determined by the reduced one.

In order to simplify the presentation, we focus on the symmetric representations $\lambda=S^r$. 

Some of the properties of reduced homology extend to the unreduced case: these include the existence and the behaviour of the positive colored differentials and the totally refined exponential growth property. By (\ref{simprel}) these properties will follow from the fact the unreduced homology of the unknot satisfies these properties.  

The unreduced colored homology of the unknot is described in Section \ref{sec:unknot} in terms of the bosonic and fermionic generators. Their expicit $(a,Q,t_r,t_c)$-gradings have been obtained, which gives the following Poincare polynomial of the tilde-version of the unreduced $S^r$-colored homology of the unknot:
\be
\tilde{\bar{\CP}}^{S^r}(U)(a,Q,t_r,t_c)=a^{-r}Q^{r}\prod_{i=1}^r \frac{1+a^2t_rt_c^{2i-1}}{1-Q^2t_c^{2i-2}}.
\ee
In ``standard" $(a,q,t_r,t_c)$-gradings we have (the same convention, $Q=q+t_r-t_c$, holds for the symmetric representations as in the case of the reduced homology):
\be
{\bar{\CP}}^{S^r}(U)(a,q,t_r,t_c)=a^{-r}q^{r}\prod_{i=1}^r \frac{1+a^2q^{2i-2}t_rt_c^{2i-1}}{1-q^{2i}t_c^{2i-2}}.
\ee

The degrees of the positive colored differentials are as in the reduced case:
\be
(a,q,t_r,t_c)[d^{col}_{1|k}]=(-2,2(1-k),-1,-2k-1),\quad 0\le k <r.
\ee
They also have the ``colored" property, i.e. we have:
\be\label{unredcol}
H^*(\bar{\CH}^{S^r}(K),d^{col}_{1|k})\cong \bar{\CH}^{S^k}(K),\quad 0\le k <r.
\ee
Indeed, the last follows from the following Lemma which shows that (\ref{unredcol}) holds for the unknot and with the re-gradings in the corresponding isomorphisms being exactly the same as for the reduced homology (Section \ref{coldif}):

\begin{lemma}
For all nonnegative integers $r$ and $k$, with $k<r$, we have that the following holds
\be
\tilde{\bar{\CP}}^{S^r}(U)(a,Q,t_r,t_c)-\tilde{\bar{\CP}}^{S^k}(U)(a,Qt_c^{r-k},t_r,t_c)=
(1+a^{-2}Q^2t_r^{-1}t_c^{-2k-1})X,
\ee
for some Laurent polynomial $X(a,Q,t_r,t_c)$ whose all coefficients are nonnegative  integers, i.e. $X\in \Z_+[a^{\pm 1},Q^{\pm 1},t_r^{\pm 1},t_c^{\pm 1}]$. 
\end{lemma}
\textit{Proof:} Straightforward computation.\\

The totally refined exponential growth property also holds for the unreduced homology of the unknot. Indeed:
\be
\tilde{\bar{\CP}}^{S^r}(U)(a,Q,t_r=t,t_c=1)=\prod_{i=1}^r \frac{1+a^2t}{1-Q^2}=\left(\tilde{\bar{\CP}}^{\tableau1}(U)(a,Q,t_r=t,t_c=1)\right)^r.
\ee

\vskip 0.25cm

On the other hand there are some differences between the reduced and unreduced homology of any knot.

The first clear distinction from the reduced case is that the unreduced homology is infinite-dimensional, since such is the unreduced homology of the unknot. In particular there can be no symmetry that inverts the $Q$, or $q$-grading and therefore the selfsymmetry is not satisfied for the unreduced homology. Consequently, there are no negative colored differentials in this case (recall that selfsymmetry exchanges positive and negative colored differentials).

The second distinction from the reduced homology case is that all differentials $d_{n|0}$, for $n>0$, are nontrivial now, even for the unknot. This is the reason why there is no such a simple relation like (\ref{simprel}) between the unreduced and reduced $sl(n)$ homology.
In the next section we compute $sl(n)$ colored HOMFLY homology for some simple knots.

\subsection{$sl(n)$ colored HOMFLY homology}

The unreduced $(sl(n),S^r)$-colored homology of a knot $K$, denoted $\bar{\CH}^{sl(n),S^r}(K)$, is obtained as the homology of $\bar{\CH}^{S^r}(K)$ with respect to $d_{n|0}$, followed by collapsing of grading by setting $a=q^n$. We observe $\bar{\CH}^{S^r}(K)$ here as triply-graded theory in $(a,q,t_c)$-gradings, and also to simplify the notation we denote $t_c$ simply by $t$.
Therefore, $\bar{\CH}^{sl(n),S^r}(K)$ becomes doubly-graded theory in $(q,t)$-degrees.

Unlike the reduced homology case, here all differentials $d_{n|0}$, with $n>0$, are highly nontrivial, even for the unknot. Below we compute $sl(2)$ homologies for the fundamental and the second-symmetric representation of the unknot, trefoil and the figure-eight knot.

The $(a,q,t)$-degrees of the $sl(n)$ differentials, $d_{n|0}$, $n>0$, are as for the reduced homology:
\be
(a,q,t)[d_{n|0}]=(-2,2n,-1),\quad n>0.
\ee

For all three knots we assume that $d_{n|0}$ is such that it cancels all pairs of generators
of $\bar{\CH}^{S^r}(K)$ whose $(a,q,t)$-degrees differ by $(-2,2n,-1)$.

\subsubsection{Unknot}
In the fundamental representation, we have that the Poincare polynomial of the unreduced uncolored HOMFLY homology is given by:
\[
{\bar{\CP}}^{\tableau1}(U)(a,q,t)=a^{-1}q \frac{1+a^2t}{1-q^{2}}.
\]
Then the homology with respect to $d_{n|0}$ is finite-dimensional  for any $n$:
\be
({\bar{\CH}}^{\tableau1}(U),d_{n|0})=a^{-1}q(1+q^2+\ldots+q^{2(n-1)}).
\ee
Finally, after collapsing the tri-grading by setting $a=q^n$, we get the familiar 
expression for the uncolored $sl(n)$ homology of the unknot:
\be
{\bar{\CP}}^{sl(n),\tableau1}(U)(q,t)=({\bar{\CH}}^{\tableau1}(U),d_{n|0})_{|_{a=q^n}}=q^{1-n}+q^{3-n}+\ldots+q^{(n-1)}.
\ee

The $S^2$-colored unreduced homology of the unknot is given by:
\be\label{s2un}
\bar{\CP}^{\tableau2}(U)=a^{-2}q^2\frac{(1+a^2t)(1+a^2q^2t^3)}{(1-q^2)(1-q^4t^2)}.
\ee
To simplify the further computation of the homology with respect to $d_{2|0}$ in terms of the Poincare polynomial, we say that $P_1$ is equal to  $P_2$ {\it modulo} $d_{2|0}$ if $P_1-P_2=(1+a^{-2}q^4t^{-1}) Y$ for some Laurent polynomial $Y$ in variables $a,q$ and $t$ whose all coefficients are nonnegative integers.

Then by using the result from the uncolored case, we have that modulo $d_{2|0}$ the expression (\ref{s2un}) equals to:
\[
a^{-2}q^2(1+q^2)\frac{1+a^2q^2t^3}{1-q^4t^2}=(1+q^2+q^4t^2\frac{1+a^2t}{1-q^4t^2})+
(1+a^{-2}q^4t^{-1})\frac{a^2q^2t^3}{1-q^4t^2}.
\]
Therefore:
\be
({\bar{\CH}}^{\tableau2}(U),d_{2|0})=a^{-2}q^2(1+q^2+q^4t^2\frac{1+a^2t}{1-q^4t^2}),
\ee
and so 
\be
{\bar{\CP}}^{sl(2),\tableau2}(U)(q,t)=({\bar{\CH}}^{\tableau2}(U),d_{2|0})_{|_{a=q^2}}=q^{-2}+1+q^2t^2\frac{1+q^4t}{1-q^4t^2}.
\ee

\subsubsection{Trefoil}
For the trefoil, we have:
\be
\bar{\CP}^{\tableau1}(3_1)={\CP}^{\tableau1}(3_1)\bar{\CP}^{\tableau1}(U)=(a^2q^{-2}+a^2q^2t^2+a^4t^3)a^{-1}q\frac{1+a^2t}{1-q^2}.
\ee
Now after taking the homology with respect to $d_{2|0}$ which cancels all pairs 
of generators which differ in $(a,q,t)$-degree by $(-2,4,-1)$, we are left with
\be
(\bar{\CH}^{\tableau1}(3_1),d_{2|0})=a(q^{-1}+q+q^3t^2)+a^3q^3t^3,
\ee
which after setting $a=q^2$ gives 
\[
\bar{\CP}^{sl(2),\tableau1}(3_1)(q,t)=(\bar{\CH}^{\tableau1}(3_1),d_{2|0})|_{a=q^2}=q+q^3+q^5t^2+q^9t^3,
\]
which equals precisely the free part of the
(unreduced) Khovanov homology of the trefoil.\\

For the second-symmetric representation, we have:
\be\label{gl31}
\bar{\CP}^{\tableau2}(3_1)={\CP}^{\tableau2}(3_1)\bar{\CP}^{\tableau2}(U),
\ee
where the reduced homology of the trefoil was computed in Section \ref{sec:examples}:
\be
{\CP}^{\tableau2}(3_1)=
a^4(q^{-4}+q^2t^4+q^4t^6+q^8t^8)+a^6(t^5+q^2t^7+q^6t^9+q^8t^{11})+a^8q^6t^{12}.
\ee
As we have seen above, the Poincare polynomial of the unknot $\bar{\CP}^{\tableau2}(U)$
modulo $d_{2|0}$ equals 
\be\label{sl2s2u}
(\bar{\CH}^{\tableau2}(U),d_{2|0})=a^{-2}(q^2+q^4+q^6t^2\frac{1+a^2t}{1-q^4t^2}).
\ee
As for the the reduced $S^2$-colored Poincare polynomial of the trefoil, ${\CP}^{\tableau2}(3_1)$, we split it according to the canceling differential $d^{col}_{0|2}$:
\be
{\CP}^{\tableau2}(3_1)=(1+a^2q^4t^5)(a^4q^{-4}+a^4q^2t^4+a^4q^4t^6+a^6q^2t^7)+a^4q^8t^8.
\ee
Multiplying the factor $(1+a^2q^4t^5)$ with the RHS of (\ref{sl2s2u}) gives:
\begin{eqnarray*}
&&(1+a^2q^4t^5)a^{-2}(q^2+q^4+q^6t^2\frac{1+a^2t}{1-q^4t^2})=\\
&&=(1+a^2q^4t^5)a^{-2}(q^2+q^4)+a^{-2}q^6t^2(1+a^2q^4t^5)\frac{1+a^2t}{1-q^4t^2})=\\
&&=(1+a^2q^4t^5)a^{-2}(q^2+q^4)+\\
&&\quad +\, a^{-2}q^6t^2(1+a^2t+q^4t^2+a^2q^4t^3+(1+a^{-2}q^4t^{-1})a^2q^4t^5\frac{1+a^2t}{1-q^4t^2})=\\
&&=a^{-2}q^2(1+q^2+q^4t^2+a^2q^4t^3+a^2q^6t^5+a^2q^8t^5)+\\
&&\quad\quad + \,(1+a^{-2}q^4t^{-1})(q^6t^5+q^{10}t^7\frac{1+a^2t}{1-q^4t^2}).
\end{eqnarray*} 
Therefore, after replacing this in (\ref{gl31}), we have that $\bar{\CP}^{\tableau2}(3_1)$ modulo $d_{2|0}$ equals:
\begin{eqnarray*}
&\!\!\!\!\!\!\!a^{-2}q^2(1+q^2+q^4t^2+a^2q^4t^3+a^2q^6t^5+a^2q^8t^5)(a^4q^{-4}+a^4q^2t^4+a^4q^4t^6+a^6q^2t^7)+\\
&\quad \quad + a^2q^8t^8 (q^2+q^4+q^6t^2\frac{1+a^2t}{1-q^4t^2}).
\end{eqnarray*}
By canceling all possible remaining pairs of generators by $d_{2|0}$ we get:
\begin{eqnarray*}
(\bar{\CH}^{\tableau2}(3_1),d_{2|0})\!\!&\!=\!&\!\! a^2(q^{-2}+1)+a^2q^2t^2+a^4q^2t^3\!+\!a^2(q^4+q^6)t^4\!+\!a^4(q^4+q^6)t^5+\\
&&\,\,+\,a^4(q^6+q^8)t^7+a^2q^{10}t^8+a^4q^{10}t^9+a^4q^{12}t^{11}+a^6q^{12}t^{12}+\\
&&\quad\quad +a^2q^6t^6(1+q^2+q^4t^2\frac{1+a^2t}{1-q^4t^2}).
\end{eqnarray*}
Finally, by setting $a=q^2$, we obtain the following $(sl(2),S^2)$-colored homology of the trefoil:
\begin{eqnarray*}
\CP^{sl(2),\tableau2}(3_1)(q,t)&=&(\bar{\CH}^{\tableau2}(3_1),d_{2|0})_{|_{a=q^2}}=\\
&&\, = (q^2+q^4)+q^6t^2+q^{10}t^3+(q^{8}+q^{10})t^4+(q^{12}+q^{14})t^5+\\
&&\quad + (q^{14}+q^{16})t^7+q^{14}t^8+q^{18}t^9+q^{20}t^{11}+q^{24}t^{12} +\\
&&\quad + q^{10}t^6 (1+q^2+q^4t^2 \frac{1+q^4t}{1-q^4t^2}).
\end{eqnarray*}
We note that this result for the free part of the homology coincides (up to an overall shift) with the one obtained in \cite{CHKrushkal}. \footnote{Note that the $so(3)$ homology corresponding to the fundamental representation that is computed in \cite{CHKrushkal} is indeed isomorphic to the homology $\CH^{sl(2),\tableau2}$, due to the well-known isomorphism between the fundamental (vector) representation of $so(3)$ and the second-symmetric representation of $sl(2)$.}\\

\subsubsection{Figure-eight knot}
Similarly, for the figure-eight knot, we have 
\be
\bar{\CP}^{\tableau1}(4_1)={\CP}^{\tableau1}(4_1)\bar{\CP}^{\tableau1}(U)=(a^2t^2+q^{-2}t^{-1}+1+q^2t+a^{-2}t^{-2})a^{-1}q\frac{1+a^2t}{1-q^2}.
\ee
After taking homology with respect to $d_{2|0}$ (i.e. modulo $(1+a^{-2}q^4t^{-1})$), we get
\be
(\bar{\CH}^{\tableau1}(4_1),d_{2|0})=a^2qt^2+q^{-1}(t^{-1}+1)+q(1+t)+a^{-2}q^{-1}t^{-2},
\ee
and after setting $a=q^2$, we are left with 
\[
(\bar{\CH}^{\tableau1}(4_1),d_{2|0})_{|_{a=q^2}}=q^5t^2+qt+q+q^{-1}+q^{-1}t^{-1}+q^{-5}t^{-2},
\]
which is the free part of the
(unreduced) Khovanov homology of the figure-eight knot.\\

For the second-symmetric representation we have: 
\be\label{gl41}
\bar{\CP}^{\tableau2}(4_1)={\CP}^{\tableau2}(4_1)\bar{\CP}^{\tableau2}(U),
\ee
where the reduced $S^2$-homology is already computed in Section \ref{sec:examples} 
\begin{eqnarray*}
{\CP}^{\tableau2}(4_1)&=&1+(1+q^2t^2)a^{-2}q^{-2}t^{-4}(1+a^{2}q^{-2}t)(1+a^{2}q^{4}t^{5})+\\
&&\quad +
a^{-4}q^{-4}t^{-8}(1+a^{2}q^{-2}t)(1+a^{2}t^{3})(1+a^{2}q^{4}t^{5})(1+a^{2}q^{6}t^{7}).
\end{eqnarray*}
In the same way as above for the trefoil, after grouping all terms with the factor $1+a^2q^4t^5$ and multiplying with the RHS of (\ref{sl2s2u}), we get that modulo $d_{2|0}$ the polynomial
$\bar{\CP}^{\tableau2}(4_1)$ equals to:
\begin{eqnarray*}
&\!\!\!\!\!\!a^{-2}q^2(1\!+\!q^2\!+\!q^4t^2\!+\!a^2q^4t^3\!+\!a^2q^6t^5\!+\!a^2q^8t^8)((1+q^2t^2)(1+a^{2}q^{-2}t)a^{-2}q^{-2}t^{-4}+\\
&\quad +a^{-4}q^{-4}t^{-8}(1+a^{2}q^{-2}t)(1+a^{2}t^{3})(1+a^{2}q^{6}t^{7}))+\\
&\quad \quad + a^{-2}(q^2+q^4+q^6t^2\frac{1+a^2t}{1-q^4t^2}).
\end{eqnarray*}
Finally, by canceling by $d_{2|0}$ all possible remaining pairs of generators we get:
\begin{eqnarray*}
(\bar{\CH}^{\tableau2}(4_1),d_{2|0})&=& 
a^{-4}(t^{-4}+(q^2+q^4)t^{-2}+q^6)+ \\
&&\quad\quad + a^{-2}(t^{-3}+2q^2t^{-1}+q^4(t^{-1}+t)+2q^6t+q^8t^3)+\\
&&\quad\quad +(q^2+(q^4+q^6)t^2+q^8t^4)+\\
&&\quad + a^{-6}q^{-2}t^{-8}+a^{-4}(q^{-2}t^{-7}+t^{-5}+q^2t^{-3}+q^4t^{-1})+\\
&&\quad\quad + a^{-2}(t^{-4}+q^2t^{-2}+2q^4+q^6t^2+q^8t^4)+\\
&&\quad\quad + (q^4t+q^6t^3+q^8t^5+q^{10}t^7)+a^2q^{10}t^8+\\
&&\quad\quad + a^{-2}q^2(1+q^2+q^4\frac{1+a^2t}{1-q^4t^2}).
\end{eqnarray*}
After setting $a=q^2$ in the above expression, we get the $(sl(2),S^2)$-colored homology of the figure-eight knot:
\begin{eqnarray*}
\CP^{sl(2),\tableau2}(4_1)(q,t)&=&(\bar{\CH}^{\tableau2}(4_1),d_{2|0})_{|_{a=q^2}}=\\
&&\quad = q^{-14}t^{-8}+q^{-10}t^{-7}+ q^{-8}t^{-5}+(q^{-8}+q^{-4})t^{-4}+\\
&&\quad\quad + (q^{-6}+q^{-4})t^{-3}+ (q^{-6}+q^{-4}+q^{-2})t^{-2} + \\
&&\quad\quad + (q^{-4}+2q^{-2}+1)t^{-1}+(q^{-2}+2+q^2)+\\
&&\quad\quad + (1+2q^2+q^4)t+ (q^2+q^4+q^6)t^2+(q^4+q^6)t^3+\\
&&\quad\quad + (q^4+q^8)t^4+q^8t^5+q^{10}t^7+q^{14}t^8+\\
&&\quad\quad\quad + q^{-2}+1+q^{2}t^2\frac{1+q^4t}{1-q^4t^2}.
\end{eqnarray*}
Again, the obtained homology matches the result from \cite{CHKrushkal}. 

\subsection{Comparison with the algebraic model}

Following Section \ref{sec:extension} and \cite{GOR}, we can compare the above computations with the corresponding algebraic models. Recall that the $S^r$-colored triply graded homology of the unknot has even generators $u_1,\ldots,u_r$ and odd generators $\xi_1,\ldots,\xi_r$ such that
$$[a,q,t_{c}]u_i=(0,2i,2i-2),\ [a,q,t_{c}]\xi_i=(2,2i-2,2i-1).$$
The differential $d_{2|0}$ is defined by the equation $d_{2|0}(\xi_i)=\sum_{i=1}^{j+1}u_ju_{i+1-j}$,
in particular,
\begin{equation}
\label{d20}
d_{2|0}(\xi_1)=u_1^2,\ d_{2|0}(\xi_2)=2u_1u_2,\ d_{2|0}(\xi_3)=u_2^2+2u_1u_3,\ d_{2|0}(\xi_4)=2u_1u_4+2u_2u_3.
\end{equation}

The uncolored HOMFLY homology of the unknot is a free algebra generated by $u_1$ and $\xi_1$, so the uncolored
$sl(2)$ homology is two-dimensional and spanned by $1$ and $u_1$. 

Using (\ref{d20}), one can check that the $(S^2, sl(2))$ homology of the unknot is spanned by
$$1,\ u_1,\ \mu_1,\ u_2^{k},\  u_2^{k}\mu_1,\qquad k\ge 1,$$
where $\mu_1=2u_2\xi_1-u_1\xi_2$. Note that the homology has nontrivial $\mathbb{Z}_2$-torsion that we do not consider here, see \cite{GOR} for more details.

The reduced triply graded homology of the trefoil is three-dimensional and spanned by $1,u_2,\xi_2$.
Therefore the uncolored $sl(2)$ homology of the trefoil can be considered as a quotient of $(S^2, sl(2))$ 
homology of the unknot by $u_2^2$ and $u_2\mu_1$, hence it is spanned by $1,u_1,u_2$ and $\mu_1$.

Finally, let us consider the $(S^2, sl(2))$ homology of the trefoil. Following (\ref{tref456}), the reduced $S^2$ HOMFLY homology of the trefoil is given by is spanned by $$1,\ u_3,\ u_4,\ u_3^2,\ \xi_3,\ \xi_4,\ u_3\xi_3, u_4\xi_3=-u_3\xi_4,\ \xi_3\xi_4.$$
For simplicity, let us focus on the bottom row. The unreduced $S^2$ HOMFLY homology of the trefoil has form
$$\langle 1,\ u_1,\ u_2^{k} \rangle\otimes\langle 1,\ u_3,\ u_4,\ u_3^2 \rangle.$$
Using (\ref{d20}), we can eliminate $u_1u_3$ and $u_1u_4$. Moreover, 
$$u_1u_3^2=d_{2|0}\left(\frac{1}{2}u_3\xi_3-\frac{1}{4}x_2\xi_4+\frac{1}{4}x_4\xi_2\right),$$
$$u_2^2u_3=d_{2|0}\left(\frac{1}{2}x_2\xi_4-\frac{1}{2}x_4\xi_2\right),$$
$$u_2^2u_4=d_{2|0}(\xi_3u_4-\xi_4u_3)+2u_2u_3^2,\ 
u_2^3=d_{2|0}(u_2\xi_3-u_3\xi_2).$$
and can be eliminated too. Therefore the bottom row of the unreduced $(S^2, sl(2))$ homology of the trefoil 
is spanned by
$$1,\ u_1,\ u_2,\ u_3,\ u_4,\ u_2^2,\ u_2u_4,\ u_2u_3,\ u_2^{k}u_3^2,\ k\ge 0.$$
The homology of $d_{2|0}$ on higher levels can be computed by similar methods.

\end{appendix}

\bibliographystyle{amsalpha}

\end{document}